\documentclass[10pt,twoside,leqno]{amsart}
\usepackage{amssymb,amsbsy,amsmath,amsfonts,amssymb,amscd,times,graphics,color,wasysym,xypic}
\newcommand{\C}                 {\mathbb{C}}

\newcommand{\K}                 {\mathbb{K}}
\newcommand{\R}                 {\mathbb{R}}
\newcommand{\N}                 {\mathbb{N}}

\newtheorem{theorem}{Theorem}
\newtheorem{lemma}{Lemma}

\newtheorem{problem}{Problem}
\newtheorem{notationalconvention}{Notational convention}

\begin{document}

%\large

\title[
Explicit prolongations of infinitesimal Lie symmetries
]{
Four explicit formulas for the prolongations
\\
of an infinitesimal Lie symmetry
\\
and multivariate Fa\`a di Bruno formulas
}

\author{Jo\"el Merker}

\address{
CNRS, Universit\'e de Provence, LATP, UMR 6632, CMI, 
39 rue Joliot-Curie, F-13453 Marseille Cedex 13, France. \ \ 
{\it Internet}:
{\tt http://www.cmi.univ-mrs.fr/$\sim$merker/index.html}}

\email{merker@cmi.univ-mrs.fr} 


\date{\number\year-\number\month-\number\day}

\begin{abstract}
In 1979, building on S.~Lie's theory of symmetries of (partial)
differential equations, P.J.~Olver formulated inductive formulas
which are appropriate for the computation of the prolongations of an
infinitesimal Lie symmetry to jet spaces, for an arbitrary
number $n\geq 1$ of independent variables $(x^1, \dots, x^n)$ and for
an arbitrary number $m \geq 1$ of dependent variables $(y^1, \dots,
y^m)$. This paper is devoted to elaborate a formalism based on
multiple Kronecker symbols which enables one to handle these
``unmanageable'' prolongations and to discover the underlying
complicated combinatorics. Proceeding progressively, we write down
closed explicit formulas in four cases: $n=m=1$; $n\geq 1$, $m=1$;
$n=1$, $m\geq 1$; general case $n\geq 1$, $m\geq 1$. As a subpart of
the obtained formulas, we recover four possible versions of the
(multivariate) Fa\`a di Bruno formula. We do not employ the classical
formalism based on the symmetric algebra ({\it cf.}~e.g. H.~Federer's
book, p.~222), because it hides several explicit sums in symbolic
compactifications, and because the presence of supplementary
complexities ({\it e.g.} splitting of indices, combinatorics of
partial derivatives) impedes us to apply such compactifications
coherently. Our method of exposition is inductive: we conduct our
reasonings by analyzing several thoroughly organized formulas, by
comparing them together and by ``drifting'' towards generality, in
homology with the classical style of L.~Euler.
\end{abstract}

\maketitle

\begin{center}
\begin{minipage}[t]{11cm}
\baselineskip =0.35cm
{\scriptsize

\centerline{\bf Table of contents}

\smallskip

{\bf 1.~Jet spaces 
and prolongations \dotfill 1.}

{\bf 2.~One independent variable and
one dependent variable \dotfill 6.}

{\bf 3.~Several independent variables and
one dependent variable \dotfill 12.}

{\bf 4.~One independent variable and
several dependent variables \dotfill 26.}

{\bf 5.~Several independent variables and
several dependent variables \dotfill 33.}

}\end{minipage}
\end{center}

\section*{\S1.~Jet spaces and prolongations}

\subsection*{1.1.~Choice of notations for the jet space variables}
Let $\K = \R$ or $\C$. Let $n\geq 1$ and $m\geq 1$ be two positive
integers and consider two sets of variables $x= (x^1, \dots, x^n) \in
\K^n$ and $y = (y^1, \dots, y^m)$. In the classical theory of Lie
symmetries of partial differential equations, one considers certain
differential systems whose (local) solutions should be mappings of the
form $y = y(x)$. We refer to~\cite{ ol1986} and to~\cite{ bk1989} for
an exposition of the fundamentals of the theory. Accordingly, the
variables $x$ are usually called {\sl independent}, whereas the
variables $y$ are called {\sl dependent}. Not to enter in subtle
regularity considerations (as in~\cite{ m2005}), we shall assume
$\mathcal{ C}^\infty$-smoothness of all functions throughout this
paper.

Let $\kappa \geq 1$ be a positive integer. For us, in a very concrete
way (without fiber bundles), the {\sl $\kappa$-th jet space}
$\mathcal{ J}_{ n, m}^\kappa$ consists of the space $\K^{n+m+ m \frac{
(n+m)!}{n! \ m!}}$ equipped with the affine coordinates
\def\theequation{1.2}\begin{equation}
\left(
x^i, y^j, y_{i_1}^j, y_{i_1,i_2}^j, 
\dots\dots, 
y_{i_1,i_2,\dots,i_\kappa}^j
\right),
\end{equation}
having the symmetries
\def\theequation{1.3}\begin{equation}
y_{i_1,i_2,\dots,i_\lambda}^j
=
y_{i_{\sigma(1)},i_{\sigma(2)},\dots,i_{\sigma(\lambda)}}^j,
\end{equation}
for every $\lambda$ with $1\leq \lambda \leq \kappa$ and for every
permutation $\sigma$ of the set $\{1, \dots, \lambda \}$. The variable
$y_{i_1, i_2, \dots, i_\lambda}^j$ is an independent coordinate
corresponding to the $\lambda$-th partial derivative $\frac{
\partial^\lambda y^j}{ \partial x^{ i_1} \partial x^{ i_2} \cdots
\partial x^{ i_\lambda }}$, which explains the symmetries~\thetag{
1.3}.

In the classical Lie theory (\cite{ ol1979},
\cite{ ol1986}, \cite{ bk1989}), all the
geometric objects: point transformations, vector fields, {\it etc.},
are local, defined in a neighborhood of some point lying in some
affine space $\K^N$. However, in this paper, the original geometric
motivations are rapidly forgotten in order to focus on combinatorial
considerations. Thus, to simplify the presentation, we shall not
introduce any special notation to speak of certain local open subsets
of $\K^{n+ m}$, or of the jet space $\mathcal{ J}_{ n, m}^\kappa =
\K^{n+m+ m \frac{ (n+m)! }{ n! \ m! }}$, {\it etc.}: we will always
work in global affine spaces $\K^N$.

\subsection*{ 1.4.~Prolongation $\varphi^{ (\kappa)}$ 
of a local diffeomorphism $\varphi$ to the $\kappa$-th jet space} In
this paragraph, we recall how the prolongation of a diffeomorphism to
the $\kappa$-th jet space is defined (\cite{ ol1979},
\cite{ ol1986}, \cite{ bk1989}).

Let $x_*\in \K^n$ be a central fixed point and let $\varphi : \K^{n+m}
\to \K^{n+m}$ be a diffeomorphism whose Jacobian matrix is close to
the identity matrix at least in a small neighborhood of $x_*$. Let
\def\theequation{1.5}\begin{equation}
J_{x_*}^\kappa
:=
\left(
x_*^i,y_{*i_1}^j,
y_{*i_1,i_2}^j,
\dots\dots,
y_{*i_1,i_2,\dots,i\kappa}^j
\right)
\in
\left.
\mathcal{J}_{n,m}^\kappa
\right\vert_{x_*}
\end{equation}
be an arbitrary $\kappa$-jet based at $x_*$. The goal is to defined
its transformation $\varphi^{(\kappa)} ( J_{x_*}^\kappa)$ by
$\varphi$.

To this aim, choose an arbitrary mapping $\K^n\ni x \mapsto g(x) \in
\K^m$ defined at least in a neighborhood of $x_*$ and representing
this $\kappa$-jet, {\it i.e.} satisfying
\def\theequation{1.6}\begin{equation}
y_{*i_1,\dots,i_\lambda}^j 
=
\frac{\partial^\lambda g^j}{
\partial x^{i_1}\cdots \partial x^{i_\lambda}}
(x_*),
\end{equation}
for every $\lambda\in \N$ with $0 \leq \lambda \leq \kappa$, for all
indices $i_1, \dots, i_\lambda$ with $1 \leq i_1, \dots, i_\lambda
\leq n$ and for every $j\in \N$ with $1\leq j\leq m$. In accordance with
the splitting $(x, y)\in \K^n\times \K^m$ of coordinates, split the
components of the diffeomorphism $\varphi$ as $\varphi = (\phi, \psi)
\in \K^n\times \K^m$. Write $\left( \overline{x}, \overline{y}
\right)$ the coordinates in the target space, so that the
diffeomorphism $\varphi$ is:
\def\theequation{1.7}\begin{equation}
\K^{n+m}
\ni
(x,y)
\longmapsto
\left(
\overline{x}, \overline{y}
\right) 
= 
(\phi(x, y), \psi (x, y))
\in
\K^{n+m}.
\end{equation}
Restrict the variables $(x, y)$ to belong to the graph of $g$, namely
put $y:= g(x)$ above, which yields
\def\theequation{1.8}\begin{equation}
\left\{
\aligned
\overline{x} 
&
=
\phi(x, g(x)), \\
\overline{y}
&
=
\psi(x,g(x)).
\endaligned\right.
\end{equation}
As the differential of $\varphi$ at $x_*$ is close to the identity,
the first family of $n$ scalar equations may be solved with respect to
$x$, by means of the implicit function theorem. Denote $x
= \overline{ \chi}(\overline{ x})$ the resulting mapping, satisfying
by definition
\def\theequation{1.9}\begin{equation}
\overline{x}
\equiv 
\phi\left(
\overline{ \chi}(\overline{ x}), 
g(\overline{ \chi}(\overline{ x}))
\right).
\end{equation}
Replace $x$ by $\overline{ \chi}(\overline{ x})$ in the second family
of $m$ scalar equations~\thetag{1.8} above, which yields:
\def\theequation{1.10}\begin{equation}
\overline{y}
=
\psi\left(
\overline{ \chi}(\overline{ x}),
g(\overline{ \chi}(\overline{ x}))
\right).
\end{equation}
Denote simply by $\overline{ y} = \overline{ g} ( \overline{ x})$ this
last relation, where $\overline{ g} ( \cdot) := \psi \left( \overline{
\chi} (\cdot), g ( \overline{ \chi} (\cdot)) \right)$.

In summary, the graph $y=g(x)$ has been transformed to the graph
$\overline{y} = \overline{ g} (\overline{ x})$ by the diffeomorphism
$\varphi$ whose first order approximation is close to the identity.

Define then the {\sl transformed jet $\varphi^{(\kappa)} \left(
J_{x_*}^\kappa \right)$} to be the $\kappa$-th jet of $\overline{ g}$ at
the point $\overline{ x}_* := \phi ( x_*)$, namely:
\def\theequation{1.11}\begin{equation}
\varphi^{(\kappa)}
\left(
J_{x_*}^\kappa
\right)
:= 
\left(
\frac{\partial^\lambda \overline{ g}^j}{
\partial \overline{ x}^{i_1} 
\cdots
\partial \overline{ x}^{i_\lambda}}
(\overline{ x}_*)
\right)_{
1\leq i_1, \dots, i_\lambda \leq n, \
0 \leq \lambda \leq \kappa
}^{
1\leq j\leq m}
\left.
\in\mathcal{J}_{n,m}^\kappa
\right\vert_{\overline{x}_*}.
\end{equation}
It may be shown that this jet does not depend on the choice of a local
graph $y = g(x)$ representing the $\kappa$-jet $J_{x_*}^\kappa$ at
$x_*$. Furthermore, if $\pi_\kappa := \mathcal{ J}_{n,m}^\kappa \to
\K^m$ denotes the canonical projection onto the first factor, the
following diagram commutes:
$$
\diagram \mathcal{J}_{n,m}^\kappa \rto^{\varphi^{(\kappa)}} 
\dto_{\pi_\kappa} 
& \mathcal{J}_{n,m}^\kappa
\dto^{\pi_\kappa} \\
\K^{n+m} \rto^{\varphi} & 
\K^{n+m}
\enddiagram.
$$

\subsection*{1.12.~Inductive formulas for the 
$\kappa$-th prolongation $\varphi^{(\kappa)}$} To present them, we
change our notations. Instead of $(\overline{ x}, \overline{ y})$, as
coordinates in the target space $\K^n\times \K^m$, we shall use
capital letters:
\def\theequation{1.13}\begin{equation}
\left(
X^1,\dots,X^n,Y^1,\dots,Y^m
\right).
\end{equation}
In the source space $\K^{n+m}$ equipped with the coordinates $(x, y)$,
we use the jet coordinates~\thetag{ 1.2} on the associated
$\kappa$-th jet space. In the target space $\K^{n+m}$ equipped with
the coordinates $(X, Y)$, we use the coordinates
\def\theequation{1.14}\begin{equation}
\left(
X^i, Y^j, Y_{X^{i_1}}^j, 
Y_{X^{i_1,i_2}}^j,
\dots\dots,
Y_{X^{i_1,i_2,\dots,i_\kappa}}^j
\right)
\end{equation}
on the associated $\kappa$-th jet space. In these notations, the
diffeomorphism $\varphi$ whose first order approximation is close to
the identity mapping in a neighborhood of $x_*$
may be written under the form:
\def\theequation{1.15}\begin{equation}
\varphi : 
\
(x^{i'},y^{j'}) 
\mapsto 
\left(
X^i,Y^j
\right) 
= 
\left(
X^i(x^{i'},y^{j'}), \ Y^j(x^{i'},y^{j'})
\right),
\end{equation}
for some $\mathcal{ C}^\infty$-smooth functions $X^i(x^{i'},y^{j'})$,
$i = 1, \dots, n$, and $Y^j(x^{i'},y^{j'})$, $j = 1, \dots, m$. The
first prolongation $\varphi^{(1)}$ of $\varphi$ may be written 
under the form:
\def\theequation{1.16}\begin{equation}
\varphi^{(1)}
: \
\left(
x^{i'},y^{j'}, y_{i_1'}^{j'}
\right) 
\longmapsto
\left(
X^i(x^{i'},y^{j'}), \
Y^j(x^{i'},y^{j'}), \
Y_{X^{i_1}}^j
\left(
x^{i'}, y^{j'}, y_{i_1'}^{j'}
\right)
\right),
\end{equation}
for some functions $Y_{X^{i_1}}^j \left( x^{i'}, y^{j'}, y_{i_1'}^{j'}
\right)$ which depend on the pure first jet variables
$y_{i_1'}^{j'}$. The way how these functions depend on the first order
partial derivatives functions $X_{x^{i'}}^i$, $X_{y^{j'}}^i$,
$Y_{x^{i'}}^j$, $Y_{y^{j'}}^j$ and on the pure first jet variables
$y_{i_1'}^{j'}$ is provided (in principle) by the following compact
formulas (\cite{ bk1989}):
\def\theequation{1.17}\begin{equation}
\left(
\begin{array}{c}
Y_{X^1}^j \\
\vdots \\
Y_{X^n}^j \\
\end{array}
\right)
=
\left(
\begin{array}{ccc}
D_1^1 X^1 & \cdots & D_1^1 X^n \\
\vdots & \cdots & \vdots \\
D_n^1 X^1 & \cdots & D_n^1 X^n \\
\end{array} 
\right)^{-1}
\left(
\begin{array}{c}
D_1^1 Y^j \\
\vdots \\
D_n^1 Y^j \\
\end{array}
\right), 
\end{equation}
where, for $i' = 1, \dots, n$, the $D_{ i'}^1$ denote the {\sl $i'$-th
first order total differentiation operators}:
\def\theequation{1.18}\begin{equation}
D_{i'}^1 
:= 
\frac{\partial}{\partial x^{i'}}
+
\sum_{j'=1}^m\,y_{i'}^{j'}\,\frac{\partial}{\partial y^{j'}}.
\end{equation}
Striclty speaking, these formulas~\thetag{ 1.17} are not explicit,
because an inverse matrix is involved and because the terms $D_{i'}^1
X^i$, $D_{i'}^1 Y^j$ are not developed. However, it would be
elementary to write down the corresponding totally explicit complete
formulas for the functions $Y_{X^{i_1}}^j = Y_{X^{i_1}}^j \left(
x^{i'}, y^{j'}, y_{i_1'}^{j'} \right)$.

Next, the second prolongation $\varphi^{ (2)}$ is of the form
\def\theequation{1.19}\begin{equation}
\varphi^{(2)}
: \
\left(
x^{i'},y^{j'}, y_{i_1'}^{j'},
y_{i_1',i_2'}^{j'}
\right) 
\longmapsto 
\left(
\varphi^{(1)}
\left(
x^{i'},y^{j'}, y_{i_1'}^{j'}
\right), 
\
Y_{X^{i_1}X^{i_2}}^j
\left(
x^{i'}, y^{j'}, y_{i_1'}^{j'}, y_{i_1',i_2'}^{j'}
\right)
\right),
\end{equation}
for some functions $Y_{X^{i_1} X^{i_2}}^j \left( x^{ i'}, y^{ j'}, y_{
i_1' }^{j'}, y_{i_1', i_2'}^{j'} \right)$ which depend on the pure
first and second jet variables. For $i = 1, \dots, n$, the
expressions of $Y_{X^{i_1}X^i}^j$ are given by the following compact
formulas (again \cite{ bk1989}):
\def\theequation{1.20}\begin{equation}
\left(
\begin{array}{c}
Y_{X^{i_1}X^1}^j \\
\vdots \\
Y_{X^{i_1}X^n}^j \\
\end{array}
\right)
=
\left(
\begin{array}{ccc}
D_1^1 X^1 & \cdots & D_1^1 X^n \\
\vdots & \cdots & \vdots \\
D_n^1 X^1 & \cdots & D_n^1 X^n \\
\end{array} 
\right)^{-1}
\left(
\begin{array}{c}
D_1^2 Y_{X^{i_1}}^j \\
\vdots \\
D_n^2 Y_{X^{i_1}}^j \\
\end{array}
\right), 
\end{equation}
where, for $i' = 1, \dots, n$, the $D_{i'}^2$ denote the {\sl 
$i'$-th second order
total differentiation operators}:
\def\theequation{1.21}\begin{equation}
D_{i'}^2
:= 
\frac{\partial}{\partial x^{i'}}
+
\sum_{j'=1}^m\,y_{i'}^{j'}\,\frac{\partial}{\partial y^{j'}}
+
\sum_{j'=1}^m\,\sum_{i_1'=1}^n\,y_{i',i_1'}^{j'}\,
\frac{\partial}{\partial y_{i_1'}^{j'}}.
\end{equation}
Again, these formulas~\thetag{ 1.20} are not explicit in the sense
that an inverse matrix is involved and the terms $D_{i'}^1 X^i$,
$D_{i'}^2 Y_{ X^{ i_1}}^j$ are not developed. It would already be a
nontrivial computational task to develope these expressions and to
find out nice satisfying combinatorial formulas.

In order to present the general inductive non-explicit formulas for
the computation of the $\kappa$-th prolongation $\varphi^{(\kappa)}$,
we need some more notation. Let $\lambda \in \N$ be an arbitrary
integer. For $i' = 1, \dots, n$, let $D_{i'}^\lambda$ denotes the {\sl
$i'$-th $\lambda$-th order total differentiation operators}, defined
precisely by:
\def\theequation{1.22}\begin{equation}
\left\{
\aligned
D_{i'}^\lambda
&
:=
\frac{\partial}{\partial x^{i'}}
+
\sum_{j'=1}^m\,y_{i'}^{j'}\,\frac{\partial}{\partial y^{j'}}
+
\sum_{j'=1}^m\,\sum_{i_1'=1}^n\,y_{i',i_1'}^{j'}\,
\frac{\partial}{\partial y_{i_1'}^{j'}}
+
\sum_{j'=1}^m\,\sum_{i_1',i_2'=1}^n\,y_{i',i_1',i_2'}^{j'}\,
\frac{\partial}{\partial y_{i_1',i_2'}^{j'}}
+ \\
& \
\ \ \ \ \
+
\cdots
+
\sum_{j'=1}^m\,\sum_{i_1',i_2',\dots,i_{\lambda-1}'=1}^n\,
y_{i',i_1',i_2',\dots,i_{\lambda-1}'}^{j'}\,
\frac{\partial}{\partial y_{i_1',i_2',\dots,i_{\lambda-1}'}^{j'}}.
\endaligned\right.
\end{equation} 
Then, for $i = 1, \dots, n$, the expressions of $Y_{X^{i_1}\cdots
X^{i_{\lambda-1}} X^i }^j$ are given by the following compact formulas
(again \cite{ bk1989}):
\def\theequation{1.23}\begin{equation}
\left(
\begin{array}{c}
Y_{X^{i_1}\cdots X^{i_{\lambda-1}}X^1}^j \\
\vdots \\
Y_{X^{i_1}\cdots X^{i_{\lambda-1}}X^n}^j \\
\end{array}
\right)
=
\left(
\begin{array}{ccc}
D_1^1 X^1 & \cdots & D_1^1 X^n \\
\vdots & \cdots & \vdots \\
D_n^1 X^1 & \cdots & D_n^1 X^n \\
\end{array} 
\right)^{-1}
\left(
\begin{array}{c}
D_1^\lambda Y_{X^{i_1}\cdots X^{i_{\lambda-1}}}^j \\
\vdots \\
D_n^\lambda Y_{X^{i_1}\cdots X^{i_{\lambda-1}}}^j \\
\end{array}
\right).
\end{equation}
Again, these inductive formulas are incomplete and unsatisfactory.

\def\theproblem{1.24}\begin{problem}
Find totally explicit complete formulas
for the $\kappa$-th prolongation $\varphi^{(\kappa)}$.
\end{problem}

Except in the cases $\kappa = 1, 2$, we have not been able to solve
this problem. The case $\kappa = 1$ is elementary. Complete formulas
in the particular cases $\kappa = 2$, $n=1$, $m\geq 1$ and $n\geq 1$,
$m=1$ are implicitely provided in~\cite{ m2004a} and in~\cite{
m2004b}, where one observes the appearance of some modifications of
the Jacobian determinant of the diffeomorphism $\varphi$, inserted in
a clearly understandable combinatorics. In fact, there is a nice
dictionary between the formulas for $\varphi^{ (2)}$ and the formulas
for the second prolongation $\mathcal{ L}^{ (2)}$ of a vector field
$\mathcal{ L}$ which were written in equation~\thetag{ 43} of~\cite{
gm2003} ({\it see} also equations~\thetag{ 2.6}, \thetag{ 3.20},
\thetag{ 4.6} and~\thetag{ 5.3} in the next paragraphs). In the
passage from $\varphi^{ (2)}$ to $\mathcal{ L}^{ (2)}$, a sort of
formal first order linearization may be observed and the reverse
passage may be easily guessed. However, for $\kappa \geq 3$, the
formulas for $\varphi^{ (\kappa)}$ explode faster than the formulas
for the $\kappa$-th prolongation $\mathcal{ L}^{ (\kappa )}$ of a
vector field $\mathcal{ L}$. Also, the dictionary between $\varphi^{
(\kappa)}$ and $\mathcal{ L}^{ ( \kappa )}$ disappears. In fact, to
elaborate an appropriate dictionary, we believe that one should
introduce before a sort of formal $(\kappa-1)$ order linearizations of
$\varphi^{ ( \kappa)}$, finer than the first order linearization
$\mathcal{ L}^{ (\kappa)}$. To be optimistic, we believe that the
final answer to Problem~1.24
is accessible.

The present article is devoted to present totally explicit complete
formulas for the $\kappa$-th prolongation $\mathcal{ L}^{ (\kappa )}$
of a vector field $\mathcal{ L}$ to $\mathcal{ J}_{n,m}^\kappa$, for
$n\geq 1$ arbitrary, for $m\geq 1$ arbitrary and for $\kappa \geq 1$
arbitrary.

\subsection*{ 1.25.~Prolongation of a vector field to the
$\kappa$-th jet space} Consider a vector field
\def\theequation{1.26}\begin{equation}
\mathcal{ L} 
= 
\sum_{i=1}^n\mathcal{X}^i\,\frac{\partial}{\partial x^i}
+
\sum_{j=1}^m\,\mathcal{Y}^j\,\frac{\partial}{\partial y^j},
\end{equation}
defined in $\K^{ n+m}$. Its flow:
\def\theequation{1.27}\begin{equation}
\varphi_t ( x, y) 
:=
\exp 
\left(
t \mathcal{ L} 
\right) 
(x, y)
\end{equation}
constitutes a one-parameter of diffeomorphisms of $\K^{ n+m} $ close
to the identity. The lift $(\varphi_t )^{ ( \kappa )}$ to the
$\kappa$-th jet space constitutes a one-parameter family of
diffeomorphisms of $\mathcal{ J}_{ n,m }^\kappa$. By definition, the
{\sl $\kappa$-th prolongation $\mathcal{ L }^{( \kappa)}$ of
$\mathcal{ L }$ to the jet space $\mathcal{ J}_{n, m }^\kappa$} is the
infinitesimal generator of $(\varphi_t )^{ (\kappa)}$, namely:
\def\theequation{1.28}\begin{equation}
\mathcal{ L}^{(\kappa)}
:= 
\left.
\frac{d}{dt}
\right\vert_{t=0}
\left[
(\varphi_t)^{(\kappa)}
\right].
\end{equation}

\subsection*{1.29.~Inductive formulas for the 
$\kappa$-th prolongation $\mathcal{ L}^{(\kappa)}$} As a vector field
defined in $\K^{n+m+ m \frac{ (n+m)! }{ n! \ m!}}$, the $\kappa$-th
prolongation $\mathcal{ L}^{ (\kappa) }$ may be written under the
general form:
\def\theequation{1.30}\begin{equation}
\left\{
\aligned
\mathcal{L}^{(\kappa)}
&
=
\sum_{i=1}^n\mathcal{X}^i\,\frac{\partial}{\partial x^i}
+
\sum_{j=1}^m\,\mathcal{Y}^j\,\frac{\partial}{\partial y^j}
+ \\
& \
\ \ \ \ \
+
\sum_{j=1}^m\,\sum_{i_1=1}^n\,{\bf Y}_{i_1}^j\,
\frac{\partial}{\partial y_{i_1}^j}
+
\sum_{j=1}^m\,\sum_{i_1,i_2=1}^n\,{\bf Y}_{i_1,i_2}^j\,
\frac{\partial}{\partial y_{i_1,i_2}^j}
+ 
\cdots
+ \\
& \
\ \ \ \ \ 
+
\sum_{j=1}^m\,\sum_{i_1,\dots,i_\kappa=1}^n\,
{\bf Y}_{i_1,\dots,i_\kappa}^j\,
\frac{\partial}{\partial y_{i_1,\dots,i_\kappa}^j}.
\endaligned\right.
\end{equation}
Here, the coefficients ${\bf Y}_{i_1}^j$, ${\bf Y}_{i_1, i_2}^j$,
$\dots$, ${\bf Y}_{i_1, i_2, \dots, i_\kappa}^j$ are uniquely
determined in terms of partial derivatives of the coefficients
$\mathcal{ X}^i$ and $\mathcal{ Y}^j$ of the original vector field
$\mathcal{ L}$, together with the pure jet variables $\left( y_{ i_1
}^j,\dots, y_{i_1, \dots, i_\kappa}^j \right)$, by means of the
following {\sl fundamental inductive formulas} (\cite{ ol1979}, 
\cite{ ol1986}, \cite{
bk1989}):
\def\theequation{1.31}\begin{equation}
\left\{
\aligned
{\bf Y}_{i_1}^j
&
:=
D_{i_1}^1
\left(
\mathcal{ Y}^j
\right)
-
\sum_{k=1}^n\,D_{i_1}^1
\left(
\mathcal{X}^k
\right)
\, y_k^j, 
\\
{\bf Y}_{i_1,i_2}^j
&
:=
D_{i_2}^2
\left(
{\bf Y}_{i_1}^j
\right)
-
\sum_{k=1}^n\,D_{i_2}^1
\left(
\mathcal{X}^k
\right)
\, y_{i_1,k}^j, 
\\
\cdots \cdots \cdots 
&
\cdots \cdots \cdots \cdots \cdots \cdots \cdots \cdots
\cdots \cdots \cdots \cdots \cdots \cdots \cdots \cdots
\\
{\bf Y}_{i_1,i_2,\dots,i_\kappa}^j
&
:=
D_{i_\kappa}^\kappa
\left(
{\bf Y}_{i_1,i_2,\dots,i_{\kappa-1}}^j
\right)
-
\sum_{k=1}^n\,D_{i_\kappa}^1
\left(
\mathcal{X}^k
\right)
\, y_{i_1,i_2,\dots,i_{\kappa-1},k}^j,
\endaligned\right.
\end{equation}
where, for every $\lambda \in \N$ with $0 \leq \lambda \leq \kappa$,
and for every $i\in \N$ with 
$1 \leq i \leq n$, the $i$-th $\lambda$-th order total
differentiation operator $D_i^\lambda$ was defined in~\thetag{ 
1.22}
above.

\def\theproblem{1.32}\begin{problem}
Applying these inductive formulas, find totally explicit complete
formulas for the $\kappa$-th prolongation $\mathcal{L}^{(\kappa)}$.
\end{problem}

The present article is devoted to provide all the desired formulas.

\subsection*{ 1.33.~Inductive methodology}
We have the intention of presenting our results in a purely inductive
style, based on several thorough visual comparisons between massive formulas
which will be written and commented in four different cases:

\smallskip

\begin{itemize}
\item[{\bf (i)}]
$n=1$ and $m=1$; $\kappa\geq 1$ arbitrary;
\item[{\bf (ii)}]
$n\geq 1$ and $m=1$; $\kappa\geq 1$ arbitrary;
\item[{\bf (iii)}]
$n=1$ and $m\geq 1$; $\kappa\geq 1$ arbitrary;
\item[{\bf (iv)}]
general case: $n\geq 1$ and $m\geq 1$; $\kappa\geq 1$ arbitrary.
\end{itemize}

\smallskip

Accordingly, we shall particularize and slightly lighten our notations
in each of the three (preliminary) cases (i)
[Section~2], (ii) [Section~3] and (iii) [Section~4].

\section*{\S2.~One independent variable and one dependent variable}

\subsection*{2.1.~Simplified adapted notations}
Assume $n=1$ and $m=1$, let $\kappa\in \N$ with $\kappa \geq 1$ and
simply denote the jet variables by:
\def\theequation{2.2}\begin{equation}
\left(
x,y,y_1,y_2,\dots,y_\kappa
\right) \in \mathcal{J}_{1,1}^\kappa.
\end{equation}
The $\kappa$-th prolongation of a vector field will be denoted by:
\def\theequation{2.3}\begin{equation}
\mathcal{L}
=
\mathcal{X}\,\frac{\partial}{\partial x}
+
\mathcal{Y}\,\frac{\partial}{\partial y}
+
{\bf Y}_1\,\frac{\partial}{\partial y_1}
+
{\bf Y}_2\,\frac{\partial}{\partial y_2}
+ 
\cdots
+
{\bf Y}_\kappa\,\frac{\partial}{\partial y_\kappa}.
\end{equation}
The coefficients 
${\bf Y}_1$, ${\bf Y}_2$, $\dots$, ${\bf Y}_\kappa$
are computed by means of the inductive formulas:
\def\theequation{2.4}\begin{equation}
\left\{
\aligned
{\bf Y}_1
&
:= 
D^1(\mathcal{Y})
-
D^1(\mathcal{X})\,y_1, \\
{\bf Y}_2
&
:= 
D^2({\bf Y}_1)
-
D^1(\mathcal{X})\,y_2, \\
\cdots 
&
\cdots \cdots \cdots \cdots \cdots \cdots \cdots \cdots
\\
{\bf Y}_\kappa
&
:= 
D^\kappa({\bf Y}_{\kappa-1})
-
D^1(\mathcal{X})\,y_\kappa, \\
\endaligned\right.
\end{equation}
where, for $1 \leq \lambda \leq \kappa$:
\def\theequation{2.5}\begin{equation}
D^\lambda
:= 
\frac{\partial}{\partial x}
+
y_1\,\frac{\partial}{\partial y}
+ 
y_2\,\frac{\partial}{\partial y_1}
+
\cdots
+
y_\lambda\,\frac{\partial}{\partial y_{\lambda-1}}.
\end{equation}
By direct elementary computations, for $\kappa = 1$ and for $\kappa =
2$, we obtain the following two very classical formulas :
\def\theequation{2.6}\begin{equation}
\left\{
\aligned
{\bf Y}_1
&
=
\mathcal{Y}_x
+
\left[
\mathcal{Y}_y
-
\mathcal{X}_x
\right]
y_1
+
\left[
-\mathcal{Y}_y
\right]
(y_1)^2, 
\\
{\bf Y}_2
&
=
\mathcal{Y}_{x^2}
+
\left[
2\,\mathcal{Y}_{xy}
-
\mathcal{X}_{x^2}
\right]
\,y_1
+
\left[
\mathcal{Y}_{y^2}
-
2\,\mathcal{X}_{xy}
\right](y_1)^2
+
\left[
-
\mathcal{X}_{y^2}
\right]
(y_1)^3
+ \\
& \
\ \ \ \ \
+
\left[
\mathcal{Y}_y
-
2\,\mathcal{X}_x
\right]\,y_2
+
\left[
-
3\,\mathcal{X}_y
\right]\,
y_1\,y_2.
\endaligned\right.
\end{equation}
Our main objective is to {\it devise the general combinatorics}. 
Thus, to
attain this aim, we have to achieve patiently formal computations of
the next coefficients ${\bf Y}_3$, ${\bf Y}_4$ and ${\bf Y}_5$. We
systematically use parentheses $\left[ \cdot \right]$ to single out
every coefficient of the polynomials ${\bf Y}_3$, ${\bf Y}_4$
and ${\bf
Y}_5$ in the pure jet variables $y_1, y_2, y_3, y_4$ and
$y_5$, putting every sign inside these parentheses. We always put
the monomials in the pure jet variables $y_1, y_2, y_3, y_4$ and
$y_5$ after the parentheses. For completeness, let us provide the
intermediate computation of the third coefficient ${\bf Y}_3$.
In detail:
$$
\small
\aligned
{\bf Y}_3
&
=
D^3
\left(
{\bf Y}_2
\right)
-
D^1
\left(
\mathcal{ X}
\right)
y_3 
\\
&
=
\left(
\frac{\partial}{\partial x}
+
y_1\,\frac{\partial}{\partial y}
+
y_2\,\frac{\partial}{\partial y_1}
+
y_3\,\frac{\partial}{\partial y_2}
\right)
\left(
\mathcal{Y}_{x^2}
+
\left[
2\,\mathcal{Y}_{xy}
-
\mathcal{X}_{x^2}
\right]
y_1
+ 
\right.
\\
& \
\ \ \ \ \
\left.
+
\left[
\mathcal{ Y}_{y^2}
-
2\,\mathcal{X}_{xy}
\right]
(y_1)^2
+
\left[
-
\mathcal{X}_{y^2}
\right]
(y_1)^3
+
\left[
\mathcal{ Y}_y
-
2\,\mathcal{X}_x
\right]
y_2
+
\left[
-
3\,\mathcal{X}_y
\right]
y_1\,y_2
\right)
\\
\endaligned
$$
\def\theequation{2.7}\begin{equation}
\aligned
&
=
\underline{
\mathcal{Y}_{x^3} 
}_{ \fbox{\tiny 1}}
+
\underline{
\left[
2\,\mathcal{Y}_{x^2y}
-
\mathcal{X}_{x^3}
\right]
y_1
}_{ \fbox{\tiny 2}}
+
\underline{
\left[
\mathcal{Y}_{xy^2}
-
2\,\mathcal{X}_{x^2y}
\right]
(y_1)^2
}_{ \fbox{\tiny 3}}
+
\underline{
\left[
-
\mathcal{X}_{xy^2}
\right]
(y_1)^3
}_{ \fbox{\tiny 4}}
+ \\
& \
\ \ \ \ \
+
\underline{ 
\left[
\mathcal{Y}_{xy}
-
2\,\mathcal{X}_{x^2}
\right]
y_2
}_{ \fbox{\tiny 6}}
+
\underline{ 
\left[
-
3\,\mathcal{X}_{xy}
\right]
y_1y_2
}_{ \fbox{\tiny 7}}
+
\underline{ 
\left[
\mathcal{Y}_{x^2y}
\right]
y_1
}_{ \fbox{\tiny 2}}
+ \\
& \
\ \ \ \ \
+
\underline{ 
\left[
2\,\mathcal{Y}_{xy^2}
-
\mathcal{X}_{x^2y}
\right]
(y_1)^2
}_{ \fbox{\tiny 3}}
+
\underline{ 
\left[
\mathcal{Y}_{y^3}
-
2\,\mathcal{X}_{xy^2}
\right]
(y_1)^3
}_{ \fbox{\tiny 4}}
+
\underline{ 
\left[
-
\mathcal{X}_{y^3}
\right]
(y_1)^4
}_{ \fbox{\tiny 5}}
+
\endaligned
\end{equation}
$$
\aligned
& \
\ \ \ \ \
+
\underline{ 
\left[
\mathcal{Y}_{y^2}
-
2\,\mathcal{X}_{xy}
\right]
y_1y_2
}_{ \fbox{\tiny 7}}
+
\underline{ 
\left[
-
3\,\mathcal{X}_{y^2}
\right]
(y_1)^2y_2
}_{ \fbox{\tiny 8}}
+
\underline{ 
\left[
2\,\mathcal{Y}_{xy}
-
\mathcal{X}_{x^2}
\right]
y_2
}_{ \fbox{\tiny 6}}
+ \\
& \
\ \ \ \ \
+
\underline{ 
\left[
\mathcal{Y}_{y^2}
-
2\,\mathcal{X}_{xy}
\right]
2\,y_1y_2
}_{ \fbox{\tiny 7}}
+
\underline{ 
\left[
-
\mathcal{X}_{y^2}
\right]
3(y_1)^2y_2
}_{ \fbox{\tiny 8}}
+
\underline{ 
\left[
-
3\,\mathcal{X}_y
\right]
(y_2)^2
}_{ \fbox{\tiny 9}}
+ \\
& \
\ \ \ \ \ 
+
\underline{ 
\left[
\mathcal{Y}_y
-
2\,\mathcal{X}_x
\right]
y_3
}_{ \fbox{\tiny 10}}
+
\underline{ 
\left[
-
3\,\mathcal{X}_y
\right]
y_1y_3
}_{ \fbox{\tiny 11}}
- \\
& \
\ \ \ \ \
-
\underline{ 
\left[
\mathcal{X}_x
\right]
y_3
}_{ \fbox{\tiny 10}}
-
\underline{ 
\left[
\mathcal{X}_y
\right]
y_1y_3
}_{ \fbox{\tiny 11}} \ .
\endaligned
$$
We have underlined all the terms with 
a number appended. Each number
refers to the order of appearance
of the terms in the final simplified
expression of ${\bf Y}_3$, also written in~\cite{ bk1989}
with different notations:
\def\theequation{2.8}\begin{equation}
\left\{
\aligned
{\bf Y}_3
& 
=
\mathcal{Y}_{x^3}
+
\left[
3\,\mathcal{Y}_{x^2y}
-
\mathcal{X}_{x^3}
\right]
y_1
+
\left[
3\,\mathcal{Y}_{xy^2}
-
3\,\mathcal{X}_{x^2y}
\right]
(y_1)^2
+ \\
& \
\ \ \ \ \
+
\left[
\mathcal{Y}_{y^3}
-
3\,\mathcal{X}_{xy^2}
\right]
(y_1)^3
+
\left[
-
\mathcal{X}_{y^3}
\right]
(y_1)^4
+
\left[
3\,\mathcal{Y}_{xy}
-
3\,\mathcal{X}_{x^2}
\right]
y_2
+ \\
& \
\ \ \ \ \
+
\left[
3\,\mathcal{Y}_{y^2}
-
9\,\mathcal{X}_{xy}
\right]
y_1y_2
+
\left[
-
6\,\mathcal{X}_{y^2}
\right]
(y_1)^2y_2
+
\left[
-
3\,\mathcal{X}_y
\right]
(y_2)^2
+ \\
& \
\ \ \ \ \
+
\left[
\mathcal{Y}_y
-
3\,\mathcal{X}_x
\right]
y_3
+
\left[
-
4\,\mathcal{X}_y
\right]
y_1y_3.
\endaligned\right.
\end{equation}
After similar manual computations, the intermediate details of which
we will not copy in this Latex file, we get the desired expressions of
${\bf Y}_4$ and of ${\bf Y}_5$.
Firstly:
\def\theequation{2.9}\begin{equation}
\small
\left\{
\aligned
{\bf Y}_4
&
=
\mathcal{Y}_{x^4}
+
\left[
4\,\mathcal{Y}_{x^3y}
-
\mathcal{X}_{x^4}
\right]
y_1
+
\left[
6\,\mathcal{Y}_{x^2y^2}
-
4\,\mathcal{X}_{x^3y}
\right]
(y_1)^2
+ \\
& \
\ \ \ \ \
+
\left[
4\,\mathcal{Y}_{xy^3}
-
6\,\mathcal{X}_{x^2y^2}
\right]
(y_1)^3
+
\left[
\mathcal{Y}_{y^4}
-
4\,\mathcal{X}_{xy^3}
\right]
(y_1)^4
+
\left[
-
\mathcal{X}_{y^4}
\right]
(y_1)^5
+ \\
& \
\ \ \ \ \
+
\left[
6\,\mathcal{Y}_{x^2y}
-
4\,\mathcal{X}_{x^3}
\right]
y_2
+
\left[
12\,\mathcal{Y}_{xy^2}
-
18\,\mathcal{X}_{x^2y}
\right]
y_1y_2
+ \\
& \
\ \ \ \ \
+
\left[
6\,\mathcal{Y}_{y^3}
-
24\,\mathcal{X}_{xy^2}
\right]
(y_1)^2y_2
+
\left[
-
10\,\mathcal{X}_{y^3}
\right]
(y_1)^3y_2
+ \\
& \
\ \ \ \ \
+
\left[
3\,\mathcal{Y}_{y^2}
-
12\,\mathcal{X}_{xy}
\right]
(y_2)^2
+
\left[
-
15\,\mathcal{X}_{y^2}
\right]
y_1(y_2)^2
+ \\
& \
\ \ \ \ \
+
\left[
4\,\mathcal{Y}_{xy}
-
6\,\mathcal{X}_{x^2}
\right]
y_3
+
\left[
4\,\mathcal{Y}_{y^2}
-
16\,\mathcal{X}_{xy}
\right]
y_1y_3
+
\left[
-
10\,\mathcal{X}_{y^2}
\right]
(y_1)^2y_3
+ \\
& \
\ \ \ \ \
+
\left[
-
10\,\mathcal{X}_y
\right]
y_2y_3
+
\left[
\mathcal{Y}_y
-
4\,\mathcal{X}_x
\right]
y_4
+
\left[
-
5\,\mathcal{X}_y
\right]
y_1y_4.
\endaligned\right.
\end{equation}
Secondly:
\def\theequation{2.10}\begin{equation}
\small
\left\{
\aligned
{\bf Y}_5
&
=
\mathcal{Y}_{x^5}
+
\left[
5\,\mathcal{Y}_{x^4y}
-
\mathcal{X}_{x^5}
\right]
y_1
+
\left[
10\,\mathcal{Y}_{x^3y^2}
-
5\,\mathcal{X}_{x^4y}
\right]
(y_1)^2
+ \\
& \
\ \ \ \ \
+
\left[
10\,\mathcal{Y}_{x^2y^3}
-
10\,\mathcal{X}_{x^3y^2}
\right]
(y_1)^3
+
\left[
5\,\mathcal{Y}_{xy^4}
-
10\,\mathcal{X}_{x^2y^3}
\right]
(y_1)^4
+ \\
& \
\ \ \ \ \
+
\left[
\mathcal{Y}_{y^5}
-
5\,\mathcal{X}_{xy^4}
\right]
(y_1)^5
+
\left[
-
\mathcal{X}_{y^5}
\right]
(y_1)^6
+
\left[
10\,\mathcal{Y}_{x^3y}
-
5\,\mathcal{X}_{x^4}
\right]
y_2
+ \\
& \
\ \ \ \ \
+
\left[
30\,\mathcal{Y}_{x^2y^2}
-
30\,\mathcal{X}_{x^3y}
\right]
y_1y_2
+
\left[
30\,\mathcal{Y}_{xy^3}
-
60\,\mathcal{X}_{x^2y^2}
\right]
(y_1)^2y_2
+ \\
& \
\ \ \ \ \
+
\left[
10\,\mathcal{Y}_{y^4}
-
50\,\mathcal{X}_{xy^3}
\right]
(y_1)^3y_2
+
\left[
-
15\,\mathcal{X}_{y^4}
\right]
(y_1)^4y_2
+ \\
& \
\ \ \ \ \
+
\left[
15\,\mathcal{Y}_{xy^2}
-
30\,\mathcal{X}_{x^2y}
\right]
(y_2)^2
+
\left[
15\,\mathcal{Y}_{y^3}
-
75\,\mathcal{X}_{xy^2}
\right]
y_1(y_2)^2
+ \\
& \
\ \ \ \ \
+
\left[
-
45\,\mathcal{X}_{y^3}
\right]
(y_1)^2(y_2)^2
+
\left[
-
15\,\mathcal{X}_{y^2}
\right]
(y_2)^3
+ \\
& \
\ \ \ \ \
+
\left[
10\,\mathcal{Y}_{x^2y}
-
10\,\mathcal{X}_{x^3}
\right]
y_3
+
\left[
20\,\mathcal{Y}_{xy^2}
-
40\,\mathcal{X}_{x^2y}
\right]
y_1y_3
+ \\
& \
\ \ \ \ \
+
\left[
10\,\mathcal{Y}_{y^3}
-
50\,\mathcal{X}_{xy^2}
\right]
(y_1)^2y_3
+
\left[
-
20\,\mathcal{X}_{y^3}
\right]
(y_1)^3y_3
+ \\
& \
\ \ \ \ \
+
\left[
10\,\mathcal{Y}_{y^2}
-
50\,\mathcal{X}_{xy}
\right]
y_2y_3
+
\left[
-
60\,\mathcal{X}_{y^2}
\right]
y_1y_2y_3
+
\left[
-
10\,\mathcal{X}_y
\right]
(y_3)^2
+ \\
& \
\ \ \ \ \
+
\left[
5\,\mathcal{Y}_{xy}
-
10\,\mathcal{X}_{x^2}
\right]
y_4
+
\left[
5\,\mathcal{Y}_{y^2}
-
25\,\mathcal{X}_{xy}
\right]
y_1y_4
+
\left[
-
15\,\mathcal{X}_{y^2}
\right]
(y_1)^2y_4
+ \\
& \
\ \ \ \ \
+
\left[
-
15\,\mathcal{X}_y
\right]
y_2y_4
+
\left[
\mathcal{Y}_y
-
5\,\mathcal{X}_y
\right]
y_5
+
\left[
-
6\,\mathcal{X}_y
\right]
y_1y_5.
\endaligned\right.
\end{equation}

\subsection*{2.11.~Formal inspection, formal intuition
and formal induction} Now, we have to comment these formulas. We have
written in length the five polynomials ${\bf Y}_1$, ${\bf Y}_2$, ${\bf
Y}_3$, ${\bf Y}_4$ and ${\bf Y}_5$ in the pure jet variables $y_1,
y_2, y_3, y_4$ and $y_5$. Except the first ``constant'' term
$\mathcal{ Y}_{x^\kappa}$, all the monomials in the expression of
${\bf Y}_\kappa$ are of the general form
\def\theequation{2.12}\begin{equation}
\left(
y_{\lambda_1}
\right)^{\mu_1}
\left(
y_{\lambda_2}
\right)^{\mu_2}
\cdots
\left(
y_{\lambda_d}
\right)^{\mu_d},
\end{equation}
for some positive integer $d\geq 1$, for some collection of
strictly increasing jets indices:
\def\theequation{2.13}\begin{equation}
1 \leq \lambda_1 < \lambda_2 < \cdots < \lambda_d \leq \kappa,
\end{equation}
and for some positive integers $\mu_1, \dots, \mu_d \geq 1$. This and
the next combinatorial facts may be confirmed by reading the formulas
giving ${\bf Y}_1$, ${\bf Y}_2$, ${\bf Y}_3$, ${\bf Y}_4$ and ${\bf
Y}_5$. It follows that the integer $d$ satisfies the inequality
$d\leq \kappa+1$. To include the first ``constant'' term $\mathcal{
Y}_{x^\kappa}$, we shall make the convention that putting $d=0$ in the
monomial~\thetag{ 2.12} yields the constant term $1$.
 
Furthermore, by inspecting the formulas giving ${\bf Y}_1$, ${\bf
Y}_2$, ${\bf Y}_3$, ${\bf Y}_4$ and ${\bf Y}_5$, we see that the
following inequality should be satisfied:
\def\theequation{2.14}\begin{equation}
\mu_1\lambda_1
+
\mu_2\lambda_2
+
\cdots
+
\mu_d\lambda_d 
\leq 
\kappa+1.
\end{equation}
For instance, in the expression of ${\bf Y}_4$, the two monomials
$(y_1)^3 y_2$ and $y_1 (y_2)^2$ do appear, but the two monomials
$(y_1)^4 y_2$ and $(y_1)^2 (y_2)^2$ cannot appear. All coefficients
of the pure jet monomials are of the general form:
\def\theequation{2.15}\begin{equation}
\left[
A\,
\mathcal{Y}_{x^\alpha y^{\beta+1}}
-
B\,
\mathcal{X}_{x^{\alpha+1}y^\beta}
\right],
\end{equation}
for some nonnegative integers $A, B, \alpha, \beta \in \N$. Sometimes
$A$ is zero, but $B$ is zero only for the (constant, with respect to
pure jet variables) term $\mathcal{ Y}_{x^\kappa}$. Importantly,
$\mathcal{ X}$ is differentiated once more with respect to $x$ and
$\mathcal{ Y}$ is differentiated once more with respect to $y$. Again,
this may be confirmed by reading all the terms in the formulas for
${\bf Y}_1$, ${\bf Y}_2$, ${\bf Y}_3$, ${\bf Y}_4$ and ${\bf Y}_5$.

In addition, we claim that there is a link between the couple
$(\alpha, \beta)$ and the collection $\{ \mu_1, \lambda_1, \dots,
\mu_d, \lambda_d \}$. To discover it, let us write some of the
monomials appearing in the expressions of ${\bf Y}_4$ (first column)
and of ${\bf Y}_5$ (second column), for instance:
\def\theequation{2.16}\begin{equation}
\left\{
\aligned
&
\left[6\,
\mathcal{Y}_{x^2y^2}
-
4\,\mathcal{X}_{x^3y}
\right]
(y_1)^2, 
\ \ \ \ \ \ \ \ \ \ \ \ \ \ \
&
\left[
5\,\mathcal{Y}_{xy^4}
-
10\,\mathcal{X}_{x^2y^3}
\right]
(y_1)^4,
&
\\
&
\left[
12\,\mathcal{Y}_{xy^2}
-
18\,\mathcal{X}_{x^2y}
\right]
y_1y_2,
\ \ \ \ \ \ \ \ \ \ \ \ \ \ \
&
\left[
30\,\mathcal{Y}_{xy^3}
-
60\,\mathcal{X}_{x^2y^2}
\right]
(y_1)^2y_2,
& 
\\
& 
\left[
-
10\,\mathcal{X}_{y^3}
\right]
(y_1)^3y_2, 
\ \ \ \ \ \ \ \ \ \ \ \ \ \ \
& 
\left[
-
15\,\mathcal{X}_{y^4}
\right]
(y_1)^4y_2,
&
\\
&
\left[
4\,\mathcal{Y}_{y^2}
-
16\,\mathcal{X}_{xy}
\right]
y_1y_3, 
\ \ \ \ \ \ \ \ \ \ \ \ \ \ \
& 
\left[
10\,\mathcal{Y}_{y^2}
-
50\,\mathcal{X}_{xy}
\right]
y_2y_3,
&
\\
&
\left[
-
10\,\mathcal{X}_{y^2}
\right]
(y_1)^2 y_3,
\ \ \ \ \ \ \ \ \ \ \ \ \ \ \
& 
\left[
-
60\,\mathcal{X}_{y^2}
\right]
y_1y_2y_3.
& \\
\endaligned\right.
\end{equation}
After some reflection, we discover the hidden intuitive rule: the
partial derivatives of $\mathcal{ Y}$ and of $\mathcal{ X}$ associated
with the monomial $(y_{\lambda_1})^{\mu_1} \cdots
(y_{\lambda_d})^{\mu_d}$ are, respectively:
\def\theequation{2.17}\begin{equation}
\left\{
\aligned
&
\mathcal{ Y}_{
x^{\kappa-\mu_1\lambda_1-\cdots-\mu_d\lambda_d}
\,
y^{\mu_1+\cdots+\mu_d}}, 
\\
&
\mathcal{ X}_{
x^{\kappa-\mu_1\lambda_1-\cdots-\mu_d\lambda_d+1}
\,
y^{\mu_1+\cdots+\mu_d-1}}. 
\\
\endaligned\right.
\end{equation}
This may be checked on each of the $10$ examples~\thetag{ 2.16}
above.

Now that we have explored and discovered the combinatorics of the pure
jet monomials, of the partial derivatives and of the complete sum
giving ${\bf Y}_\kappa$, we may express that it is of the following
general form:
\def\theequation{2.18}\begin{equation}
\left\{
\aligned
{\bf Y}_\kappa
&
=
\mathcal{ Y}_{x^\kappa}
+
\sum_{d=1}^{\kappa+1}
\ \
\sum_{1\leq\lambda_1<\cdots<\lambda_d\leq\kappa}
\ \
\sum_{\mu_1\geq 1,\dots,\mu_d\geq 1}
\
\sum_{
\mu_1\lambda_1
+
\cdots
+
\mu_d\lambda_d\leq \kappa+1} 
\\
& \
\ \ \ \ \
\left[
A_\kappa^{
(\mu_1, \lambda_1), \dots, (\mu_d, \lambda_d) }
\cdot
\mathcal{Y}_{
x^{\kappa-\mu_1\lambda_1-\cdots-\mu_d\lambda_d}
\,
y^{\mu_1+\cdots+\mu_d}
}
-
\right. \\
& \
\ \ \ \ \ \ \ \ \ \
\left.
-
B_\kappa^{
(\mu_1, \lambda_1), \dots, (\mu_d, \lambda_d) }
\cdot
\mathcal{X}_{
x^{\kappa-\mu_1\lambda_1-\cdots-\mu_d\lambda_d+1}
\,
y^{\mu_1+\cdots+\mu_d-1}
}
\right]
\cdot
\\
& \
\ \ \ \ \ \ \ \ \ \ \ \ \ \ \ \ \ \ \ \
\ \ \ \ \ \ \ \ \ \ \ \ \ \ \ \ \ \ \ \
\ \ \ \ \ \ \ \ \ \ \ \ \
\cdot
(y_{\lambda_1})^{\mu_1}
\cdots
(y_{\lambda_d})^{\mu_d}.
\endaligned\right.
\end{equation} 
Here, we separate the first term $\mathcal{ Y}_{x^\kappa}$ from the
general sum; it is the constant term in ${\bf Y}_\kappa$, which is a
polynomial with respect to the jet variables $y_\lambda$. In this
general formula, the only remaining unknowns are the nonnegative
integer coefficients $A_\kappa^{ (\mu_1, \lambda_1), \dots, (\mu_d,
\lambda_d) } \in \N$ and $B_\kappa^{ (\mu_1, \lambda_1), \dots,
(\mu_d, \lambda_d) } \in \N$. In Section~3 below, we shall explain how
we have discovered their exact value.

At present, even if we are unable to devise their explicit
expression, we may observe that the value of the special integer
coefficients $A^{(\mu_1, 1)}_{ \mu_1}$ and $B^{( \mu_1, 1)}_{ \mu_1}$
which are attached to the monomials ${\rm ct.}$, $y_1$, $(y_1)^2$,
$(y_1)^3$, $(y_1)^4$ and $(y_1)^5$ are simple. Indeed, by
inspecting the first terms in the expressions of ${\bf Y}_1$, ${\bf
Y}_2$, ${\bf Y}_3$, ${\bf Y}_4$ and ${\bf Y}_5$, we of course
recognize the binomial coefficients. In general:

\def\thelemma{2.19}\begin{lemma} For $\kappa \geq 1$, 
\def\theequation{2.20}\begin{equation}
\left\{
\aligned
{\bf Y}_\kappa
&
=
\mathcal{Y}_{x^\kappa}
+
\sum_{\lambda=1}^\kappa
\left[
\binom{\kappa}{\lambda}
\,\mathcal{Y}_{x^{\kappa-\lambda}y^\lambda}
-
\binom{\kappa}{\lambda-1}
\,\mathcal{X}_{x^{\kappa-\lambda+1}y^{\lambda-1}}
\right]
(y_1)^\lambda
+ \\
& \
\ \ \ \ \
+
\left[
-
\mathcal{X}_{y^\kappa}
\right]
(y_1)^\kappa
+
{\sf remainder},
\endaligned\right.
\end{equation}
where the term {\sf remainder} collects all remaining monomials in
the pure jet variables.
\end{lemma}

In addition, let us remind what we have observed and used in a
previous co-signed work.

\def\thelemma{2.21}\begin{lemma}
\text{\rm (\cite{ gm2003}, p.~536)} For $\kappa \geq 4$, nine among
the monomials of ${\bf Y}_\kappa$ are of the following general
form{\rm :}
\def\theequation{2.22}\begin{equation}
\left\{
\aligned
{\bf Y}_\kappa
&
=
\mathcal{Y}_{x^\kappa}
+
\left[
C_\kappa^1 \, \mathcal{Y}_{x^{\kappa-1}y}
-
\mathcal{X}_{x^\kappa}
\right]
y_1
+
\left[
C_\kappa^2\,
\mathcal{Y}_{x^{\kappa-2}y}
-
C_\kappa^1 \,\mathcal{X}_{x^{\kappa-1}}
\right]
y_2
+ \\
& \
\ \ \ \ \
+
\left[
C_\kappa^2\,
\mathcal{Y}_{x^2y}
-
C_\kappa^3\,
\mathcal{X}_{x^3}
\right]
y_{\kappa-2}
+
\left[
C_\kappa^1 \,\mathcal{Y}_{xy}
-
C_\kappa^2\,
\mathcal{X}_{x^2}
\right]
y_{\kappa-1}
+ \\
& \
\ \ \ \ \
+
\left[
C_\kappa^1 \, \mathcal{Y}_{y^2}
-
\kappa^2 \, \mathcal{ X}_{xy}
\right]
y_1 y_{\kappa-1}
+
\left[
-C_\kappa^2 \, \mathcal{ X}_y
\right]
y_2y_{\kappa-1}
+ \\
& \
\ \ \ \ \
+
\left[
\mathcal{Y}_y
- 
C_\kappa^1\,\mathcal{X}_x
\right]
+
\left[
-C_{\kappa+1}^1 \, \mathcal{ X}_y
\right]
y_1 y_\kappa
+
{\sf remainder},
\endaligned\right.
\end{equation}
where the term {\sf remainder} denotes all the remaining monomials,
and where $C_\kappa^\lambda := \frac{ \kappa!}{(\kappa - \lambda)! \
\lambda !}$ is a notation for the binomial coefficient which occupies
less space in Latex ``equation mode'' than the classical notation
\def\theequation{2.23}\begin{equation}
\binom{\kappa}{\lambda}.
\end{equation}
\end{lemma}

Now, we state directly the final theorem, without further inductive or
intuitive information.

\def\thetheorem{2.24}\begin{theorem} 
For $\kappa \geq 1$, we have{\rm :}
\def\theequation{2.25}\begin{equation}
\boxed{
\aligned
{\bf Y}_\kappa
&
=
\mathcal{ Y}_{x^\kappa}
+
\sum_{d=1}^{\kappa+1}
\ \
\sum_{1\leq\lambda_1<\cdots<\lambda_d\leq\kappa}
\ \
\sum_{\mu_1\geq 1,\dots,\mu_d\geq 1} 
\
\sum_{
\mu_1\lambda_1
+
\cdots
+
\mu_d\lambda_d\leq \kappa+1} 
\\
& \
\ \ 
\left[
\frac{\kappa\cdots(\kappa-\mu_1\lambda_1-\cdots-\mu_d\lambda_d+1)}
{(\lambda_1!)^{\mu_1}\,\mu_1!
\cdots
(\lambda_d!)^{\mu_d}\,\mu_d!
}
\cdot
\mathcal{Y}_{
x^{\kappa-\mu_1\lambda_1-\cdots-\mu_d\lambda_d}
\,
y^{\mu_1+\cdots+\mu_d}
}
-
\right. \\
& \
\ \ \ \ \ \ \ \ \ \
\left.
-
\frac{\kappa\cdots(\kappa-\mu_1\lambda_1-\cdots-\mu_d\lambda_d+2)
(\mu_1\lambda_1+\cdots+\mu_d\lambda_d)}
{(\lambda_1!)^{\mu_1}\,\mu_1!
\cdots
(\lambda_d!)^{\mu_d}\,\mu_d!
}
\cdot
\right.
\\
& \
\ \ \ \ \ \ \ \ \ \ \ \ \ \ \ \
\cdot
\mathcal{X}_{
x^{\kappa-\mu_1\lambda_1-\cdots-\mu_d\lambda_d+1}
\,
y^{\mu_1+\cdots+\mu_d-1}
}
\Big]
(y_{\lambda_1})^{\mu_1}
\cdots
(y_{\lambda_d})^{\mu_d}.
\endaligned
}
\end{equation}
\end{theorem}

Once the correct theorem is formulated, its proof follows by
accessible induction arguments which will not be developed here. It is
better to continue through and to examine thorougly the case of
several variables, since it will help us considerably to explain how
we discovered the exact values of the integer coefficients $A_\kappa^{
(\mu_1, \lambda_1), \dots, (\mu_d, \lambda_d) }$ and $B_\kappa^{
(\mu_1, \lambda_1), \dots, (\mu_d, \lambda_d) }$.

\subsection*{2.26.~Verification and application} 
Before proceeding further, let us rapidly verify that the above
general formula~\thetag{ 2.25} is correct by inspecting two instances
extracted from ${\bf Y}_5$.

Firstly, the coefficient of $(y_1)^3 y_3$ in ${\bf Y}_5$ is obtained
by putting $\kappa = 5$, $d = 2$, $\lambda_1 = 1$, $\mu_1 = 3$,
$\lambda_2 = 3$ and $\mu_2 = 1$ in the general formula~\thetag{ 2.25},
which yields:
\def\theequation{2.27}\begin{equation} 
\left[ 
0
-
\frac{5\cdot 4\cdot 3\cdot 2\cdot 1\cdot 6}{ 
(1!)^3\ 3!\ (3!)^1\ 1!}
\,\mathcal{X}_{y^3} 
\right] = 
\left[ 
-
20\,\mathcal{X}_{y^3} 
\right].
\end{equation}
This value is the same as in the original formula~\thetag{ 2.10}:
confirmation.

Secondly, the coefficient of $y_1 (y_2)^2$ in ${\bf Y}_5$ is obtained
by $\kappa = 5$, $d = 2$, $\lambda_1 = 1$, $\mu_1 = 1$, $\lambda_2 =
2$ and $\mu_2 = 2$ in the general formula~\thetag{ 2.25}, which yields:
\def\theequation{2.28}\begin{equation}
\left[
\frac{5\cdot 4\cdot 3\cdot 2\cdot 1}{
(1!)^1\ 1!\ (2!)^2\ 2!}
\,\mathcal{Y}_{y^3}
-
\frac{5\cdot 4\cdot 3\cdot 2\cdot 5}{
(1!)^1\ 1!\ (2!)^2\ 2!}
\,\mathcal{X}_{xy^2}
\right]
= 
\left[
15\,\mathcal{Y}_{y^3}
-
75\,\mathcal{X}_{xy^2}
\right].
\end{equation}
This value is the same as in the original formula~\thetag{ 2.10};
again: confirmation.

Finally, applying our general formula~\thetag{ 2.25}, we deduce the value of
${\bf Y}_6$ {\it without having to use ${\bf Y}_5$ and the induction
formulas~\thetag{ 2.4}}, which shortens substantially the
computations. For the pleasure, we obtain:
\def\theequation{2.29}\begin{equation}
\small
\left\{
\aligned
{\bf Y}_6
&
=
\mathcal{Y}_{x^6}
+
\left[
6\,\mathcal{Y}_{x^5y}
-
\mathcal{X}_{x^6}
\right]
y_1
+
\left[
15\,\mathcal{Y}_{x^4y^2}
-
6\,\mathcal{X}_{x^5y}
\right]
(y_1)^2
+ \\
& \
\ \ \ \ \
+
\left[
20\,\mathcal{Y}_{x^3y^3}
-
15\,\mathcal{X}_{x^4y^2}
\right]
(y_1)^3
+
\left[
15\,\mathcal{Y}_{x^2y^4}
-
20\,\mathcal{X}_{x^3y^3}
\right]
(y_1)^4
+ \\
& \
\ \ \ \ \
+
\left[
6\,\mathcal{Y}_{xy^5}
-
15\,\mathcal{X}_{x^2y^4}
\right]
(y_1)^5
+
\left[
\mathcal{Y}_{y^6}
-
6\,\mathcal{X}_{xy^5}
\right]
(y_1)^6
+
\left[
-
\mathcal{X}_{y^6}
\right]
(y_1)^7
+ \\
& \
\ \ \ \ \
+
\left[
15\,\mathcal{Y}_{x^4y}
-
6\,\mathcal{X}_{x^5}
\right]
y_2
+
\left[
60\,\mathcal{Y}_{x^3y^2}
-
45\,\mathcal{X}_{x^4y}
\right]
y_1y_2
+ \\
& \
\ \ \ \ \
+
\left[
90\,\mathcal{Y}_{x^2y^3}
-
120\,\mathcal{X}_{x^3y^2}
\right]
(y_1)^2y_2
+
\left[
60\,\mathcal{Y}_{xy^4}
-
150\,\mathcal{X}_{x^2y^3}
\right]
(y_1)^3y_2
+ \\
& \
\ \ \ \ \ 
+
\left[
15\,\mathcal{Y}_{y^5}
-
90\,\mathcal{X}_{xy^4}
\right]
(y_1)^4y_2
+
\left[
-
21\,\mathcal{X}_{y^5}
\right]
(y_1)^5y_2
+ \\
& \
\ \ \ \ \
+
\left[
45\,\mathcal{Y}_{x^2y^2}
-
60\,\mathcal{X}_{x^3y}
\right]
(y_2)^2
+
\left[
90\,\mathcal{Y}_{xy^3}
-
225\,\mathcal{X}_{x^2y^2}
\right]
y_1(y_2)^2
+ \\
& \
\ \ \ \ \
+
\left[
45\,\mathcal{Y}_{y^4}
-
270\,\mathcal{X}_{xy^3}
\right]
(y_1)^2(y_2)^2
+
\left[
-
210\,\mathcal{X}_{y^4}
\right]
(y_1)^3(y_2)^2
+ \\
& \
\ \ \ \ \
+
\left[
15\,\mathcal{Y}_{y^3}
-
90\,\mathcal{X}_{xy^2}
\right]
(y_2)^3
+
\left[
-
105\,\mathcal{X}_{y^3}
\right]
y_1(y_2)^3
+ \\
& \
\ \ \ \ \
+
\left[
20\,\mathcal{Y}_{x^3y}
-
15\,\mathcal{X}_{x^4}
\right]
y_3
+
\left[
60\,\mathcal{Y}_{x^2y^2}
-
80\,\mathcal{X}_{x^3y}
\right]
y_1y_3
+ \\
& \
\ \ \ \ \
+
\left[
60\,\mathcal{Y}_{xy^3}
-
150\,\mathcal{X}_{x^2y^2}
\right]
(y_1)^2y_3
+
\left[
20\,\mathcal{Y}_{y^4}
-
120\,\mathcal{X}_{xy^3}
\right]
(y_1)^3y_3
+ \\
& \
\ \ \ \ \
+
\left[
-
35\,\mathcal{X}_{y^4}
\right]
(y_1)^4y_3
+
\left[
60\,\mathcal{Y}_{xy^2}
-
150\,\mathcal{X}_{x^2y}
\right]
y_2y_3
+ \\
& \
\ \ \ \ \
+
\left[
60\,\mathcal{Y}_{y^3}
-
360\,\mathcal{X}_{xy^2}
\right]
y_1y_2y_3
+
\left[
-
210\,\mathcal{X}_{y^3}
\right]
(y_1)^2y_2y_3
+ \\
& \
\ \ \ \ \
+
\left[
-
105\,\mathcal{X}_{y^2}
\right]
(y_2)^2y_3
+
\left[
10\,\mathcal{Y}_{y^2}
-
60\,\mathcal{X}_{xy}
\right]
(y_3)^2
+ \\
& \
\ \ \ \ \
+
\left[
-
70\,\mathcal{X}_{y^2}
\right]
y_1(y_3)^2
+
\left[
15\,\mathcal{Y}_{x^2y}
-
20\,\mathcal{X}_{x^3}
\right]
y_4
+ \\
& \
\ \ \ \ \
+
\left[
30\,\mathcal{Y}_{xy^2}
-
75\,\mathcal{X}_{x^2y}
\right]
y_1y_4
+
\left[
15\,\mathcal{Y}_{y^3}
-
90\,\mathcal{X}_{xy^2}
\right]
(y_1)^2y_4
+ \\
& \
\ \ \ \ \
+
\left[
-
35\,\mathcal{X}_{y^3}
\right]
(y_1)^3y_4
+
\left[
15\,\mathcal{Y}_{y^2}
-
90\,\mathcal{X}_{xy}
\right]
y_2y_4
+ \\
& \
\ \ \ \ \
+ 
\left[
-
105\,\mathcal{X}_{y^2}
\right]
y_1y_2y_4
+ 
\left[
-
35\,\mathcal{X}_y
\right]
y_3y_4
+
\left[
6\,\mathcal{Y}_{xy}
-
15\,\mathcal{X}_{x^2}
\right]
y_5
+ \\
& \
\ \ \ \ \
+
\left[
6\,\mathcal{Y}_{y^2}
-
36\,\mathcal{X}_{xy}
\right]
y_1y_5+
\left[
-
21\,\mathcal{X}_{y^2}
\right]
(y_1)^2y_5
+
\left[
-
21\,\mathcal{X}_y
\right]
y_2y_5
+ \\
& \
\ \ \ \ \
+
\left[
\mathcal{Y}_y
-
6\,\mathcal{X}_y
\right]
y_6
+
\left[
-
7\,\mathcal{X}_y
\right]
y_1y_6.
\endaligned\right.
\end{equation}

\subsection*{ 2.30.~Deduction of the classical Fa\`a 
di Bruno formula} Let $x,y \in \K$ and let $g = g(x)$, $f = f ( y)$ be
two $\mathcal{ C }^\infty$-smooth functions $\K \to \K$. Consider the
composition $h := f \circ g$, namely $h(x) = f (g (x))$. For $\lambda
\in \N$ with $\lambda \geq 1$, simply denote by $g_\lambda$ the
$\lambda$-th derivative $\frac{ d^\lambda g}{d x^\lambda}$ and
similarly for $h_\lambda$. Also, abbreviate $f_\lambda := \frac{
d^\lambda f}{ d y^\lambda}$.

By the classical formula for the derivative of a composite function,
we have $h_1 = f_1 \, g_1$. Further computations provide the following
list of subsequent derivatives of $h$:
\def\theequation{2.31}\begin{equation}
\left\{
\aligned
h_1 
& 
=
f_1 \, g_1, 
\\
h_2
& 
=
f_2 \, (g_1)^2
+ 
f_1 \, g_2, 
\\
h_3
&
=
f_3\,(g_1)^3+3\,f_2\,g_1\,g_2
+
f_1\,g_3,
\\
h_4
&
=
f_4\,(g_1)^4
+
6\,f_3\,(g_1)^2\,g_2
+
3\,f_2\,(g_2)^2
+
4\,f_2\,g_1\,g_3
+
f_1\,g_4,
\\
h_5
&
=
f_5\,(g_1)^5\,
+
10\,f_4\,(g_1)^3\,g_2
+
15\,f_3\,(g_1)^2\,g_3
+
10\,f_3\,g_1\,(g_2)^2
+ \\
& \
\ \ \ \ \
+
10\,f_2\,g_2\,g_3
+
5\,f_2\,g_1\,g_4
+
f_1\,g_5, 
\\
h_6
&
=
f_6\,(g_1)^6\,
+
15\,f_5\,(g_1)^4\,g_2
+
45\,f_4\,(g_1)^2\,(g_2)^2
+
15\,f_3\,(g_2)^3
+ \\
& \
\ \ \ \ \
+
20\,f_4\,(g_1)^3\,g_3
+
60\,f_3\,g_1\,g_2\,g_3
+
10\,f_2\,(g_3)^2
+
15\,f_3\,(g_1)^2\,g_4
+ \\
& \
\ \ \ \ \
+
15\,f_2\,g_2\,g_4
+
6\,f_2\,g_1\,g_5
+
f_1\,g_6.
\endaligned\right.
\end{equation}

\def\thetheorem{2.32}\begin{theorem}
For every integer $\kappa \geq 1$, the $\kappa$-th derivative of the
composite function $h = f\circ g$ may be expressed as an explicit
polynomial in the partial derivatives of $f$ and of $g$ having integer
coefficients{\rm :}
\def\theequation{2.33}\begin{equation}
\boxed{
\aligned
\frac{ d^\kappa h}{dx^\kappa}
& 
=
\sum_{d=1}^\kappa
\
\sum_{1\leq\lambda_1<\cdots<\lambda_d\leq\kappa}
\
\sum_{\mu_1\geq 1,\dots,\mu_d\geq 1}
\
\sum_{\mu_1\lambda_1+\cdots+\mu_d\lambda_d=\kappa}
\\
& \
\ \ \ \ \ 
\frac{\kappa !}{(\lambda_1!)^{\mu_1}\ \mu_1! 
\cdots
(\lambda_d!)^{\mu_d}\ \mu_d!}
\
\frac{d^{\mu_1+\cdots+\mu_d} f}{
dy^{\mu_1+\cdots+\mu_d}}
\
\left(
\frac{d^{\lambda_1}g}{dx^{\lambda_1}}
\right)^{\mu_1}
\cdots\cdots
\left(
\frac{d^{\lambda_d}g}{dx^{\lambda_d}}
\right)^{\mu_d}.
\endaligned
}
\end{equation}
\end{theorem}

This is the classical {\it Fa\`a di Bruno formula}. Interestingly, we
observe that this formula is included as subpart of the general
formula for ${\bf Y}_\kappa$, after a suitable translation. Indeed,
in the formulas for ${\bf Y}_1$, ${\bf Y}_2$, ${\bf Y}_3$, ${\bf
Y}_4$, ${\bf Y}_5$, ${\bf Y}_6$ and in the general sum for ${\bf
Y}_\kappa$, pick only the terms for which $\mu_1\lambda_1 + \cdots +
\mu_d \lambda_d = \kappa$ and drop $\mathcal{ X}$, which yields:
\def\theequation{2.34}\begin{equation}
\aligned
& \
\sum_{d=1}^\kappa
\
\sum_{1\leq\lambda_1<\cdots<\lambda_d\leq \kappa}
\
\sum_{\mu_1\geq 1,\dots,\mu_d\geq 1}
\
\sum_{\mu_1\lambda_1+\cdots+\mu_d\lambda_d=\kappa}
\\
& \
\left[
\frac{\kappa!}{
\mu_1!(\lambda_1!)^{\mu_1}
\cdots
\mu_d!(\lambda_d!)^{\mu_d}
}
\
\mathcal{Y}_{y^{\mu_1+\cdots+\mu_d}}
\right]
\left(
y_{\lambda_1}
\right)^{\mu_1}
\cdots
\left(
y_{\lambda_d}\right)^{\mu_d}.
\endaligned
\end{equation}
The similarity between the two formulas~\thetag{ 2.33} and~\thetag{
2.34} is now clearly visible.

The Fa\`a di Bruno formula may be established by means of
substitutions of power series (\cite{ f1969}, p.~222), by means of the
umbral calculus (\cite{ cs1996}), or by means of some induction
formulas, which we write for completeness. Define the differential
operators
\def\theequation{2.35}\begin{equation}
\small
\aligned
F_2
&
:=
g_2\,\frac{\partial}{\partial g_1}
+
g_1\left(
f_2\,\frac{\partial}{\partial f_1}
\right), 
\\
F_3
&
:=
g_2\,\frac{\partial}{\partial g_1}
+
g_3\,\frac{\partial}{\partial g_2}
+
g_1\left(
f_2\,\frac{\partial}{\partial f_1}
+
f_3\,\frac{\partial}{\partial f_2}
\right), 
\\
\cdots
&
\cdots\cdots
\cdots\cdots
\cdots\cdots
\cdots\cdots
\cdots\cdots
\cdots\cdots
\cdots\cdots
\\
F_\lambda
&
:=
g_2\,\frac{\partial}{\partial g_1}
+
g_3\,\frac{\partial}{\partial g_2}
+
\cdots
+
g_\lambda\,\frac{\partial}{\partial g_{\lambda-1}}
+
g_1\left(
f_2\,\frac{\partial}{\partial f_1}
+
f_3\,\frac{\partial}{\partial f_2}
+
\cdots
+
f_\lambda\,\frac{\partial}{\partial f_{\lambda-1}}
\right).
\endaligned
\end{equation}
Then we have 
\def\theequation{2.36}\begin{equation}
\aligned
h_2
&
=
F^2(h_1),
\\
h_3
&
=
F^3(h_2),
\\
\cdots
&
\cdots
\cdots
\cdots
\cdots
\\
h_\lambda
&
=
F^\lambda(h_{\lambda-1}).
\endaligned
\end{equation}

\section*{\S3.~Several independent variables and
one dependent variable}

\subsection*{3.1.~Simplified adapted notations}
As announced after the statement of Theorem~2.24, it is only after we
have treated the case of several independent variables that we will
understand perfectly the general formula~\thetag{ 2.25}, valid in the
case of one independent variable and one dependent variable. We will
discover massive formal computations, exciting our computational
intuition.

Thus, assume $n\geq 1$ and $m=1$, let $\kappa\in \N$ with $\kappa \geq
1$ and simply denote (instead of~\thetag{ 1.2}) 
the jet variables by:
\def\theequation{3.2}\begin{equation}
\left(
x^i,y,y_{i_1},y_{i_1,i_2},\dots,y_{i_1,i_2,\dots,i_\kappa}
\right).
\end{equation}
Also, instead of~\thetag{ 1.30}, denote
the $\kappa$-th prolongation of a vector field by:
\def\theequation{3.3}\begin{equation}
\left\{
\aligned
\mathcal{L}^{(\kappa)}
&
=
\sum_{i=1}^n\,\mathcal{X}^i\,\frac{\partial}{\partial x^i}
+
\mathcal{Y}\,\frac{\partial}{\partial y}
+
\sum_{i_1=1}^n\,{\bf Y}_{i_1}\,\frac{\partial}{\partial y_{i_1}}
+
\sum_{i_1,i_2=1}^n\,{\bf Y}_{i_1,i_2}\,
\frac{\partial}{\partial y_{i_1,i_2}}
+ \\
& \
\ \ \ \ \
+
\cdots
+
\sum_{i_1,i_2,\dots,i_\kappa=1}^n\,
{\bf Y}_{i_1,i_2,\dots,i_\kappa}\,
\frac{\partial}{\partial y_{i_1,i_2,\dots,i_\kappa}}.
\endaligned\right.
\end{equation}
The induction formulas are
\def\theequation{3.4}\begin{equation}
\left\{
\aligned
{\bf Y}_{i_1}
&
:=
D_{i_1}^1
\left(
\mathcal{ Y}
\right)
-
\sum_{k=1}^n\,D_{i_1}^1
\left(
\mathcal{X}^k
\right)
\, y_k, 
\\
{\bf Y}_{i_1,i_2}
&
:=
D_{i_2}^2
\left(
{\bf Y}_{i_1}
\right)
-
\sum_{k=1}^n\,D_{i_2}^1
\left(
\mathcal{X}^k
\right)
\, y_{i_1,k}, 
\\
\cdots \cdots \cdots 
&
\cdots \cdots \cdots \cdots \cdots \cdots \cdots \cdots
\cdots \cdots \cdots \cdots \cdots \cdots \cdots \cdots
\\
{\bf Y}_{i_1,i_2,\dots,i_\kappa}
&
:=
D_{i_\kappa}^\kappa
\left(
{\bf Y}_{i_1,i_2,\dots,i_{\kappa-1}}
\right)
-
\sum_{k=1}^n\,D_{i_\kappa}^1
\left(
\mathcal{X}^k
\right)
\, y_{i_1,i_2,\dots,i_{\kappa-1},k},
\endaligned\right.
\end{equation}
where the total differentiation operators $D_{ i' }^\lambda$ are
defined as in~\thetag{ 1.22}, dropping the sums $\sum_{j ' = 1}^m$ and
the indices $j'$.

\subsection*{3.5.~Two instructing explicit computations}
To begin with, let us compute ${\bf Y}_{i_1}$. With $D_{i_1}^1 =
\frac{ \partial }{\partial x^{i_1}} + y_{i_1}\, \frac{\partial
}{\partial y}$, we have:
\def\theequation{3.6}\begin{equation}
\aligned
{\bf Y}_{i_1}
&
=
D_{i_1}
\left(
\mathcal{Y}
\right)
-
\sum_{k_1=1}^n\,D_{i_1}^1
\left(
\mathcal{X}^{k_1}
\right)
y_{k_1} 
\\
& 
=
\mathcal{Y}_{x^{i_1}}
+
\mathcal{Y}_y\,y_{i_1}
-
\sum_{k_1=1}^n\,\mathcal{X}_{x^{i_1}}^{k_1}\,y_{k_1}
-
\sum_{k_1=1}^n\,\mathcal{X}_y^{k_1}\,y_{i_1}\,y_{k_1}.
\endaligned
\end{equation}
Searching for formal harmony and for coherence with the formula
$(2.6)_1$, we must include the term $\mathcal{ Y}_y\, y_{i_1}$ inside
the sum $\sum_{ k_1 =1 }^n\, \left[ \cdot \right] y_{ k_1}$. Using
the Kronecker symbol, we may write:
\def\theequation{3.7}\begin{equation}
\mathcal{Y}_y\,y_{i_1}
\equiv
\sum_{k_1=1}^n\,
\left[
\delta_{i_1}^{k_1}\,\mathcal{Y}_y
\right]
y_{k_1}.
\end{equation}
Also, we may rewrite the last term of~\thetag{ 3.6} with a double sum:
\def\theequation{3.8}\begin{equation}
-
\sum_{k_1=1}^n\,
\mathcal{X}_y^{k_1}\,y_{i_1}\,y_{k_1}
\equiv
\sum_{k_1,k_2=1}^n\,
\left[
-
\delta_{i_1}^{k_1}\,\mathcal{X}_y^{k_2}
\right]
y_{k_1}y_{k_2}.
\end{equation}
From now on and up to equation~\thetag{ 3.39}, we shall abbreviate any
sum $\sum_{k=1}^n$ from $1$ to $n$ as $\sum_k$. Putting everything
together, we get the final desired perfect expression of ${\bf
Y}_{i_1}$:
\def\theequation{3.9}\begin{equation}
{\bf Y}_{i_1}
=
\mathcal{Y}_{x^{i_1}}
+
\sum_{k_1}\,
\left[
\delta_{i_1}^{k_1}\,\mathcal{Y}_y
-
\mathcal{X}_{x^{i_1}}^{k_1}
\right]
y_{k_1}
+
\sum_{k_1,k_2}\,
\left[
-
\delta_{i_1}^{k_1}\,\mathcal{X}_y^{k_2}
\right]
y_{k_1}y_{k_2}.
\end{equation}
This completes the first explicit computation.

The second one is about ${\bf Y}_{ i_1, i_2}$. It becomes more
delicate, because several algebraic transformations must be achieved
until the final satisfying formula is obtained. Our goal is to
present each step very carefully, explaining every tiny
detail. Without such a care, it would be impossible to claim that some
of our subsequent computations, for which we will not provide the
intermediate steps, may be redone and verified. Consequently, we will
expose our rules of formal computation thoroughly.

Replacing the value of ${\bf Y}_1$ just obtained in the induction
formula $(3.4)_2$ and developing, we may conduct the very
first steps of the computation:
$$
\small
\aligned
{\bf Y}_{i_1,i_2}
&
=
D_{i_2}^2
\left(
{\bf Y}_{i_1}
\right)
-
\sum_{k_1}\,D_{i_2}^1
\left(
\mathcal{X}^{k_1}
\right)
y_{i_1,k_1} 
\\
& 
=
\left(
\frac{\partial}{\partial x^{i_2}}
+
y_{i_2}\,\frac{\partial}{\partial y}
+
\sum_{k_1}\,y_{i_2,k_1}\,
\frac{\partial}{\partial y_{k_1}}
\right)
\left(
\mathcal{Y}_{x^{i_1}}
+
\sum_{k_1}\,
\left[
\delta_{i_1}^{k_1}\,\mathcal{Y}_y
-
\mathcal{X}_{x^{i_1}}^{k_1}
\right]
y_{k_1}
+
\right.
\\
& \
\ \ \ \ \ \ \ \ \ \ \ \ \ \ \
\ \ \ \ \ \
\left.
+
\sum_{k_1,k_2}\,
\left[
-
\delta_{i_1}^{k_1}\,\mathcal{X}_y^{k_2}
\right]
y_{k_1}y_{k_2}
\right)
-
\sum_{k_1}\,
\left[
\mathcal{X}_{x^{i_2}}^{k_1}
+
y_{i_2}\,\mathcal{X}_y^{k_1}
\right]
y_{i_1,k_1}
\endaligned
$$
\def\theequation{3.10}\begin{equation}
\small
\aligned
&
=
\left(
\frac{\partial}{\partial x^{i_2}}
\right)
\left(
\mathcal{Y}_{x^{i_1}}
+
\sum_{k_1}\,
\left[
\delta_{i_1}^{k_1}\,\mathcal{Y}_y
-
\mathcal{X}_{x^{i_1}}^{k_1}
\right]
y_{k_1}
+
\sum_{k_1,k_2}\,
\left[
-
\delta_{i_1}^{k_1}\,\mathcal{X}_y^{k_2}
\right]
y_{k_1}y_{k_2}
\right)
+ \\
& \
\ \ \ \ \
+
\left(
y_{i_2}\,\frac{\partial}{\partial y}
\right)
\left(
\mathcal{Y}_{x^{i_1}}
+
\sum_{k_1}\,
\left[
\delta_{i_1}^{k_1}\,\mathcal{Y}_y
-
\mathcal{X}_{x^{i_1}}^{k_1}
\right]
y_{k_1}
+
\sum_{k_1,k_2}\,
\left[
-
\delta_{i_1}^{k_1}\,\mathcal{X}_y^{k_2}
\right]
y_{k_1}y_{k_2}
\right)
+ \\
& \
\ \ \ \ \
+
\left(
\sum_{k_1}\,y_{i_2,k_1}\,
\frac{\partial}{\partial y_{k_1}}
\right)
\left(
\mathcal{Y}_{x^{i_1}}
+
\sum_{k_1}\,
\left[
\delta_{i_1}^{k_1}\,\mathcal{Y}_y
-
\mathcal{X}_{x^{i_1}}^{k_1}
\right]
y_{k_1}
+
\sum_{k_1,k_2}\,
\left[
-
\delta_{i_1}^{k_1}\,\mathcal{X}_y^{k_2}
\right]
y_{k_1}y_{k_2}
\right) 
+ \\
& \
\ \ \ \ \
+
\sum_{k_1}\,
\left[
-
\mathcal{X}_{x^{i_2}}^{k_1}
\right]
y_{k_1,i_1}
+
\sum_{k_1}\,
\left[
-
\mathcal{X}_y^{k_1}
\right]
y_{i_2}y_{i_1,k_1}
\endaligned
\end{equation}
$$
\small
\aligned
& 
=
\mathcal{Y}_{x^{i_1}x^{i_2}}
+
\sum_{k_1}\,
\left[
\delta_{i_1}^{k_1}\,\mathcal{Y}_{x^{i_2}y}
-
\mathcal{X}_{x^{i_1}x^{i_2}}^{k_1}
\right]
y_{k_1}
+
\sum_{k_1,k_2}\,
\left[
-
\delta_{i_1}^{k_1}\,\mathcal{X}_{x^{i_2}y}^{k_2}
\right]
y_{k_1}y_{k_2}
+ \\
& \
\ \ \ \ \
+
\mathcal{Y}_{x^{i_1}y}\,y_{i_2}
+
\sum_{k_1}\,
\left[
\delta_{i_1}^{k_1}\,\mathcal{Y}_{yy}
-
\mathcal{X}_{x^{i_1}y}^{k_1}
\right]
y_{k_1}y_{i_2}
+
\sum_{k_1,k_2}\,
\left[
-
\delta_{i_1}^{k_1}\,\mathcal{X}_{yy}^{k_2}
\right]
y_{k_1}y_{k_2}y_{i_2}
+ \\
& \
\ \ \ \ \
+
\sum_{k_1}\,
\left[
\delta_{i_1}^{k_1}\,\mathcal{Y}_y
-
\mathcal{X}_{x^{i_1}}^{k_1}
\right]
y_{i_2,k_1}
+
\sum_{k_1,k_2}\,
\left[
-
\delta_{i_1}^{k_1}\,\mathcal{X}_y^{k_2}
\right]
y_{k_2}y_{i_2,k_1}
+
\sum_{k_1,k_2}\,
\left[
-
\delta_{i_1}^{k_1}\,\mathcal{X}_y^{k_2}
\right]
y_{k_1}y_{i_2,k_2}
+ \\
& \
\ \ \ \ \
+
\sum_{k_1}\,
\left[
-
\mathcal{X}_{x^{i_2}}^{k_1}
\right]
y_{k_1,i_1}
+
\sum_{k_1}
\left[
-
\mathcal{X}_y^{k_1}
\right]
y_{i_2}y_{i_1,k_1}.
\endaligned
$$
Some explanations are needed about the computation of the last two
terms of line 11, {\it i.e.} about the passage from line 7 of~\thetag{
3.10} just above to line 11. We have to compute:
\def\theequation{3.11}\begin{equation}
\footnotesize
\left(
\sum_{k_1}\,y_{i_2,k_1}\,
\frac{\partial}{\partial y_{k_1}}
\right)
\left(
\sum_{k_1,k_2}\,
\left[
-
\delta_{i_1}^{k_1}\,
\mathcal{X}_y^{k_2}
\right]
y_{k_1}y_{k_2}
\right).
\end{equation}
This term is of the form
\def\theequation{3.12}\begin{equation}
\footnotesize
\left(
\sum_{k_1}\,A_{k_1}\,
\frac{\partial}{\partial y_{k_1}}
\right)
\left(
\sum_{k_1,k_2}\,
\left[
B_{k_1,k_2}
\right]
y_{k_1}y_{k_2}
\right),
\end{equation}
where the terms $B_{k_1,k_2}$ are independent of the pure first jet
variables $y_{x^k}$.
By the rule of Leibniz for the differentiation 
of a product, we may write
\def\theequation{3.13}\begin{equation}
\footnotesize
\aligned
&
\left(
\sum_{k_1}\,A_{k_1}\,
\frac{\partial}{\partial y_{k_1}}
\right)
\left(
\sum_{k_1,k_2}\,
\left[
B_{k_1,k_2}
\right]
y_{k_1}y_{k_2}
\right)
= 
\\
&
=
\sum_{k_1,k_2}\,
\left[
B_{k_1,k_2}
\right]
y_{k_2}
\left(
\sum_{k_1'}\,A_{k_1'}\,
\frac{\partial}{\partial y_{k_1'}}
(y_{k_1})
\right)
+
\sum_{k_1,k_2}\,
\left[
B_{k_1,k_2}
\right]
y_{k_1}
\left(
\sum_{k_2'}\,A_{k_2'}\,
\frac{\partial}{\partial y_{k_2'}}
(y_{k_2})
\right)
\\
& 
=
\sum_{k_1,k_2}\,
\left[
B_{k_1,k_2}
\right]
y_{k_2}\,A_{k_1}
+
\sum_{k_1,k_2}\,
\left[
B_{k_1,k_2}
\right]
y_{k_1}\,A_{k_2}.
\endaligned
\end{equation}
This is how we have written line 11 of~\thetag{ 3.10}.

Next, the first term $\mathcal{ Y}_{ x^{i_1}y} \, y_{ i_2}$ in line 10
of~\thetag{ 3.10} is not in a suitable shape. For reasons of harmony
and coherence, we must insert it inside a sum of the form
$\sum_{k_1}\, \left[ \cdot \right] y_{k_1}$. Hence, using the
Kronecker symbol, we transform:
\def\theequation{3.14}\begin{equation}
\mathcal{ Y}_{x^{i_1}y} \, y_{i_2}
\equiv
\sum_{k_1}\,
\left[
\delta_{i_2}^{k_1}\,
\mathcal{Y}_{x^{i_1}y}
\right]
y_{k_1}.
\end{equation}
Also, we must ``summify'' the seven other terms, remaining in lines 10,
11 and 12 of~\thetag{ 3.10}. Sometimes, we use the symmetry $y_{i_2,
k_1} \equiv y_{ k_1, i_2}$ without mention. Similarly, we get:
$$
\aligned
\sum_{k_1}\,
\left[
\delta_{i_1}^{k_1}\,\mathcal{Y}_{yy}
-
\mathcal{X}_{x^{i_1y}}^{k_1}
\right]
y_{k_1}y_{i_2}
&
\equiv
\sum_{k_1,k_2}\,
\left[
\delta_{i_1}^{k_1}\,\delta_{i_2}^{k_2}\,
\mathcal{Y}_{yy}
-
\delta_{i_2}^{k_2}\,\mathcal{X}_{x^{i_1}y}^{k_1}
\right]
y_{k_1}y_{k_2}, \\
\sum_{k_1,k_2}\,
\left[
-
\delta_{i_1}^{k_1}\,
\mathcal{X}_{yy}^{k_2}
\right]
y_{k_1}y_{k_2}y_{i_2}
&
\equiv
\sum_{k_1,k_2,k_3}\,
\left[
-
\delta_{i_1}^{k_1}\,
\delta_{i_2}^{k_3}
\mathcal{X}_{yy}^{k_2}
\right]
y_{k_1}y_{k_2}y_{k_3}, \\
\sum_{k_1}
\left[
\delta_{i_1}^{k_1}\,
\mathcal{Y}_y
-
\mathcal{X}_{x^{i_1}}^{k_1}
\right]
y_{k_1,i_2}
&
\equiv
\sum_{k_1,k_2}
\left[
\delta_{i_1}^{k_1}\,
\delta_{i_2}^{k_2}\,
\mathcal{Y}_y
-
\delta_{i_2}^{k_2}\,
\mathcal{X}_{x^{i_1}}^{k_1}
\right]
y_{k_1,k_2},
\endaligned
$$
\def\theequation{3.15}\begin{equation}
\aligned
\sum_{k_1,k_2}\,
\left[
-
\delta_{i_1}^{k_1}\,\mathcal{X}_y^{k_2}
\right]
y_{k_2}y_{k_1,i_2}
&
=
\sum_{k_1,k_2}\,
\left[
-
\delta_{i_1}^{k_2}\,\mathcal{X}_y^{k_1}
\right]
y_{k_1}y_{k_2,i_2}
\\
& 
\equiv
\sum_{k_1,k_2,k_3}\,
\left[
-
\delta_{i_1}^{k_2}\,\delta_{i_2}^{k_3}\,
\mathcal{X}_y^{k_1}
\right]
y_{k_1}y_{k_2,k_3},
\endaligned
\end{equation}
$$
\aligned
\sum_{k_1,k_2}\,
\left[
-
\delta_{i_1}^{k_1}\,\mathcal{X}_y^{k_2}
\right]
y_{k_1}y_{k_2,i_2}
&
\equiv
\sum_{k_1,k_2,k_3}\,
\left[
-
\delta_{i_1}^{k_1}\,
\delta_{i_2}^{k_3}\,
\mathcal{X}_y^{k_2}
\right]
y_{k_1}y_{k_2,k_3}, 
\\
\sum_{k_1}\,
\left[
-
\mathcal{X}_{x^{i_2}}^{k_1}
\right]
y_{k_1,i_1}
&
\equiv
\sum_{k_1,k_2}\,
\left[
-
\delta_{i_1}^{k_2}
\mathcal{X}_{x^{i_2}}^{k_1}
\right]
y_{k_1,k_2}, 
\\
\sum_{k_1}\,
\left[
-
\mathcal{X}_y^{k_1}
\right]
y_{i_2}y_{k_1,i_1}
&
=
\sum_{k_2}\,
\left[
-
\mathcal{X}_y^{k_2}
\right]
y_{i_2}y_{k_2,i_1}
\\
& 
\equiv
\sum_{k_1,k_2,k_3}\,
\left[
-
\delta_{i_2}^{k_1}\,
\delta_{i_1}^{k_3}\,
\mathcal{X}_y^{k_2}
\right]
y_{k_1}y_{k_2,k_3}.
\endaligned
$$
In the sequel, for products of
Kronecker symbols, it will be convenient
to adopt the following self-evident contracted
notation:
\def\theequation{3.16}\begin{equation}
\delta_{i_1}^{k_1}\,
\delta_{i_2}^{k_2}
\equiv
\delta_{i_1,\,i_2}^{k_1,k_2}; 
\ \ \ \ \ \
{\rm generally:} 
\ \ \
\delta_{i_1}^{k_1}\,
\delta_{i_2}^{k_2}
\cdots
\delta_{i_\lambda}^{k_\lambda}
\equiv
\delta_{i_1,\,i_2,\,\cdots\,,i_\lambda}^{
k_1,k_2,\cdots,k_\lambda}.
\end{equation}
Re-inserting plainly these eight summified terms~\thetag{ 3.14},
\thetag{ 3.15} in the last expression~\thetag{ 3.10} of ${\bf Y}_{i_1,
i_2}$ (lines 10, 11 and 12), we get:
\def\theequation{3.17}\begin{equation}
\small
\aligned
{\bf Y}_{i_1,i_2}
&
=
\underline{
\mathcal{Y}_{x^{i_1}x^{i_2}} 
}_{ \fbox{\tiny 1}}
+
\underline{
\sum_{k_1}\,
\left[
\delta_{i_1}^{k_1}\,
\mathcal{Y}_{x^{i_2}y}
-
\mathcal{X}_{x^{i_1}x^{i_2}}^{k_1}
\right]
y_{k_1}
}_{ \fbox{\tiny 2}}
+
\underline{
\sum_{k_1,k_2}\,
\left[
-
\delta_{i_1}^{k_1}\,
\mathcal{X}_{x^{i_2}y}^{k_2}
\right]
y_{k_1}y_{k_2}
}_{ \fbox{\tiny 3}}
+ \\
& \
\ \ \ \ \
+
\underline{
\sum_{k_1}\,
\left[
\delta_{i_2}^{k_1}\,
\mathcal{Y}_{x^{i_1}y}
\right]
y_{k_1}
}_{ \fbox{\tiny 2}}
+ 
\underline{
\sum_{k_1,k_2}\,
\left[
\delta_{i_1,\,i_2}^{k_1,k_2}\,
\mathcal{Y}_{yy}
-
\delta_{i_2}^{k_2}\,
\mathcal{X}_{x^{i_1}y}^{k_1}
\right]
y_{k_1}y_{k_2}
}_{ \fbox{\tiny 3}}
+ \\
& \
\ \ \ \ \
+
\underline{
\sum_{k_1,k_2,k_3}\,
\left[
-
\delta_{i_1,\,i_2}^{k_1,k_3}\,
\mathcal{X}_{yy}^{k_2}
\right]
y_{k_1}y_{k_2}y_{k_3}
}_{ \fbox{\tiny 4}}
+
\underline{
\sum_{k_1,k_2}\,
\left[
\delta_{i_1,\,i_2}^{k_1,k_2}\,
\mathcal{Y}_y
-
\delta_{i_2}^{k_2}\,
\mathcal{X}_{x^{i_1}}^{k_1}
\right]
y_{k_1,k_2}
}_{ \fbox{\tiny 5}}
+ \\
& \
\ \ \ \ \
+
\underline{
\sum_{k_1,k_2,k_3}\,
\left[
-
\delta_{i_1,\,i_2}^{k_2,k_3}\,
\mathcal{X}_y^{k_1}
\right]
y_{k_1}y_{k_2,k_3}
}_{ \fbox{\tiny 6}}
+
\underline{
\sum_{k_1,k_2,k_3}\,
\left[
-
\delta_{i_1,\,i_2}^{k_1,k_3}\,
\mathcal{X}_y^{k_2}
\right]
y_{k_1}y_{k_2,k_3}
}_{ \fbox{\tiny 6}}
+ \\
& \
\ \ \ \ \
+
\underline{
\sum_{k_1,k_2}\,
\left[
-
\delta_{i_1}^{k_2}\,
\mathcal{X}_{x^{i_2}}^{k_1}
\right]
y_{k_1,k_2}
}_{ \fbox{\tiny 5}}
+
\underline{
\sum_{k_1,k_2,k_3}\,
\left[
-
\delta_{i_2,\,i_1}^{k_1,k_3}\,
\mathcal{X}_y^{k_2}
\right]
y_{k_1}y_{k_2,k_3}
}_{ \fbox{\tiny 6}}.
\endaligned
\end{equation}
Next, we gather the underlined terms, ordering them according to their
number. This yields 6 collections of
sums of monomials in the 
pure jet variables:
\def\theequation{3.18}\begin{equation}
\small
\aligned
{\bf Y}_{i_1,i_2}
&
=
\mathcal{Y}_{x^{i_1}x^{i_2}} 
+
\sum_{k_1}\,
\left[
\delta_{i_1}^{k_1}\,\mathcal{Y}_{x^{i_2}y}
+
\delta_{i_2}^{k_1}\,\mathcal{Y}_{x^{i_1}y}
-
\mathcal{X}_{x^{i_1}x^{i_2}}^{k_1}
\right]
y_{k_1}
+ \\
& \
\ \ \ \ \
+
\sum_{k_1,k_2}\,
\left[
\delta_{i_1,\,i_2}^{k_1,k_2}\,
\mathcal{Y}_{yy}
-
\delta_{i_1}^{k_1}\,
\mathcal{X}_{x^{i_2}y}^{k_2}
-
\delta_{i_2}^{k_2}\,
\mathcal{X}_{x^{i_1}y}^{k_1}
\right]
y_{k_1}y_{k_2}
+ \\
& \
\ \ \ \ \
+
\sum_{k_1,k_2,k_3}\,
\left[
-
\delta_{i_1,\,i_2}^{k_1,k_3}\,
\mathcal{X}_{yy}^{k_2}
\right]
y_{k_1}y_{k_2}y_{k_3}
+ \\
& \
\ \ \ \ \
+
\sum_{k_1,k_2}\,
\left[
\delta_{i_1,\,i_2}^{k_1,k_2}\,
\mathcal{Y}_y
-
\delta_{i_2}^{k_2}\,\mathcal{X}_{x^{i_1}}^{k_1}
-
\delta_{i_1}^{k_2}\,\mathcal{X}_{x^{i_2}}^{k_1}
\right]
y_{k_1,k_2}
+ \\
& \
\ \ \ \ \
+
\sum_{k_1,k_2,k_3}\,
\left[
-
\delta_{i_1,\,i_2}^{k_2,k_3}\,
\mathcal{X}_y^{k_1}
-
\delta_{i_1,\,i_2}^{k_1,k_3}\,
\mathcal{X}_y^{k_2}
-
\delta_{i_2,\,i_1}^{k_1,k_3}\,
\mathcal{X}_y^{k_2}
\right]
y_{k_1}y_{k_2,k_3}.
\endaligned
\end{equation}
To attain the real perfect harmony, this last expression has still to
be worked out a little bit.

\def\thelemma{3.19}\begin{lemma}
The final expression of ${\bf Y}_{i_1,i_2}$ is
as follows{\rm :}
\def\theequation{3.20}\begin{equation}
\left\{
\aligned
{\bf Y}_{i_1,i_2}
&
=
\mathcal{Y}_{x^{i_1}x^{i_2}} 
+
\sum_{k_1}\,
\left[
\delta_{i_1}^{k_1}\,\mathcal{Y}_{x^{i_2}y}
+
\delta_{i_2}^{k_1}\,\mathcal{Y}_{x^{i_1}y}
-
\mathcal{X}_{x^{i_1}x^{i_2}}^{k_1}
\right]
y_{k_1}
+ \\
& \
\ \ \ \ \
+
\sum_{k_1,k_2}\,
\left[
\delta_{i_1,\,i_2}^{k_1,k_2}\,
\mathcal{Y}_{yy}
-
\delta_{i_1}^{k_1}\,
\mathcal{X}_{x^{i_2}y}^{k_2}
-
\delta_{i_2}^{k_1}\,
\mathcal{X}_{x^{i_1}y}^{k_2}
\right]
y_{k_1}y_{k_2}
+ \\
& \
\ \ \ \ \
+
\sum_{k_1,k_2,k_3}\,
\left[
-
\delta_{i_1,\,i_2}^{k_1,k_2}\,
\mathcal{X}_{yy}^{k_3}
\right]
y_{k_1}y_{k_2}y_{k_3}
+ \\
& \
\ \ \ \ \
+
\sum_{k_1,k_2}\,
\left[
\delta_{i_1,\,i_2}^{k_1,k_2}\,
\mathcal{Y}_y
-
\delta_{i_1}^{k_1}\,\mathcal{X}_{x^{i_2}}^{k_2}
-
\delta_{i_2}^{k_1}\,\mathcal{X}_{x^{i_1}}^{k_2}
\right]
y_{k_1,k_2}
+ \\
& \
\ \ \ \ \
+
\sum_{k_1,k_2,k_3}\,
\left[
-
\delta_{i_1,\,i_2}^{k_1,k_2}\,
\mathcal{X}_y^{k_3}
-
\delta_{i_1,\,i_2}^{k_3,k_1}\,
\mathcal{X}_y^{k_2}
-
\delta_{i_1,\,i_2}^{k_2,k_3}\,
\mathcal{X}_y^{k_1}
\right]
y_{k_1}y_{k_2,k_3}.
\endaligned\right.
\end{equation}
\end{lemma}

\proof
As promised, we explain every tiny detail.

The first lines of~\thetag{ 3.18} and of~\thetag{ 3.20} are exactly
the same. For the transformations of terms in the second, in the third
and in the fourth lines, we use the following device. Let $\Upsilon_{
k_1, k_2}$ be an indexed quantity which is symmetric: $\Upsilon_{ k_1,
k_2} =\Upsilon_{ k_2, k_1}$. Let $A_{k_1, k_2}$ be an arbitrary
indexed quantity. Then obviously:
\def\theequation{3.21}\begin{equation}
\sum_{k_1,k_2}\,
A_{k_1, k_2}\,
\Upsilon_{k_1, k_2}
=
\sum_{k_1,k_2}\,
A_{k_2,k_1}\,
\Upsilon_{k_1,k_2}.
\end{equation}
Similar relations hold with a quantity $\Upsilon_{i_1, i_2, \dots,
i_\lambda}$ which is symmetric with respect to its $\lambda$ indices.
Consequently, in the second, in the third and in the fourth lines
of~\thetag{ 3.18}, we may permute freely certain indices in some of
the terms inside the brackets. This yields
the passage from lines 2, 3 and 4 of~\thetag{ 3.18}
to lines 2, 3 and 4 of~\thetag{ 3.20}.

It remains to explain how we pass from the fifth (last) line
of~\thetag{ 3.18} to the fifth (last) line of~\thetag{ 3.20}. The
bracket in the fifth line of~\thetag{ 3.18} contains three terms:
$\left[ -T_1 -T_2 -T_3 \right]$. The term $T_3$ involves the product
$\delta_{i_2,\,i_1}^{k_1,k_3}$, which we rewrite as
$\delta_{i_1,\,i_2}^{k_3,k_1}$, in order that $i_1$ appears before
$i_2$. Then, we rewrite the three terms in the new order $\left[ -T_2
-T_3 -T_1 \right]$, which yields:
\def\theequation{3.22}\begin{equation}
\sum_{k_1,k_2,k_3}\,
\left[
-
\delta_{i_1,\,i_2}^{k_1,k_3}\,
\mathcal{X}_y^{k_2}
-
\delta_{i_1,\,i_2}^{k_3,k_1}\,
\mathcal{X}_y^{k_2}
-
\delta_{i_1,\,i_2}^{k_2,k_3}\,
\mathcal{X}_y^{k_1}
\right]
y_{k_1}y_{k_2,k_3}.
\end{equation}
It remains to observe that we can permute $k_2$ and $k_3$ in the first
term $-T_2$, which yields the last
line of~\thetag{ 3.20}. The detailed proof is complete.
\endproof

\subsection*{3.23.~Final perfect expression of ${\bf Y}_{i_1, i_2,
i_3}$} Thanks to similar (longer) computations, we have obtained an
expression of ${\bf Y}_{i_1, i_2, i_3}$ which we consider to be in
final harmonious shape. Without copying the intermediate steps, let
us write down the result. The comments which 
are necessary to read it and to interpret it 
start just below.
$$
\small
\aligned
{\bf Y}_{i_1,i_2,i_3}
&
=
\mathcal{Y}_{x^{i_1}x^{i_2}x^{i_3}}
+
\sum_{k_1}\,
\left[
\delta_{i_1}^{k_1}\,\mathcal{Y}_{x^{i_2}x^{i_3}y}
+
\delta_{i_2}^{k_1}\,\mathcal{Y}_{x^{i_1}x^{i_3}y}
+
\delta_{i_3}^{k_1}\,\mathcal{Y}_{x^{i_1}x^{i_2}y}
-
\mathcal{X}_{x^{i_1}x^{i_2}x^{i_3}}^{k_1}
\right]
y_{k_1}
+ \\
& \
\ \ \ \ \
+
\sum_{k_1,k_2}\,
\left[
\delta_{i_1,\,i_2}^{k_1,k_2}\,\mathcal{Y}_{x^{i_3}y^2}
+
\delta_{i_3,\,i_1}^{k_1,k_2}\,\mathcal{Y}_{x^{i_2}y^2}
+
\delta_{i_2,\,i_3}^{k_1,k_2}\,\mathcal{Y}_{x^{i_1}y^2}
-
\right. \\
& \
\ \ \ \ \ \ \ \ \ \ \ \ \ \ \ \ \ \ \ \ 
\ \ \ \ \ \ \ \ \ \ \
\left.
-
\delta_{i_1}^{k_1}\,\mathcal{X}_{x^{i_2}x^{i_3}y}^{k_2}
-
\delta_{i_2}^{k_1}\,\mathcal{X}_{x^{i_1}x^{i_3}y}^{k_2}
-
\delta_{i_3}^{k_1}\,\mathcal{X}_{x^{i_1}x^{i_2}y}^{k_2}
\right]
y_{k_1}y_{k_2}
+ \\
& \
\ \ \ \ \
+
\sum_{k_1,k_2,k_3}\,
\left[
\delta_{i_1,\,i_2,\,i_3}^{k_1,k_2,k_3}\,\mathcal{Y}_{y^3}
-
\delta_{i_1,\,i_2}^{k_1,k_2}\,\mathcal{X}_{x^{i_3}y^2}^{k_3}
-
\delta_{i_1,\,i_3}^{k_1,k_2}\,\mathcal{X}_{x^{i_2}y^2}^{k_3}
-
\delta_{i_2,\,i_3}^{k_1,k_2}\,\mathcal{X}_{x^{i_1}y^2}^{k_3}
\right]
y_{k_1}y_{k_2}y_{k_3}
+ \\
& \
\ \ \ \ \
+
\sum_{k_1,k_2,k_3,k_4}\,
\left[
-\delta_{i_1,\,i_2,\,i_3}^{k_1,k_2,k_3}\,\mathcal{X}_{y^3}^{k_4}
\right]
y_{k_1}y_{k_2}y_{k_3}y_{k_4}
+ \\
\endaligned
$$
\def\theequation{3.24}\begin{equation}
\small
\aligned
& \
\ \ \ \ \
+
\sum_{k_1,k_2}\,
\left[
\delta_{i_1,\,i_2}^{k_1,k_2}\,\mathcal{Y}_{x^{i_3}y}
+
\delta_{i_3,\,i_1}^{k_1,k_2}\,\mathcal{Y}_{x^{i_2}y}
+
\delta_{i_2,\,i_3}^{k_1,k_2}\,\mathcal{Y}_{x^{i_1}y}
-
\right. 
\\
& \
\ \ \ \ \ \ \ \ \ \ \ \ \ \ \ \ \ \ \ \ \ \ 
\left.
-
\delta_{i_1}^{k_1}\,\mathcal{X}_{x^{i_2}x^{i_3}}^{k_2}
-
\delta_{i_2}^{k_1}\,\mathcal{X}_{x^{i_1}x^{i_3}}^{k_2}
-
\delta_{i_3}^{k_1}\,\mathcal{X}_{x^{i_1}x^{i_2}}^{k_2}
\right]
y_{k_1,k_2}
+ \\
& \
\ \ \ \ \
+
\sum_{k_1,k_2,k_3}\,
\left[
\delta_{i_1,\,i_2,\,i_3}^{k_1,k_2,k_3}\,\mathcal{Y}_{y^2}
+
\delta_{i_1,\,i_2,\,i_3}^{k_3,k_1,k_2}\,\mathcal{Y}_{y^2}
+
\delta_{i_1,\,i_2,\,i_3}^{k_2,k_3,k_1}\,\mathcal{Y}_{y^2}
-
\right.
\\
& \
\ \ \ \ \ \ \ \ \ \ \ \ \ \ \ \ \ \ \ \ \ \ \ \ \
\left.
-
\delta_{i_1,\,i_2}^{k_1,k_2}\,\mathcal{X}_{x^{i_3}y}^{k_3}
-
\delta_{i_1,\,i_2}^{k_3,k_1}\,\mathcal{X}_{x^{i_3}y}^{k_2}
-
\delta_{i_1,\,i_2}^{k_2,k_3}\,\mathcal{X}_{x^{i_3}y}^{k_1}
-
\right.
\\
& \
\ \ \ \ \ \ \ \ \ \ \ \ \ \ \ \ \ \ \ \ \ \ \ \ \
\left.
-
\delta_{i_1,\,i_3}^{k_1,k_2}\,\mathcal{X}_{x^{i_2}y}^{k_3}
-
\delta_{i_1,\,i_3}^{k_3,k_1}\,\mathcal{X}_{x^{i_2}y}^{k_2}
-
\delta_{i_1,\,i_3}^{k_2,k_3}\,\mathcal{X}_{x^{i_2}y}^{k_1}
-
\right.
\\
& \
\ \ \ \ \ \ \ \ \ \ \ \ \ \ \ \ \ \ \ \ \ \ \ \ \
\left.
-
\delta_{i_2,\,i_3}^{k_1,k_2}\,\mathcal{X}_{x^{i_1}y}^{k_3}
-
\delta_{i_2,\,i_3}^{k_3,k_1}\,\mathcal{X}_{x^{i_1}y}^{k_2}
-
\delta_{i_2,\,i_3}^{k_2,k_3}\,\mathcal{X}_{x^{i_1}y}^{k_1}
\right]
y_{k_1}y_{k_2,k_3}
+
\\
\endaligned
\end{equation}
$$
\small
\aligned
& \
\ \ \ \ \
+ 
\sum_{k_1,k_2,k_3,k_4}\,
\left[
-
\delta_{i_1,\,i_2,\,i_3}^{k_1,k_2,k_3}\,\mathcal{X}_{y^2}^{k_4}
-
\delta_{i_1,\,i_2,\,i_3}^{k_2,k_3,k_1}\,\mathcal{X}_{y^2}^{k_4}
-
\delta_{i_1,\,i_2,\,i_3}^{k_3,k_2,k_1}\,\mathcal{X}_{y^2}^{k_4}
- 
\right.
\\
& \
\ \ \ \ \ \ \ \ \ \ \ \ \ \ \ \ \ \ \ \ \ \ \ \ \
\left.
-
\delta_{i_1,\,i_2,\,i_3}^{k_3,k_4,k_1}\,\mathcal{X}_{y^2}^{k_2}
-
\delta_{i_1,\,i_2,\,i_3}^{k_3,k_1,k_4}\,\mathcal{X}_{y^2}^{k_2}
-
\delta_{i_1,\,i_2,\,i_3}^{k_1,k_3,k_4}\,\mathcal{X}_{y^2}^{k_2}
\right]
y_{k_1}y_{k_2}y_{k_3,k_4}
+ \\
& \
\ \ \ \ \
+
\sum_{k_1,k_2,k_3,k_4}\,
\left[
-
\delta_{i_1,\,i_2,\,i_3}^{k_1,k_2,k_3}\,\mathcal{X}_y^{k_4}
-
\delta_{i_1,\,i_2,\,i_3}^{k_2,k_3,k_1}\,\mathcal{X}_y^{k_4}
-
\delta_{i_1,\,i_2,\,i_3}^{k_3,k_1,k_2}\,\mathcal{X}_y^{k_4}
\right]
y_{k_1,k_2}y_{k_3,k_4}
+ \\
\endaligned
$$
$$
\small
\aligned
& \
\ \ \ \ \
+
\sum_{k_1,k_2,k_3}\,
\left[
\delta_{i_1,\,i_2,\,i_3}^{k_1,k_2,k_3}\,\mathcal{Y}_y
-
\delta_{i_1,\,i_2}^{k_1,k_2}\,\mathcal{X}_{x^{i_3}}^{k_3}
-
\delta_{i_1,\,i_3}^{k_1,k_2}\,\mathcal{X}_{x^{i_2}}^{k_3}
-
\delta_{i_2,\,i_3}^{k_1,k_2}\,\mathcal{X}_{x^{i_1}}^{k_3}
\right]
y_{k_1,k_2,k_3}
+ \\
& \
\ \ \ \ \
+
\sum_{k_1,k_2,k_3,k_4}\,
\left[
-
\delta_{i_1,\,i_2,\,i_3}^{k_1,k_2,k_3}\,\mathcal{X}_y^{k_4}
-
\delta_{i_1,\,i_2,\,i_3}^{k_4,k_1,k_2}\,\mathcal{X}_y^{k_3}
-
\delta_{i_1,\,i_2,\,i_3}^{k_3,k_4,k_1}\,\mathcal{X}_y^{k_2}
-
\delta_{i_1,\,i_2,\,i_3}^{k_2,k_3,k_4}\,\mathcal{X}_y^{k_1}
\right]
y_{k_1}y_{k_2,k_3,k_4}.
\endaligned
$$

\subsection*{3.25.~Comments, analysis and induction}
First of all, by comparing this expression of ${\bf Y}_{i_1, i_2,
i_3}$ with the expression~\thetag{ 2.8} of ${\bf Y}_3$, we easily
guess a part of the (inductional) dictionary beween the cases $n=1$
and the case $n \geq 1$. For instance, the three monomials $[\cdot ]
(y_1)^3$, $[\cdot]\, y_1 y_2$ and $[\cdot]\, (y_1 )^2\, y_2$ in ${\bf
Y }_3$ are replaced in ${\bf Y}_{ i_1, i_2, i_3}$ by the following
three sums:
\def\theequation{3.26}\begin{equation}
\small
\sum_{k_1,k_2,k_3}\,
\left[
\cdot
\right]
y_{k_1}y_{k_2}y_{k_3}, 
\ \ \ \ \ \ \ \ \ \ 
\sum_{k_1,k_2,k_3}\,
\left[
\cdot
\right]
y_{k_1}y_{k_2,k_3}, 
\ \ \ \ \ 
{\rm and}
\ \ \ \ \
\sum_{k_1,k_2,k_3,k_4}\,
\left[
\cdot
\right]
y_{k_1}y_{k_2}y_{k_3,k_4}. 
\end{equation}
Similar formal correspondences may be observed for all the monomials
of ${\bf Y}_1$, ${\bf Y}_{i_1}$, of ${\bf Y}_2$, ${\bf Y}_{i_1,i_2}$
and of ${\bf Y}_3$, ${\bf Y}_{i_1,i_2,i_3}$. Generally and
inductively speaking, the monomial
\def\theequation{3.27}\begin{equation}
\left[
\cdot
\right]
\left(
y_{\lambda_1}
\right)^{\mu_1}
\cdots
\left(
y_{\lambda_d}
\right)^{\mu_d}
\end{equation}
appearing in the expression~\thetag{ 2.25} of ${\bf Y }_\kappa$ should
be replaced by a certain multiple sum generalizing~\thetag{ 3.26}.
However, it is necessary to think, to pause and to search for an
appropriate formalism before writing down the desired multiple sum.

The jet variable $y_{ \lambda_1}$ should be replaced by a jet variable
corresponding to a $\lambda_1$-th partial derivative, say $y_{ k_1,
\dots, k_{ \lambda_1}}$, where $k_1, \dots,k_{ \lambda_1}= 1, \dots,
n$. For the moment, to simplify the discussion, we leave out the
presence of a sum of the form $\sum_{ k_1, \dots, k_{ \lambda_1}}$.
The $\mu_1$-th power $\left( y_{ \lambda_1 } \right)^{ \mu_1}$ should
be replaced {\it not}\, by $\left( y_{k_1, \dots, k_{ \lambda_1}}
\right)^{\mu_1}$, but by a product of $\mu_1$ different jet variables
$y_{k_1, \dots, k_{ \lambda_1}}$ of length $\lambda$, {\it with all
indices $k_\alpha = 1, \dots, n$ being distinct}. This rule may be
confirmed by inspecting the expressions of ${\bf Y}_{i_1}$, of ${\bf
Y}_{ i_1, i_2}$ and of ${\bf Y}_{i_1, i_2, i_3}$. So $y_{ k_1, \dots,
k_{ \lambda_1}}$ should be developed as a product of the form
\def\theequation{3.28}\begin{equation}
y_{k_1,\dots,k_{\lambda_1}}\,
y_{k_{\lambda_1+1},\dots,k_{2\lambda_1}}
\cdots \
y_{k_{(\mu_1-1)\lambda_1+1},\dots,k_{\mu_1\lambda_1}},
\end{equation}
where 
\def\theequation{3.29}\begin{equation}
k_1,\dots,k_{\lambda_1},\dots,k_{\mu_1\lambda_1}
=
1,\dots,n.
\end{equation}
Consider now the product $\left( y_{\lambda_1 } \right)^{ \mu_1}\left(
y_{\lambda_2 } \right)^{ \mu_2}$. How should it develope in the case
of several independent variables? For instance, in the expression of
${\bf Y}_{i_1,i_2,i_3}$, we have developed the product $(y_1)^2\,y_2$
as $y_{k_1} y_{k_2} y_{k_3,k_4}$. Thus, a reasonable proposal of
formalism would be that the product $\left( y_{\lambda_1 } \right)^{
\mu_1}\left( y_{\lambda_2 } \right)^{ \mu_2}$ should be developed as a
product of the form
\def\theequation{3.30}\begin{equation}
\aligned
&
y_{k_1,\dots,k_{\lambda_1}}\,
y_{k_{\lambda_1+1},\dots,k_{2\lambda_1}}
\cdots \
y_{k_{(\mu_1-1)\lambda_1+1},\dots,k_{\mu_1\lambda_1}} \\
&
y_{k_{\mu_1\lambda_1+1},\dots,k_{\mu_1\lambda_1+\lambda_2}}
\cdots \
y_{k_{\mu_1\lambda_1+(\mu_2-1)\lambda_2+1},\dots,
k_{\mu_1\lambda_1+\mu_2\lambda_2}},
\endaligned
\end{equation}
where
\def\theequation{3.31}\begin{equation}
k_1,\dots, k_{\lambda_1},\dots,k_{\mu_1\lambda_1},\dots,
k_{\mu_1\lambda_1+\mu_2\lambda_2}
=
1,\dots,n.
\end{equation}
However, when trying to write down the development of the general
monomial $\left( y_{\lambda_1 } \right)^{ \mu_1}\left( y_{\lambda_2 }
\right)^{ \mu_2} \cdots \left( y_{\lambda_d } \right)^{ \mu_d}$, we
would obtain the complicated product
\def\theequation{3.32}\begin{equation}
\small
\aligned
&
y_{k_1,\dots,k_{\lambda_1}}\,
y_{k_{\lambda_1+1},\dots,k_{2\lambda_1}}
\cdots \
y_{k_{(\mu_1-1)\lambda_1+1},\dots,k_{\mu_1\lambda_1}} 
\\
&
y_{k_{\mu_1\lambda_1+1},\dots,k_{\mu_1\lambda_1+\lambda_2}}
\dots
y_{k_{\mu_1\lambda_1+(\mu_2-1)\lambda_2+1},\dots,
k_{\mu_1\lambda_1+\mu_2\lambda_2}}
\\
&
\dots\dots\dots\dots\dots\dots
\dots\dots\dots\dots\dots\dots
\dots\dots\dots\dots\dots\dots
\\
&
y_{k_{\mu_1\lambda_1+\cdots+\mu_{d-1}\lambda_{d-1}+1},
\dots,k_{\mu_1\lambda_1+\cdots+\mu_{d-1}\lambda_{d-1}+\lambda_d}}
\cdots 
\\
& \
\ \ \ \ \ \ \ \ \ \ \ \ \ \ \ 
\cdots \
y_{k_{\mu_1\lambda_1+\cdots+\mu_{d-1}\lambda_{d-1}+
(\mu_d-1)\lambda_d+1},\dots,
k_{\mu_1\lambda_1+\cdots+\mu_d\lambda_d}}.
\endaligned
\end{equation} 
Essentially, this product is still readable. However, in it, some of
the integers $k_\alpha$ have a too long index $\alpha$, often
involving a sum. Such a length of $\alpha$ would be very inconvenient
in writing down and in reading the general Kronecker symbols $\delta_{
i_1, \dots \dots, i_\lambda }^{ k_{ \alpha_1}, \dots, k_{
\alpha_\lambda}}$ which should appear in the final expression of ${\bf
Y}_{ i_1, \dots, i_\kappa}$. One should read in advance Theorem~3.73
below to observe the presence of such multiple Kronecker symbols.
{\sf Consequently, for $\alpha = 1, \dots, \mu_1 \lambda_1, \dots,
\mu_1 \lambda_1 + \cdots + \mu_d \lambda_d$, we have to denote the
indices $k_\alpha$ differently}.

\def\thenotationalconvention{3.33}\begin{notationalconvention}
We denote $d$ collection of $\mu_d$ groups of $\lambda_d$ {\rm (a
priori distinct)} integers $k_\alpha = 1, \dots, n$ by
\def\theequation{3.34}\begin{equation}
\aligned
&
\underbrace{
\underbrace{
k_{1:1:1},\dots,k_{1:1:\lambda_1}}_{\lambda_1},
\dots,
\underbrace{
k_{1:\mu_1:1},\dots,k_{1:\mu_1:\lambda_1}}_{\lambda_1}}_{
\mu_1}, 
\\
&
\underbrace{
\underbrace{
k_{2:1:1},\dots,k_{2:1:\lambda_2}}_{\lambda_2},
\dots,
\underbrace{
k_{2:\mu_2:1},\dots,k_{2:\mu_2:\lambda_2}}_{\lambda_2}}_{
\mu_2},
\\
&
\text{\bf \ \
\dots\dots\dots\dots\dots\dots\dots\dots\dots\dots
\dots\dots\dots\dots\dots
}
\\
& 
\underbrace{
\underbrace{
k_{d:1:1},\dots,k_{d:1:\lambda_d}}_{\lambda_d},
\dots,
\underbrace{
k_{d:\mu_d:1},\dots,k_{d:\mu_d:\lambda_d}}_{\lambda_d}}_{
\mu_d}.
\endaligned
\end{equation}
Correspondingly, we identify the set 
\def\theequation{3.35}\begin{equation}
\small
\left\{
1,\dots,\lambda_1,\dots,\mu_1\lambda_1,
\dots\dots,
\mu_1\lambda_1+\mu_2\lambda_2,
\dots\dots,
\mu_1\lambda_1+\mu_2\lambda_2
+
\cdots
+
\mu_d\lambda_d
\right\}
\end{equation}
of all integers $\alpha$ from $1$ to $\mu_1 \lambda_1 + \mu_2
\lambda_2 + \cdots + \mu_d \lambda_d$ with the following specific set
\def\theequation{3.36}\begin{equation}
\small
\{
\underbrace{
\underbrace{
\underbrace{
\underbrace{
1\!\!:\!\!1\!\!:\!\!1,
\dots,
1\!\!:\!\!1\!\!:\!\!\lambda_1}_{\lambda_1},
\dots,
1\!\!:\!\!\mu_1\!\!:\!\!\lambda_1}_{\mu_1\lambda_1},
\dots,
2\!:\!\mu_2\!:\!\lambda_2}_{\mu_1\lambda_1+\mu_2\lambda_2},
\dots,
d\!:\!\mu_d\!:\!\lambda_d}_{\mu_1\lambda_1+\mu_2\lambda_2
+\cdots+\mu_d\lambda_d}
\},
\end{equation}
written in a lexicographic way which emphasizes clearly the
subdivision in $d$ collections of $\mu_d$ groups of $\lambda_d$
integers.
\end{notationalconvention}

With this notation at hand, we see that the development, in
several independent variables, of the general monomial $\left(
y_{\lambda_1 } \right)^{ \mu_1}
\cdots \left( y_{ \lambda_d } \right)^{ \mu_d }$, may be written as
follows:
\def\theequation{3.37}\begin{equation}
\aligned
y_{k_{1:1:1},\dots,k_{1:1:\lambda_1}}
\cdots \
y_{k_{1:\mu_1:1},\dots,k_{1:\mu_1:\lambda_1}}
\cdots \
y_{k_{d:1:1},\dots,k_{d:1:\lambda_d}}
\cdots\cdots\,
y_{k_{d:\mu_d:1},\dots,k_{d:\mu_d:\lambda_d}}.
\endaligned
\end{equation} 
Formally speaking, this expression is better than~\thetag{ 3.32}. Using
product symbols, we may even write it under the 
slightly more compact form
\def\theequation{3.38}\begin{equation}
\prod_{1\leq\nu_1\leq\mu_1}\,
y_{k_{1:\nu_1:1},\dots,k_{1:\nu_1:\lambda_1}}
\cdots \
\prod_{1\leq\nu_d\leq\mu_d}\,
y_{k_{d:\nu_d:1},\dots,k_{d:\nu_d:\lambda_d}}.
\end{equation}

Now that we have translated the monomial, we may add all the summation
symbols: the general expression of ${\bf Y}_\kappa$ (which generalizes
our three previous examples~\thetag{ 3.26}) will be of the form:
\def\theequation{3.39}\begin{equation}
\aligned
{\bf Y}_\kappa
&
=
\mathcal{Y}_{x^{i_1}\cdots x^{i_\kappa}}
+
\sum_{d=1}^{\kappa+1}
\ \
\sum_{1\leq\lambda_1<\cdots<\lambda_d\leq\kappa}
\ \
\sum_{\mu_1\geq 1,\dots,\mu_d\geq 1} 
\
\sum_{
\mu_1\lambda_1
+
\cdots
+
\mu_d\lambda_d\leq \kappa+1} 
\\
&
\sum_{k_{1:1:1},\dots,k_{1:1:\lambda_1}=1}^n
\cdots \
\sum_{k_{1:\mu_1:1},\dots,k_{1:\mu_1:\lambda_1}=1}^n
\cdots\cdots \
\sum_{k_{d:1:1},\dots,k_{d:1:\lambda_d}=1}^n
\cdots \
\sum_{k_{d:\mu_d:1},\dots,k_{d:\mu_d:\lambda_d}=1}^n
\\
&
\text{\bf[?]}
\prod_{1\leq\nu_1\leq\mu_1}\,
y_{k_{1:\nu_1:1},\dots,k_{1:\nu_1:\lambda_1}}
\ \cdots \
\prod_{1\leq\nu_d\leq\mu_d}\,
y_{k_{d:\nu_d:1},\dots,k_{d:\nu_d:\lambda_d}}.
\endaligned
\end{equation}
From now on, 
up to the end of the article, to be very precise, we will restitute
the bounds $\sum_{ k = 1 }^n$ of all the previously abbreviated sums
$\sum_k$. This is justified by the fact that, since we shall deal in
Section~5 below simultaneously with several independent variables
$(x^1, \dots, x^n)$ and with several dependent variables $(y^1, \dots,
y^m)$, we shall encounter sums $\sum_{ l = 1 }^m$, not to be confused
with sums $\sum_{ k = 1 }^n$.

\subsection*{3.40.~Combinatorics of the Kronecker symbols} 
Our next task is to determine what appears inside the brackets {\bf
[?]} of the above equation. We will treat this rather delicate
question very progressively. Inductively, we have to guess how we may
pass from the bracketed term of~\thetag{ 2.25}, namely from
\def\theequation{3.41}\begin{equation}
\aligned
&
\left[
\frac{\kappa\cdots(\kappa-\mu_1\lambda_1-\cdots-\mu_d\lambda_d+1)}
{(\lambda_1!)^{\mu_1}\,\mu_1!
\cdots
(\lambda_d!)^{\mu_d}\,\mu_d!
}
\cdot
\mathcal{Y}_{
x^{\kappa-\mu_1\lambda_1-\cdots-\mu_d\lambda_d}
\,
y^{\mu_1+\cdots+\mu_d}
}
-
\right. \\
& \
\ \ \ \ \ \ \ \ \ \
\left.
-
\frac{\kappa\cdots( 
\kappa-\mu_1\lambda_1-\cdots-\mu_d\lambda_d+2)
(\mu_1\lambda_1+\cdots+\mu_d\lambda_d)}
{(\lambda_1!)^{\mu_1}\,\mu_1!
\cdots
(\lambda_d!)^{\mu_d}\,\mu_d!
}
\cdot
\right.
\\
& \
\ \ \ \ \ \ \ \ \ \ \ \ \ \ \ \ \
\cdot
\mathcal{X}_{
x^{\kappa-\mu_1\lambda_1-\cdots-\mu_d\lambda_d+1}
\,
y^{\mu_1+\cdots+\mu_d-1}
}
\Big],
\endaligned
\end{equation}
to the corresponding (still unknown) bracketed term 
{\bf [?]}.

First of all, we examine the following term, extracted from the
complete expression of ${\bf Y}_{ i_1, i_2, i_3}$ (first line
of~\thetag{ 3.24}):
\def\theequation{3.42}\begin{equation}
\sum_{k_1=1}^n\,
\left[
\delta_{i_1}^{k_1}\,\mathcal{Y}_{x^{i_2}x^{i_3}y}
+
\delta_{i_2}^{k_1}\,\mathcal{Y}_{x^{i_1}x^{i_3}y}
+
\delta_{i_3}^{k_1}\,\mathcal{Y}_{x^{i_1}x^{i_2}y}
-
\mathcal{ X}_{x^{i_1}x^{i_2}x^{i_3}}^{k_1}
\right]
y_{k_1}.
\end{equation}
Here, the coefficient $\left[ 3\, \mathcal{ Y}_{ x^2 y} - \mathcal{
X}_{ x^3 } \right]$ of the monomial $y_1$ in ${\bf Y}_3$ is replaced
by the above bracketed terms. 

Let us precisely analyze the combinatorics. Here, $\mathcal{ X}_{x^3}$ is
replaced by $\mathcal{ X}_{x^{ i_1 }x^{ i_2}x^{i_3}}^{k_1}$, where the
lower indices $i_1, i_2, i_3$ come from ${\bf Y}_{i_1, i_2, i_3}$ and
where the upper index $k_1$ is the summation index. Also, the integer
$3$ in $3\, \mathcal{ Y}_{x^2 y}$ is replaced by a sum of exactly
three terms, each involving a single Kronecker symbol $\delta_i^k$, in
which the lower index is always an index $i= i_1, i_2, i_3$ and in
which the upper index is always equal to the summation index $k_1$.
By the way, more generally, we immediately observe that
all the successive positive integers 
\def\theequation{3.43}\begin{equation}
1, 3, 1, 3, 3, 1, 3, 1, 3, 3, 3, 9, 6, 3, 1, 3, 4
\end{equation}
appearing in the formula~\thetag{ 2.8} for ${\bf Y}_3$ are replaced,
in the formula~\thetag{ 3.24} for ${\bf Y}_{i_1, i_2, i_3}$, by sums
of exactly the same number of terms involving Kronecker
symbols. This observation will be a precious guide. Finally, in the
symbol $\delta_i^{k_1}$, if $i$ is chosen among the set $\{ i_1, i_2,
i_3\}$, for instance if $i = i_1$, it follows that the development of
$\mathcal{ Y}_{x^2y}$ necessarily involves the remaining indices, for
instance $\mathcal{ Y}_{x^{i_2}x^{i_3}y}$. Since there are three
choices for $i = i_1, i_2, i_3$, we recover the number $3$.

Next, comparing $\left[ \mathcal{ Y}_{yy}
-2\,\mathcal{ X}_{ xy} \right] (y_1)^2$ with 
the term
\def\theequation{3.44}\begin{equation}
\sum_{k_1,k_2=1}^n\,
\left[
\delta_{i_1,\,i_2}^{k_1,k_2}\,\mathcal{Y}_{yy}
-
\delta_{i_1}^{k_1}\,\mathcal{X}_{x^{i_2}y}^{k_1}
-
\delta_{i_2}^{k_1}\,\mathcal{X}_{x^{i_1}y}^{k_1}
\right]
y_{k_1}y_{k_2},
\end{equation}
extracted from the complete expression of ${\bf Y}_{i_1,i_2}$ (second
line of~\thetag{ 3.18}), we learn and we guess that the number of
Kronecker symbols before $\mathcal{ Y}_{x^\gamma y^\delta}$ must be
equal to the number of indices $k_\alpha$ minus $\gamma$. This rule
is confirmed by examining the term (second and third line of~\thetag{
3.24})
\def\theequation{3.45}\begin{equation}
\aligned
\sum_{k_1,k_2}
&
\left[
\delta_{i_1,\,i_2}^{k_1,k_2}\,\mathcal{Y}_{x^{i_3}y^2}
+
\delta_{i_3,\,i_1}^{k_1,k_2}\,\mathcal{Y}_{x^{i_2}y^2}
+
\delta_{i_2,\,i_3}^{k_1,k_2}\,\mathcal{Y}_{x^{i_1}y^2}
-
\right.
\\
& \
\ \ \ \ \
\left.
-
\delta_{i_1}^{k_1}\,\mathcal{X}_{x^{i_2}x^{i_3}y}^{k_2}
-
\delta_{i_2}^{k_1}\,\mathcal{X}_{x^{i_1}x^{i_3}y}^{k_2}
-
\delta_{i_3}^{k_1}\,\mathcal{X}_{x^{i_1}x^{i_2}y}^{k_2}
\right]
y_{k_1}y_{k_2},
\endaligned
\end{equation}
developing $\left[ 3\,\mathcal{ Y}_{xy^2} - 3\, \mathcal{
X}_{x^2y} \right] (y_1)^2$.

Also, we may examine the following term
\def\theequation{3.46}\begin{equation}
\aligned
\sum_{k_1,k_2=1}^n\,
&
\left[
\delta_{i_1,\,i_2}^{k_1,k_2}\,\mathcal{Y}_{x^{i_3}x^{i_4}y^2}
+
\delta_{i_1,\,i_3}^{k_1,k_2}\,\mathcal{Y}_{x^{i_2}x^{i_4}y^2}
+
\delta_{i_1,\,i_4}^{k_1,k_2}\,\mathcal{Y}_{x^{i_2}x^{i_3}y^2}
+ 
\right.
\\
& \
\ \ \ \ \
\left.
+
\delta_{i_2,\,i_3}^{k_1,k_2}\,\mathcal{Y}_{x^{i_1}x^{i_4}y^2}
+
\delta_{i_2,\,i_4}^{k_1,k_2}\,\mathcal{Y}_{x^{i_1}x^{i_3}y^2}
+
\delta_{i_3,\,i_4}^{k_1,k_2}\,\mathcal{Y}_{x^{i_1}x^{i_2}y^2}
-
\right.
\\
& \
\ \ \ \ \
\left.
-
\delta_{i_1}^{k_1}\,\mathcal{X}_{x^{i_2}x^{i_3}x^{i_4}y}^{k_1}
-
\delta_{i_2}^{k_1}\,\mathcal{X}_{x^{i_1}x^{i_2}x^{i_3}y}^{k_1}
-
\delta_{i_3}^{k_1}\,\mathcal{X}_{x^{i_1}x^{i_2}x^{i_4}y}^{k_1}
-
\right.
\\
& \
\ \ \ \ \
\left.
-
\delta_{i_4}^{k_1}\,\mathcal{X}_{x^{i_1}x^{i_2}x^{i_3}y}^{k_1}
\right]
y_{k_1}y_{k_2},
\endaligned
\end{equation} 
extracted from ${\bf Y}_{ i_1, i_2, i_3, i_4}$ and developing
$\left[ 6\, \mathcal{ Y}_{ x^2 y^2} - 4\, \mathcal{ X}_{ x^3 y}
\right] (y_1)^2$. We would like to mention that we have not written
the complete expression of ${\bf Y}_{ i_1, i_2, i_3, i_4}$, because it
would cover two and a half printed pages. 

By inspecting the way how the indices are permuted in the multiple
Kronecker symbols of the first two lines of this expression~\thetag{
3.46}, we observe that the six terms correspond exactly to the six
possible choices of two complementary ordered couples of integers in
the set $\{ 1, 2, 3, 4\}$, namely
\def\theequation{3.47}\begin{equation}
\aligned
&
\{1,2\}\cup\{3,4\}, 
\ \ \ \ \ \ \ \
\{1,3\}\cup\{2,4\}, 
\ \ \ \ \ \ \ \
\{1,4\}\cup\{2,3\}, 
\\
&
\{2,3\}\cup\{1,4\}, 
\ \ \ \ \ \ \ \
\{2,4\}\cup\{1,3\}, 
\ \ \ \ \ \ \ \
\{3,4\}\cup\{1,2\}.
\endaligned
\end{equation}
At this point, we start to devise the general combinatorics. Before
proceeding further, we need some notation.

\subsection*{ 3.48.~Permutation groups}
For every $p \in \N$ with $p \geq 1$, we denote by $\mathfrak{ S}_p$ the
full permutation group of the set $\{ 1, 2, \dots, p-1, p\}$. Its
cardinal equals $p!$. The letters $\sigma$ and
$\tau$ will be used to denote an
element of $\mathfrak{ S}_p$. If $p \geq 2$, and if $q \in \N$
satisfies $1\leq q \leq p-1$, we denote by $\mathfrak{ S}_p^q$ the
subset of permutations $\sigma \in \mathfrak{ S}_p$ satisfying
the two collections of inequalities
\def\theequation{3.49}\begin{equation}
\sigma(1)<\sigma(2)<\cdots<\sigma(q)
\ \ \ \ \ \ \ \ \ 
{\rm and}
\ \ \ \ \ \ \ \ \ 
\sigma(q+1)<\sigma(q+2)<\cdots<\sigma(p).
\end{equation}
The cardinal of $\mathfrak{ S}_p^q$ 
equals $C_p^q = \frac{ p!}{ q! \ (p-q) !}$.

\def\thelemma{3.50}\begin{lemma}
For $\kappa \geq 1$, the development of~\thetag{ 2.20} to several
independent variables $(x^1, \dots, x^n)$ is{\rm :}
\def\theequation{3.51}\begin{equation}
\small
\aligned
&
{\bf Y}_{i_1,i_2,\dots,i_\kappa}
=
\mathcal{Y}_{x^{i_1}x^{i_2}\cdots x^{i_\kappa}}
+
\sum_{k_1=1}^n\,
\left[
\sum_{\tau\in\mathfrak{S}_\kappa^1}\,
\delta_{i_{\tau(1)}}^{k_1}\,\mathcal{Y}_{x^{i_{\tau(2)}}
\cdots x^{i_{\tau(\kappa)}}y}
-
\mathcal{X}_{x^{i_1}x^{i_2}\cdots x^{i_{\kappa}}}^{k_1}
\right]
y_{k_1}
+ \\
& \
\ \ \ \ \
+
\sum_{k_1,k_2=1}^n\,
\left[
\sum_{\tau\in\mathfrak{S}_\kappa^2}\,
\delta_{i_{\tau(1)},i_{\tau(2)}}^{k_1,\ \ \ k_2}\,
\mathcal{Y}_{x^{i_{\tau(3)}}\cdots x^{i_{\tau(\kappa)}}y^2}
-
\sum_{\tau\in\mathfrak{S}_\kappa^1}\,
\delta_{i_{\tau(1)}}^{k_1}\,
\mathcal{X}_{x^{i_{\tau(2)}}\cdots x^{i_{\tau(\kappa)}}y}^{k_2}
\right]
y_{k_1}y_{k_2}
+
\\
& \
\ \ \ \ \
+
\sum_{k_1, k_2, k_3=1}^n
\left[
\sum_{\tau\in\mathfrak{S}_\kappa^3}\,
\delta_{i_{\tau(1)},i_{\tau(2)},i_{\tau(3)}}^{
k_1,\ \ \ k_2,\ \ \ k_3}\,
\mathcal{Y}_{x^{i_{\tau(4)}}\cdots x^{i_{\tau(\kappa)}}y^3}
-
\right.
\\
& \
\ \ \ \ \ \ \ \ \ \ \ \ \ \ \ 
\ \ \ \ \ \ \ \ \ \ \ \ \ \ \ 
\ \ \ \ \ \ 
\left.
-
\sum_{\tau\in\mathfrak{S}_\kappa^2}\,
\delta_{i_{\tau(1),i_{\tau(2)}}}^{k_1,\ \ \ k_2}
\mathcal{X}_{x^{i_{\tau(3)}}\cdots x^{i_{\tau(\kappa)}}y^2}^{k_3}
\right]
y_{k_1}y_{k_2}y_{k_3}
+ \\
& \
\ \ \ \ \
+
\cdots\cdots
+
\\
& \
\ \ \ \ \
+
\sum_{k_1,\dots,k_\kappa=1}^n
\left[
\delta_{i_1,\dots,\ i_\kappa}^{k_1,\dots,k_\kappa}\,
\mathcal{Y}_{y^\kappa}
-
\sum_{\tau\in\mathfrak{S}_\kappa^{\kappa-1}}\,
\delta_{i_{\tau(1)},\dots,i_{\tau(\kappa-1)}}^{
k_1,\dots\dots,k_{\kappa-1}}\,
\mathcal{X}_{x^{i_{\tau(\kappa)}}y^{\kappa-1}}^{k_\kappa}
\right]
y_{k_1}\cdots y_{k_\kappa}
+ \\
& \
\ \ \ \ \
+
\sum_{k_1,\dots,k_\kappa,k_{\kappa+1}=1}^n
\left[
-
\delta_{i_1,\dots,\ i_\kappa}^{k_1,\dots,k_\kappa}\,
\mathcal{X}_{y^\kappa}^{k_{\kappa+1}}
\right]
y_{k_1}\cdots y_{k_\kappa} y_{k_{\kappa+1}}
+
{\sf remainder}.
\endaligned
\end{equation}
Here, the term {\sf remainder} collects all
remaining monomials in the pure jet variables
$y_{ k_1, \dots, k_\lambda }$.
\end{lemma}

\subsection*{3.52.~Continuation}
Thus, we have devised how the part of ${\bf Y}_{i_1,\dots, i_\kappa}$
which involves only the jet variables $y_{k_\alpha}$ must be written.
To proceed further, we shall examine the following term, extracted
from ${\bf Y}_{i_1,i_2,i_3}$ (lines 12 and 13 of~\thetag{ 3.24})
\def\theequation{3.53}\begin{equation}
\aligned
& \
\ \ \ \ \
\sum_{k_1,k_2,k_3,k_4}\,
\left[
-
\delta_{i_1,\,i_2,\,i_3}^{k_1,k_2,k_3}\,\mathcal{X}_{y^2}^{k_4}
-
\delta_{i_1,\,i_2,\,i_3}^{k_2,k_3,k_1}\,\mathcal{X}_{y^2}^{k_4}
-
\delta_{i_1,\,i_2,\,i_3}^{k_3,k_2,k_1}\,\mathcal{X}_{y^2}^{k_4}
- 
\right.
\\
& \
\ \ \ \ \ \ \ \ \ \ \ \ \ \ \ \ \ \ \ \ \ \ \ \ \ \ \
\left.
-
\delta_{i_1,\,i_2,\,i_3}^{k_3,k_4,k_1}\,\mathcal{X}_{y^2}^{k_2}
-
\delta_{i_1,\,i_2,\,i_3}^{k_3,k_1,k_4}\,\mathcal{X}_{y^2}^{k_2}
-
\delta_{i_1,\,i_2,\,i_3}^{k_1,k_3,k_4}\,\mathcal{X}_{y^2}^{k_2}
\right]
y_{k_1}y_{k_2}y_{k_3,k_4},
\endaligned
\end{equation}
which developes the term $\left[ - 6\,\mathcal{ X}_{ y^2} \right]
(y_1)^2 y_2$ of ${\bf Y}_3$ (third line of~\thetag{ 2.8}). During the
computation which led us to the final expression~\thetag{ 3.24}, we
organized the formula in order that, in the six Kronecker symbols, the
lower indices $i_1,i_2,i_3$ are all written in the same order. But
then, {\it what is the rule for the appearance of the four upper
indices $k_1, k_2, k_3, k_4$}?

In April 2001, we discovered the rule by inspecting both~\thetag{
3.53} and the following complicated term, extracted from the complete
expression of ${\bf Y}_{i_1,i_2,i_3,i_4}$ written in one of our
manuscripts:
\def\theequation{3.54}\begin{equation}
\aligned
\sum_{k_1,k_2,k_3}\,
&
\left[
\delta_{i_1, \ i_2, \ i_3}^{k_1,k_2,k_3}\,\mathcal{Y}_{x^{i_4}y^2}
+
\delta_{i_1, \ i_2, \ i_3}^{k_2,k_1,k_3}\,\mathcal{Y}_{x^{i_4}y^2}
+
\delta_{i_1, \ i_2, \ i_3}^{k_2,k_3,k_1}\,\mathcal{Y}_{x^{i_4}y^2}
+
\right. 
\\
& \
\ \ 
+
\delta_{i_1, \ i_2, \ i_4}^{k_1,k_2,k_3}\,\mathcal{Y}_{x^{i_3}y^2}
+
\delta_{i_1, \ i_2, \ i_4}^{k_2,k_1,k_3}\,\mathcal{Y}_{x^{i_3}y^2}
+
\delta_{i_1, \ i_2, \ i_4}^{k_2,k_3,k_1}\,\mathcal{Y}_{x^{i_3}y^2}
+ \\
& \
\ \ 
+
\delta_{i_1, \ i_3, \ i_4}^{k_1,k_2,k_3}\,\mathcal{Y}_{x^{i_2}y^2}
+
\delta_{i_1, \ i_3, \ i_4}^{k_2,k_1,k_3}\,\mathcal{Y}_{x^{i_2}y^2}
+
\delta_{i_1, \ i_3, \ i_4}^{k_2,k_3,k_1}\,\mathcal{Y}_{x^{i_2}y^2}
+ \\
& \
\ \ 
+
\delta_{i_2, \ i_3, \ i_4}^{k_1,k_2,k_3}\,\mathcal{Y}_{x^{i_1}y^2}
+
\delta_{i_2, \ i_3, \ i_4}^{k_2,k_1,k_3}\,\mathcal{Y}_{x^{i_1}y^2}
+
\delta_{i_2, \ i_3, \ i_4}^{k_2,k_3,k_1}\,\mathcal{Y}_{x^{i_1}y^2}
- \\
& \
\ \ 
-
\delta_{i_1,\ i_2}^{k_1,k_2}\,\mathcal{X}_{x^{i_3}x^{i_4}y}^{k_3}
-
\delta_{i_1,\ i_2}^{k_2,k_1}\,\mathcal{X}_{x^{i_3}x^{i_4}y}^{k_3}
-
\delta_{i_1,\ i_2}^{k_2,k_3}\,\mathcal{X}_{x^{i_3}x^{i_4}y}^{k_1}
- \\
& \
\ \ 
-
\delta_{i_1,\ i_3}^{k_1,k_2}\,\mathcal{X}_{x^{i_2}x^{i_4}y}^{k_3}
-
\delta_{i_1,\ i_3}^{k_2,k_1}\,\mathcal{X}_{x^{i_2}x^{i_4}y}^{k_3}
-
\delta_{i_1,\ i_3}^{k_2,k_3}\,\mathcal{X}_{x^{i_2}x^{i_4}y}^{k_1}
- \\
& \
\ \ 
-
\delta_{i_1,\ i_4}^{k_1,k_2}\,\mathcal{X}_{x^{i_2}x^{i_3}y}^{k_3}
-
\delta_{i_1,\ i_4}^{k_2,k_1}\,\mathcal{X}_{x^{i_2}x^{i_3}y}^{k_3}
-
\delta_{i_1,\ i_4}^{k_2,k_3}\,\mathcal{X}_{x^{i_2}x^{i_3}y}^{k_1}
- \\
& \
\ \ 
-
\delta_{i_2,\ i_3}^{k_1,k_2}\,\mathcal{X}_{x^{i_1}x^{i_4}y}^{k_3}
-
\delta_{i_2,\ i_3}^{k_2,k_1}\,\mathcal{X}_{x^{i_1}x^{i_4}y}^{k_3}
-
\delta_{i_2,\ i_3}^{k_2,k_3}\,\mathcal{X}_{x^{i_1}x^{i_4}y}^{k_1}
- \\
& \
\ \ 
-
\delta_{i_2,\ i_4}^{k_1,k_2}\,\mathcal{X}_{x^{i_1}x^{i_3}y}^{k_3}
-
\delta_{i_2,\ i_4}^{k_2,k_1}\,\mathcal{X}_{x^{i_1}x^{i_3}y}^{k_3}
-
\delta_{i_2,\ i_4}^{k_2,k_3}\,\mathcal{X}_{x^{i_1}x^{i_3}y}^{k_1}
- \\
& \
\ \
\left. 
-
\delta_{i_3,\ i_4}^{k_1,k_2}\,\mathcal{X}_{x^{i_1}x^{i_2}y}^{k_3}
-
\delta_{i_3,\ i_4}^{k_2,k_1}\,\mathcal{X}_{x^{i_1}x^{i_2}y}^{k_3}
-
\delta_{i_3,\ i_4}^{k_2,k_3}\,\mathcal{X}_{x^{i_1}x^{i_2}y}^{k_1}
\right]
y_{k_1}y_{k_2,k_3}.
\endaligned
\end{equation}
This sum developes the term $\left[ 12\, \mathcal{ Y}_{xy^2} -
18\,\mathcal{ X}_{ x^2 y} \right]y_1y_2$ of
${\bf Y}_3$ (third line of~\thetag{
2.9}). Let us explain what are the formal rules.

In the bracketed terms of~\thetag{ 3.53}, there are no permutation of
the indices $i_1,i_2,i_3$, but there is a certain unknown subset of
all the permutations of the four indices $k_1,k_2,k_3,k_4$. In the
bracketed terms of~\thetag{ 3.54}, two combinatorics are present:

\smallskip

\begin{itemize}
\item[$\bullet$]
there are some permutations of the indices $i_1,i_2,i_3,i_4$ and
\item[$\bullet$]
there are some permutations of the indices $k_1,k_2,k_3$. 
\end{itemize}

\smallskip

Here, the permutations of the indices $i_1,i_2,i_3,i_4$ are easily
guessed, since they are the same as the permutations which were
introduced in \S3.48 above. Indeed, in the first four
lines of~\thetag{ 3.54}, we
see the four decompositions 
\def\theequation{3.55}\begin{equation}
\{i_1,i_2,i_3\}\cup\{i_4\},
\ \ \ \ \ \ \ 
\{i_1,i_2,i_4\}\cup\{i_3\},
\ \ \ \ \ \ \ 
\{i_1,i_3,i_4\}\cup\{i_2\},
\ \ \ \ \ \ \ 
\{i_2,i_3,i_4\}\cup\{i_1\},
\end{equation}
of the set $\{i_1,i_2,i_3,i_4\}$, and in the 
last six lines of~\thetag{ 3.54}, we see the
six decompositions
\def\theequation{3.56}\begin{equation}
\aligned
&
\{i_1,i_2\}\cup\{i_3,i_4\},
\ \ \ \ \ \ \ 
\{i_1,i_3\}\cup\{i_2,i_4\},
\ \ \ \ \ \ \ 
\{i_1,i_4\}\cup\{i_2,i_3\},
\\
&
\{i_2,i_3\}\cup\{i_1,i_4\},
\ \ \ \ \ \ \ 
\{i_2,i_4\}\cup\{i_1,i_3\},
\ \ \ \ \ \ \ 
\{i_3,i_4\}\cup\{i_1,i_2\},
\endaligned
\end{equation}
so that~\thetag{ 3.54} may be written under the form
\def\theequation{3.57}\begin{equation}
\small
\aligned
\sum_{k_1,k_2,k_3}\,
\left[
\sum_{\tau\in\mathfrak{S}_4^3}\,
\sum_{\sigma\in\text{\bf ?}}\,
\delta_{i_{\tau(1)},i_{\tau(2)},i_{\tau(3)}}^{
k_{\tau(1)},k_{\tau(2)},k_{\tau(3)}}\,
\mathcal{Y}_{x^{i_{\tau(4)}}y^2}
-
\sum_{\tau\in\mathfrak{S}_4^2}\,
\sum_{\sigma\in\text{\bf ?}}\,
\delta_{i_{\tau(1)},i_{\tau(2)}}^{
k_{\tau(1)},k_{\tau(2)}}\,
\mathcal{X}_{x^{i_{\tau(3)}}x^{i_{\tau(4)}}y}^{k_{\tau(3)}}
\right]
y_{k_1}y_{k_2,k_3},
\endaligned
\end{equation}
where in the two above sums $\sum_{ \sigma \in \text{\bf ?}}$, the
letter $\sigma$ denotes a permutation of the set $\{1,2,3\}$ and where
the sign {\bf ?} refers to two (still unknown) subset of the full
permutation group $\mathfrak{ S}_3$. {\it The only remaining question
is to determine how the indices $k_\alpha$ are permuted in~\thetag{
3.53} and in~\thetag{ 3.54}}.

The answer may be guessed by looking at the permutations of the set
$\{k_1,k_2,k_3,k_4\}$ which stabilize the monomial
$y_{k_1}y_{k_2}y_{k_3,k_4}$ in~\thetag{ 3.53}: we clearly have the
following four symmetry relations between monomials:
\def\theequation{3.58}\begin{equation}
y_{k_1}y_{k_2}y_{k_3,k_4}
\equiv
y_{k_2}y_{k_1}y_{k_3,k_4}
\equiv
y_{k_1}y_{k_2}y_{k_4,k_3}
\equiv
y_{k_2}y_{k_1}y_{k_4,k_3},
\end{equation}
and nothing more. 
Then the number $6$ of bracketed terms in~\thetag{ 3.53}
is exactly equal to the cardinal $24 = 4!$ of the full permutation
group of the set $\{k_1,k_2,k_3,k_4\}$ divided by the number $4$ of
these symmetry relations. The set of permutations $\sigma$ of
$\{1,2,3,4\}$ satisfying these symmetry relations
\def\theequation{3.59}\begin{equation}
y_{k_{\sigma(1)}}y_{k_{\sigma(2)}}
y_{k_{\sigma(3)},k_{\sigma(4)}}
\equiv
y_{k_1}y_{k_2}y_{k_3,k_4}
\end{equation}
consitutes a subgroup of $\mathfrak{ S}_4$ which we will denote by
$\mathfrak{ H}_4^{(2,1),(1,2)}$. Furthermore, the coset
\def\theequation{3.60}\begin{equation}
\mathfrak{ F}_4^{(2,1),(1,2)} 
:=
\mathfrak{ S}_4 / \mathfrak{ H}_4^{(2,1),(1,2)}
\end{equation}
possesses the six representatives 
\def\theequation{3.61}\begin{equation}
\small
\aligned
\left(
\begin{array}{cccc}
1 & 2 & 3 & 4 \\
1 & 2 & 3 & 4 \\
\end{array}
\right), 
\ \ \ \ \ \ \ \ \ \ 
\left(
\begin{array}{cccc}
1 & 2 & 3 & 4 \\
2 & 3 & 1 & 4 \\
\end{array}
\right), 
\ \ \ \ \ \ \ \ \ \ 
\left(
\begin{array}{cccc}
1 & 2 & 3 & 4 \\
3 & 2 & 1 & 4 \\
\end{array}
\right), 
\\
\left(
\begin{array}{cccc}
1 & 2 & 3 & 4 \\
3 & 4 & 1 & 2 \\
\end{array}
\right), 
\ \ \ \ \ \ \ \ \ \ 
\left(
\begin{array}{cccc}
1 & 2 & 3 & 4 \\
3 & 1 & 4 & 2 \\
\end{array}
\right), 
\ \ \ \ \ \ \ \ \ \ 
\left(
\begin{array}{cccc}
1 & 2 & 3 & 4 \\
1 & 3 & 4 & 2 \\
\end{array}
\right), \\
\endaligned
\end{equation}
which exactly appear as the permutations of the upper indices of our
example~\thetag{ 3.53}. Of course, the question arises whether the
choice of such six representatives in the quotient $\mathfrak{ S}_4 /
\mathfrak{ H}_4^{ (2,1), (1,2)}$ is legitimate.

Fortunately, we observe that after conjugation by any permutation
$\sigma \in \mathfrak{ H }_4^{ (2,1), (1,2)}$, we do not perturb any of
the six terms of~\thetag{ 3.53}, for instance the third term
of~\thetag{ 3.53} is not perturbed, as shown by the following
computation
\def\theequation{3.62}\begin{equation}
\small
\aligned
&
\sum_{k_1,k_2,k_3,k_4}
\left[
-
\delta_{i_1, \ \ \ \ i_2, \ \ \ \ i_3}^{
k_{\sigma(3)}, k_{\sigma(2)}, k_{\sigma(1)}}\,
\mathcal{X}_{y^2}^{k_{\sigma(4)}}
\right]
y_{k_1}y_{k_2}y_{k_3,k_4}
= \\
& \
\ \ \ \ \
=
\sum_{k_1,k_2,k_3,k_4}
\left[
-
\delta_{i_1, \ i_2, \ i_3}^{
k_3, k_2, k_1}\,
\mathcal{X}_{y^2}^{k_{\sigma(4)}}
\right]
y_{k_{\sigma^{-1}(1)}}
y_{k_{\sigma^{-1}(2)}}
y_{k_{\sigma^{-1}(3)},k_{\sigma^{-1}(4)}}
\\
& \
\ \ \ \ \
=
\sum_{k_1,k_2,k_3,k_4}
\left[
-
\delta_{i_1, \ i_2, \ i_3}^{
k_3, k_2, k_1}\,
\mathcal{X}_{y^2}^{k_{\sigma(4)}}
\right]
y_{k_1}y_{k_2}y_{k_3,k_4}
\endaligned
\end{equation} 
thanks to the symmetry~\thetag{ 3.59}. Thus, as expected, the choice of
$6$ arbitrary representatives $\sigma \in \mathfrak{ F}_4^{(2,1), (1,
2)}$ in the bracketed terms of~\thetag{ 3.53} is free.
In conclusion, we have shown that~\thetag{ 3.53} may 
be written under the form:
\def\theequation{3.63}\begin{equation}
\aligned
\sum_{k_1,k_2,k_3,k_4}\,
\left[
-
\sum_{\sigma\in\mathfrak{F}_4^{(2,1),(1,2)}}\,
\delta_{i_1,\ \ \ \ i_2,\ \ \ \ i_3}^{
k_{\sigma(1)},k_{\sigma(2)},k_{\sigma(3)}}\,
\mathcal{X}_{y^2}^{k_{\sigma(4)}}
\right]
y_{k_1}y_{k_2}y_{k_3,k_4},
\endaligned
\end{equation}

This rule is confirmed by inspecting~\thetag{ 3.54} (as
well as all the other terms of ${\bf Y}_{i_1, i_2, i_3}$ and
of ${\bf Y}_{ i_1, i_2, i_3, i_4}$). Indeed, the
permutations $\sigma$ of the set $\{k_1, k_2, k_3\}$ which stabilize the
monomial $y_{k_1} y_{k_2, k_3}$ consist just of the identity
permutation and the transposition of $k_2$ and $k_3$. The coset
$\mathfrak{ S }_3 / \mathfrak{ H }_3^{ (1,1), (1,2)}$ has the three
representatives
\def\theequation{3.64}\begin{equation}
\small
\aligned
\left(
\begin{array}{ccc}
1 & 2 & 3 \\
1 & 2 & 3 \\
\end{array}
\right), 
\ \ \ \ \ \ \ \ \ \ 
\left(
\begin{array}{ccc}
1 & 2 & 3 \\
2 & 1 & 3 \\
\end{array}
\right), 
\ \ \ \ \ \ \ \ \ \ 
\left(
\begin{array}{cccc}
1 & 2 & 3 \\
2 & 3 & 1 \\
\end{array}
\right), 
\\
\endaligned
\end{equation}
which appear in the upper index position of each of the ten lines
of~\thetag{ 3.54}. It follows that~\thetag{ 3.54}
may be written under the form
\def\theequation{3.65}\begin{equation}
\aligned
\sum_{k_1,k_2,k_3}\,
&
\left[
\sum_{\tau\in\mathfrak{S}_4^3}\,
\sum_{\sigma\in\mathfrak{F}_3^{(1,1),(1,2)}}\,
\delta_{i_{\tau(1)},i_{\tau(2)},i_{\tau(3)}}^{
k_{\sigma(1)},k_{\sigma(2)},k_{\sigma(3)}}\,
\mathcal{Y}_{x^{i_{\tau(4)}}y^2}
- 
\right.
\\
& \
\left.
-
\sum_{\sigma\in\mathfrak{S}_4^2}\,
\sum_{\tau\in\mathfrak{F}_3^{(1,1),(1,2)}}\,
\delta_{i_{\tau(1)},i_{\tau(2)}}^{
k_{\sigma(1)},k_{\sigma(2)}}\,
\mathcal{X}_{x^{i_{\tau(3)}}
x^{i_{\tau(4)}}y}^{k_{\sigma(3)}}
\right]
y_{k_1}y_{k_2,k_3}.
\endaligned
\end{equation}

\subsection*{ 3.66.~General complete expression of ${\bf Y}_{i_1, \dots,
i_\kappa}$} As in the incomplete expression~\thetag{ 3.39} of ${\bf
Y}_{i_1, \dots, i_\kappa}$, consider integers $1\leq \lambda_1 <
\cdots < \lambda_d \leq \kappa$ and $\mu_1\geq 1, \dots, \mu_d\geq 1$
satisfying $\mu_1 \lambda_1 + \cdots + \mu_d \lambda_d \leq \kappa +
1$. By $\mathfrak{ H}_{ \mu_1 \lambda_1 + \cdots + \mathfrak{ H}_{
\mu_d \lambda_d}}$, we denote the subgroup of permutations $\tau \in
\mathfrak{ S}_{\mu_1 \lambda_1 + \cdots + \mathfrak{ H}_{ \mu_d
\lambda_d}}$ that leave unchanged the
general monomial~\thetag{ 3.38}, namely 
that satisfy
\def\theequation{3.67}\begin{equation}
\aligned
&
\prod_{1\leq\nu_1\leq\mu_1}\,
y_{k_{\sigma(1:\nu_1:1)},\dots,k_{\sigma(1:\nu_1:\lambda_1)}}
\cdots \
\prod_{1\leq\nu_d\leq\mu_d}\,
y_{k_{\sigma(d:\nu_d:1)},\dots,k_{\sigma(d:\nu_d:\lambda_d)}}
= \\
&
=
\prod_{1\leq\nu_1\leq\mu_1}\,
y_{k_{1:\nu_1:1},\dots,k_{1:\nu_1:\lambda_1}}
\cdots \
\prod_{1\leq\nu_d\leq\mu_d}\,
y_{k_{d:\nu_d:1},\dots,k_{d:\nu_d:\lambda_d}}.
\endaligned
\end{equation}
The structure of this group may be described as follows. For every $e
= 1, \dots, d$, an arbitrary permutation $\sigma$ of the set
\def\theequation{3.68}\begin{equation}
\{
\underbrace{
\underbrace{
e\!:1\!:\!1, \dots, e\!:1\!:\!\lambda_e}_{\lambda_e},
\underbrace{
e\!:2\!:\!1, \dots, e\!:2\!:\!\lambda_e}_{\lambda_e},
\cdots,
\underbrace{
e\!:\mu_e\!:\!1, \dots, e\!:\mu_e\!:\!\lambda_e}_{\lambda_e}}_{
\mu_e}
\}
\end{equation}
which leaves unchanged the monomial
\def\theequation{3.69}\begin{equation}
\prod_{1\leq\nu_e\leq\mu_e}\,
y_{k_{\sigma(e:\nu_e:1)},\dots,k_{\sigma(e:\nu_e:\lambda_e)}}=
\prod_{1\leq\nu_e\leq\mu_e}\,
y_{k_{e:\nu_e:1},\dots,k_{e:\nu_e:\lambda_e}}.
\end{equation}
uniquely decomposes as the composition of

\smallskip

\begin{itemize}
\item[$\bullet$]
$\mu_e$ arbitrary permutations of the $\mu_e$ groups
of $\lambda_e$ integers 
$\{e\!:\!\nu_e\!:\!1,
\dots,e\!:\!\nu_e\!:\!\lambda_e\}$, 
of total cardinal $(\lambda_e!)^{\mu_e}$;
\item[$\bullet$]
an arbitrary permutation
between these $\mu_e$ groups, of total
cardinal $\mu_e!$.
\end{itemize}

\smallskip

Consequently
\def\theequation{3.70}\begin{equation}
{\rm Card}
\left(
\mathfrak{H}_{\mu_1\lambda_1+\cdots+\mu_d\lambda_d}^{
(\mu_1,\lambda_1),\dots,(\mu_d,\lambda_d)}
\right)
=
\mu_1!(\lambda_1!)^{\mu_1}\cdots
\mu_d!(\lambda_d!)^{\mu_d},
\end{equation}
Finally, define the coset
\def\theequation{3.71}\begin{equation}
\mathfrak{F}_{\mu_1\lambda_1+\cdots+\mu_d\lambda_d}^{
(\mu_1,\lambda_1),\dots,(\mu_d,\lambda_d)}
:=
\mathfrak{S}_{\mu_1\lambda_1+\cdots+\mu_d\lambda_d}
/
\mathfrak{H}_{\mu_1\lambda_1+\cdots+\mu_d\lambda_d}^{
(\mu_1,\lambda_1),\dots,(\mu_d,\lambda_d)}
\end{equation}
with
\def\theequation{3.72}\begin{equation}
\aligned
{\rm Card}
\left(
\mathfrak{F}_{\mu_1\lambda_1+\cdots+\mu_d\lambda_d}^{
(\mu_1,\lambda_1),\dots,(\mu_d,\lambda_d)}
\right)
&
=
\frac{
{\rm Card}
\left(
\mathfrak{S}_{\mu_1\lambda_1+\cdots+\mu_d\lambda_d}
\right)}{
{\rm Card}
\left(
\mathfrak{H}_{\mu_1\lambda_1+\cdots+\mu_d\lambda_d}^{
(\mu_1,\lambda_1),\dots,(\mu_d,\lambda_d)}
\right)}
\\
&
=
\frac{(\mu_1\lambda_1+\cdots+\mu_d\lambda_d)!}{
\mu_1!(\lambda_1!)^{\mu_1}\cdots
\mu_d!(\lambda_d!)^{\mu_d}}.
\endaligned
\end{equation}
In conclusion, by means of this formalism, we may write down the
complete generalization of Theorem~2.24 to several independent
variables.

\def\thetheorem{3.73}\begin{theorem}
For every $\kappa \geq 1$ and for every choice of $\kappa$ indices
$i_1,\dots, i_\kappa$ in the set $\{ 1, 2, \dots, n\}$, the general
expression of ${\bf Y}_{i_1, \dots, i_\kappa}$ is as follows{\rm :}
\def\theequation{3.74}\begin{equation}
\small
\boxed{
\aligned
{\bf Y}_{i_1, \dots, i_\kappa}
&
=
\mathcal{Y}_{x^{i_1}\cdots x^{i_\kappa}}
+
\sum_{d=1}^{\kappa+1}
\ \
\sum_{1\leq\lambda_1<\cdots<\lambda_d\leq\kappa}
\ \
\sum_{\mu_1\geq 1,\dots,\mu_d\geq 1} 
\
\sum_{
\mu_1\lambda_1
+
\cdots
+
\mu_d\lambda_d\leq \kappa+1} 
\\
&
\sum_{k_{1:1:1},\dots,k_{1:1:\lambda_1}=1}^n
\cdots \
\sum_{k_{1:\mu_1:1},\dots,k_{1:\mu_1:\lambda_1}=1}^n
\cdots\cdots \
\sum_{k_{d:1:1},\dots,k_{d:1:\lambda_d}=1}^n
\cdots \
\sum_{k_{d:\mu_d:1},\dots,k_{d:\mu_d:\lambda_d}=1}^n
\\
&
\left[
\aligned
& 
\sum_{\sigma\in\mathfrak{F}_{\mu_1\lambda_1+\cdots+\mu_d\lambda_d}^{
(\mu_1,\lambda_1),\dots,(\mu_d,\lambda_d)}}
\
\sum_{\tau\in\mathfrak{S}_\kappa^{
\mu_1\lambda_1+\cdots+\mu_d\lambda_d}}
\\
& \
\ \ \ \ \
\delta_{i_{\tau(1)},\dots,i_{\tau(\mu_1\lambda_1)},\dots,
i_{\tau(\mu_1\lambda_1+\cdots+\mu_d\lambda_d)}}^{
k_{\sigma(1:1:1)},\dots,k_{\sigma(1:\mu_1:\lambda_1)},
\dots,k_{\sigma(d:\mu_d:\lambda_d)}}\,
\frac{\partial^{\kappa-\mu_1\lambda_1-\cdots-\mu_d\lambda_d+
\mu_1+\cdots+\mu_d}
\mathcal{Y}}{
\partial x^{i_{\tau(\mu_1\lambda_1+\cdots+\mu_d\lambda_d+1)}}\cdots
\partial x^{i_{\tau(\kappa)}}
\left(\partial y\right)^{\mu_1+\cdots+\mu_d}}
- \\
& 
-
\sum_{\sigma\in\mathfrak{F}_{\mu_1\lambda_1+\cdots+\mu_d\lambda_d}^{
(\mu_1,\lambda_1),\dots,(\mu_d,\lambda_d)}}
\
\sum_{\tau\in\mathfrak{S}_\kappa^{
\mu_1\lambda_1+\cdots+\mu_d\lambda_d-1}}
\\
& \
\ \ \ \ \
\delta_{i_{\tau(1)},\dots,i_{\tau(\mu_1\lambda_1)},\dots,
i_{\tau(\mu_1\lambda_1+\cdots+\mu_d\lambda_d-1)}}^{
k_{\sigma(1:1:1)},\dots,k_{\sigma(1:\mu_1:\lambda_1)},
\dots,k_{\sigma(d:\mu_d:\lambda_d-1)}}\,
\frac{\partial^{\kappa-\mu_1\lambda_1-\cdots-\mu_d\lambda_d
+\mu_1+\cdots+\mu_d}\mathcal{X}^{k_{\sigma(d:\mu_d:\lambda_d)}}}{
\partial x^{i_{\tau(\mu_1\lambda_1+\cdots+\mu_d\lambda_d)}}\cdots
\partial x^{i_{\tau(\kappa)}}
\left(\partial y\right)^{\mu_1+\cdots+\mu_d-1}}
\endaligned
\right]
\cdot
\\
&
\ \ \ \ \ \ \ \ \ \ \ \ \ \ \ \ \ \ \ \ \ \ \ \ 
\cdot
\prod_{1\leq\nu_1\leq\mu_1}\,
y_{k_{1:\nu_1:1},\dots,k_{1:\nu_1:\lambda_1}}
\ \cdots \
\prod_{1\leq\nu_d\leq\mu_d}\,
y_{k_{d:\nu_d:1},\dots,k_{d:\nu_d:\lambda_d}}.
\endaligned
}
\end{equation}
\end{theorem}

\subsection*{ 3.75.~Deduction of a multivariate Fa\`a
di Bruno formula} Let $n \in \N$ with $n\geq 1$, let $x = (x^1,\dots,
x^n) \in \K^n$, let $g = g( x^1, \dots, x^n)$ be a $\mathcal{
C}^\infty$-smooth function from $\K^n$ to $\K$, let $y \in \K$ and let
$f = f(y)$ be a $\mathcal{ C}^\infty$ function from $\K$ to $\K$. The
goal is to obtain an explicit formula for the partial derivatives of
the composition $h := f\circ g$, namely $h(x^1, \dots, x^n) := f (
g(x^1,\dots, x^n))$. For $\lambda \in \N$ with $\lambda \geq 1$ and
for arbitrary indices $i_1, \dots, i_\lambda = 1, \dots, n$, we
shall abbreviate the partial derivative $\frac{ \partial^\lambda g}{
\partial x^{i_1} \cdots \partial x^{i_\lambda}}$ by $g_{i_1,\dots,
i_\lambda}$ and similarly for $h_{i_1, \dots, i_\lambda}$. The
derivative $\frac{ d^\lambda f}{ d y^\lambda}$ will be
abbreviated by $f_\lambda$.

Appying the chain rule, we may compute:
\def\theequation{3.76}\begin{equation}
\small
\aligned
h_{i_1}
&
=
f_1 
\left[
g_{i_1}
\right], 
\\
h_{i_1,i_2}
&
=
f_2
\left[
g_{i_1}\,g_{i_2}
\right]
+
f_1
\left[
g_{i_1,i_2}
\right],
\\
h_{i_1,i_2,i_3}
&
=
f_3
\left[
g_{i_1}\,g_{i_2}\,g_{i_3}
\right]
+
f_2
\left[
g_{i_1}\,g_{i_2,i_3}
+
g_{i_2}\,g_{i_1,i_3}
+
g_{i_3}\,g_{i_1,i_2}
\right]
+
f_1
\left[
g_{i_1,i_2,i_3}
\right],
\\
h_{i_1,i_2,i_3,i_4}
&
=
f_4
\left[
g_{i_1}\,g_{i_2}\,g_{i_3}\,g_{i_4}
\right]
+
f_3
\left[
g_{i_2}\,g_{i_3}\,g_{i_1,i_4}
+
g_{i_3}\,g_{i_1}\,g_{i_2,i_4}
+
g_{i_1}\,g_{i_2}\,g_{i_3,i_4}
+
\right. 
\\
& \
\ \ \ \ \ \ \ \ \ \ \ \ \ \ \ \ \ \ \ \
\ \ \ \ \ \ \ \ \ \ \ \ \ \ \ \ \ \ \ \
\ \
\left.
+
g_{i_1}\,g_{i_4}\,g_{i_2,i_3}
+
g_{i_2}\,g_{i_4}\,g_{i_1,i_3}
+
g_{i_3}\,g_{i_4}\,g_{i_1,i_2}
\right]
+ \\
& \
\ \ \ \ \
+
f_2
\left[
g_{i_1,i_2}\,g_{i_3,i_4}
+
g_{i_1,i_3}\,g_{i_2,i_4}
+
g_{i_1,i_4}\,g_{i_2,i_3}
\right]
+ \\
& \
\ \ \ \ \
+
f_2
\left[
g_{i_1}\,g_{i_2,i_3,i_4}
+
g_{i_2}\,g_{i_1,i_3,i_4}
+
g_{i_3}\,g_{i_1,i_2,i_4}
+
g_{i_4}\,g_{i_1,i_2,i_3}
\right]
+ \\
& \
\ \ \ \ \
+
f_1
\left[
g_{i_1,i_2,i_3,i_4}
\right].
\endaligned
\end{equation}
Introducing the derivations
\def\theequation{3.77}\begin{equation}
\small
\aligned
F_i^2
&
:= 
\sum_{k_1=1}^n\,g_{k_1,i}\,
\frac{\partial}{\partial g_{k_1}}
+
g_i\left(
f_2\,\frac{\partial}{\partial f_1}
\right), \\
F_i^3
&
:=
\sum_{k_1=1}^n\,g_{k_1,i}\,
\frac{\partial}{\partial g_{k_1}}
+
\sum_{k_1,k_2=1}^n\,g_{k_1,k_2,i}\,
\frac{\partial}{\partial g_{k_1,k_2}}
+
g_i\left(
f_2\,\frac{\partial}{\partial f_1}
+
f_3\,\frac{\partial}{\partial f_2}
\right), \\
\text{\bf 
\dots\dots}
& \
\ \ \ \ \
\text{\bf 
\dots\dots\dots\dots\dots\dots\dots\dots\dots
\dots\dots\dots\dots\dots\dots\dots\dots\dots
\dots\dots\dots\dots\dots\dots\dots
} \\
\endaligned
\end{equation}
$$
\small
\aligned
F_i^\lambda
&
:=
\sum_{k_1=1}^n\,g_{k_1,i}\,
\frac{\partial}{\partial g_{k_1}}
+
\sum_{k_1,k_2=1}^n\,g_{k_1,k_2,i}\,
\frac{\partial}{\partial g_{k_1,k_2}}
+
\cdots
+ \\
& \
\ \ \ \ \ \ \ \ \ \ 
+
\sum_{k_1,\dots,k_{\lambda-1}=1}^n\,g_{k_1,\dots,k_{\lambda-1},i}\,
\frac{\partial}{\partial g_{k_1,\dots,k_{\lambda-1}}}
+ \\
& \
\ \ \ \ \ \ \ \ \ \ 
+
g_i\left(
f_2\,\frac{\partial}{\partial f_1}
+
f_3\,\frac{\partial}{\partial f_2}
+
\cdots
+
f_\lambda\,\frac{\partial}{\partial f_{\lambda-1}}
\right), \\
\endaligned
$$
we observe that the following
induction relations hold:
\def\theequation{3.78}\begin{equation}
\aligned
h_{i_1,i_2}
&
=
F_{i_2}^2
\left(
h_{i_1}
\right), 
\\
h_{i_1,i_2,i_3}
&
=
F_{i_3}^3
\left(
h_{i_1,i_2}
\right), 
\\
\text{\bf 
\dots\dots
}
& \
\ \ \ \ \
\text{\bf 
\dots\dots\dots\dots\dots
}
\\
h_{i_1,i_2,\dots,i_\lambda}
&
=
F_{i_\lambda}^\lambda
\left(
h_{i_1,i_2,\dots,i_{\lambda-1}}
\right).
\endaligned
\end{equation}

To obtain the explicit version of the Fa\`a di Bruno in the case of
several variables $(x^1, \dots, x^n)$ and one variable $y$, it
suffices to extract from the expression of ${\bf Y }_{ i_1,\dots,
i_\kappa}$ provided by Theorem~3.73 only the terms corresponding to
$\mu_1 \lambda_1 + \cdots + \mu_d\lambda_d = \kappa$, dropping all the
$\mathcal{ X}$ terms. After some simplifications and after a
translation by means of an elementary dictionary, we obtain a
statement.

\def\thetheorem{3.79}\begin{theorem} 
For every integer $\kappa \geq 1$ and for every choice of indices
$i_1, \dots, i_\kappa$ in the set $\{ 1, 2, \dots, n\}$, the
$\kappa$-th partial derivative of the composite function $h = h( x^1,
\dots, x^n) = f( g(x^1, \dots, x^n))$ with respect to the variables
$x^{i_1}, \dots, x^{i_\kappa}$ may be expressed as an explicit
polynomial depending on the derivatives of $f$, on the partial
derivatives of $g$ and having integer coefficients{\rm :}
\def\theequation{3.80}\begin{equation}
\boxed{
\aligned
\frac{\partial^\kappa h}{\partial x^{i_1}\cdots
\partial x^{i_\kappa}}
&
=
\sum_{d=1}^\kappa
\
\sum_{1\leq \lambda_1 < \cdots < \lambda_d \leq \kappa}
\
\sum_{\mu_1\geq 1,\dots,\mu_d\geq 1}
\
\sum_{\mu_1\lambda_1+\cdots+\mu_d\lambda_d=\kappa}
\
\frac{d^{\mu_1+\cdots+\mu_d} f}{
dy^{\mu_1+\cdots+\mu_d}}
\\
& 
\left[
\aligned
&
\sum_{\sigma\in\mathfrak{F}_\kappa^{
(\mu_1,\lambda_1),\dots,(\mu_d,\lambda_d)}}
\
\prod_{1\leq\nu_1\leq\mu_1}
\
\frac{\partial^{\lambda_1} g}{\partial
x^{i_{\sigma(1:\nu_1:1)}}\cdots
\partial x^{i_{\sigma(1:\nu_1:\lambda_1)}}}
\
\text{\bf \dots} 
\\
& \
\ \ \ \ \ \ \ \ \ \ \ \ \ \ \ \ \ \ \ \ \ \ \ \ \ \ \ \ \ \
\text{\bf \dots} 
\prod_{1\leq\nu_d\leq\mu_d}
\
\frac{\partial^{\lambda_d} g}{\partial
x^{i_{\sigma(d:\nu_d:1)}}\cdots
\partial x^{i_{\sigma(d:\nu_ds:\lambda_d)}}}
\endaligned
\right].
\endaligned
}
\end{equation}
\end{theorem}

In this formula, the coset $\mathfrak{ F }_\kappa^{
(\mu_1, \lambda_1 ),\dots, ( \mu_d, \lambda_d)}$
was defined in equation~\thetag{ 3.71}; we have made the identification{\rm :}
\def\theequation{3.81}\begin{equation}
\{1,\dots,\kappa\}
\equiv
\{
1\!:\!1\!:\!1,
\dots,
1\!:\!\mu_1\!:\!\lambda_1,
\text{\rm \dots\dots}, 
d\!:\!1\!:\!1,
\dots,
d\!:\!\mu_d\!:\!\lambda_d
\};
\end{equation}
and also, for the sake of clarity, we have restituted the complete
(not abbreviated) notation for the (partial) derivatives of $f$ and of
$g$. 

\section*{\S4.~Several independent variables and
one dependent variable}

\subsection*{4.1.~Simplified adapted notations}
Assume $n = 1$ and $m \geq 1$, let $\kappa\in \N$ with $\kappa \geq
1$ and simply denote the jet variables by
(instead of~\thetag{ 1.2}):
\def\theequation{4.2}\begin{equation}
\left(
x,y^j,y_1^j,y_2^j,\dots,y_\kappa^j
\right)
\in
\mathcal{J}_{1,m}^\kappa.
\end{equation}
Instead of~\thetag{ 1.30}, denote the $\kappa$-th prolongation
of a vector field by:
\def\theequation{4.3}\begin{equation}
\left\{
\aligned
\mathcal{L}^{(\kappa)}
&
=
\mathcal{X}\,\frac{\partial}{\partial x}
+
\sum_{j=1}^m\,\mathcal{Y}^j\,\frac{\partial}{\partial y^j}
+
\sum_{j=1}^m\,{\bf Y}_1^j\,\frac{\partial}{\partial y_1^j}
+
\sum_{j=1}^m\,{\bf Y}_2^j\,\frac{\partial}{\partial y_2^j}
+ \\
& \
\ \ \ \ \
+
\cdots
+
\sum_{j=1}^m\,{\bf Y}_\kappa^j\,\frac{\partial}{\partial y_\kappa^j}.
\endaligned\right.
\end{equation}
The induction formulas are:
\def\theequation{4.4}\begin{equation}
\left\{
\aligned
{\bf Y}_1^j
&
:=
D^1
\left(
\mathcal{Y}^j
\right)
-
D^1
\left(
\mathcal{ X}
\right)
y_1^j, 
\\
{\bf Y}_2^j
&
:=
D^2
\left(
{\bf Y}_1^j
\right)
-
D^1
\left(
\mathcal{ X}
\right)
y_2^j, 
\\
\cdots\cdots
&
\cdots\cdots
\\
{\bf Y}_\lambda^j
&
:=
D^\lambda
\left(
{\bf Y}_{\lambda-1}^j
\right)
-
D^1
\left(
\mathcal{ X}
\right)
y_\lambda^j, 
\endaligned\right.
\end{equation}
where the total
differentiation operators $D^\lambda$
are denoted by (instead of~\thetag{ 1.22}):
\def\theequation{4.5}\begin{equation}
\aligned
D^\lambda
:=
\frac{\partial}{\partial x}
+
\sum_{l=1}^m\,y_1^l\,
\frac{\partial}{\partial y^l}
+
\sum_{l=1}^m\,y_2^l\,
\frac{\partial}{\partial y_1^l}
+
\cdots
+ 
\sum_{l=1}^m\,y_\lambda^l\,
\frac{\partial}{\partial y_{\lambda-1}^l}.
\endaligned
\end{equation}
Applying the definitions in the first two lines of~\thetag{ 4.4}, we
compute, we simplify and we organize the results in a harmonious way,
using in an essential way the Kronecker symbol. Here, the
computations are more elementary than the computations of ${\bf
Y}_{i_1}$ and of ${\bf Y}_{i_1, i_2}$ achieved thoroughly 
in the previous
Section~3, so that we do not provide a Latex track of the
details. Firstly and secondly:
\def\theequation{4.6}\begin{equation}
\small
\left\{
\aligned
{\bf Y}_1^j
&
=
\mathcal{ Y}_x^j
+
\sum_{l_1=1}^m
\left[
\mathcal{Y}_{y^{l_1}}^j
-
\delta_{l_1}^j\,
\mathcal{X}_x
\right]
y_1^{l_1}
+
\sum_{l_1,l_2=1}^m
\left[
-
\delta_{l_1}^j\,
\mathcal{X}_{y^{l_2}}
\right]
y_1^{l_1}y_1^{l_2}, \\
{\bf Y}_2^j
&
=
\mathcal{ Y}_{x^2}^j
+
\sum_{l_1=1}^m
\left[
2\,\mathcal{Y}_{xy^{l_1}}^j
-
\delta_{l_1}^j\,
\mathcal{X}_{x^2}
\right]
y_1^{l_1}
+
\sum_{l_1,l_2=1}^m
\left[
\mathcal{Y}_{y^{l_1}y^{l_2}}^j
-
\delta_{l_1}^j\,2\,
\mathcal{X}_{xy^{l_2}}
\right]
y_1^{l_1}y_1^{l_2}
+ \\
& \
\ \ \ \ \
+
\sum_{l_1,l_2,l_3}
\left[
-
\delta_{l_1}^j\,
\mathcal{X}_{y^{l_2}y^{l_3}}
\right]
y_1^{l_1}y_1^{l_2}y_1^{l_3}
+
\sum_{l_1}
\left[
\mathcal{Y}_{y^{l_1}}^j
-
\delta_{l_1}^j\,2\,
\mathcal{X}_x
\right]
y_2^{l_1}
+ \\
& \
\ \ \ \ \ 
+
\sum_{l_1,l_2=1}^m
\left[
-
\delta_{l_1}^j\,
\mathcal{X}_{y^{l_2}}
-
\delta_{l_2}^j\,2\,
\mathcal{X}_{y^{l_1}}
\right]
y_1^{l_1}y_2^{l_2}.
\endaligned\right.
\end{equation}
Thirdly:
\def\theequation{4.7}\begin{equation}
\small
\aligned
{\bf Y}_3^j
&
=
\mathcal{ Y}_{x^3}^j
+
\sum_{l_1=1}^m
\left[
3\,\mathcal{Y}_{x^2y^{l_1}}^j
-
\delta_{l_1}^j\,
\mathcal{X}_{x^3}
\right]
y_1^{l_1}
+
\sum_{l_1,l_2=1}^m
\left[
3\,\mathcal{Y}_{xy^{l_1}y^{l_2}}^j
-
\delta_{l_1}^j\,
3\,
\mathcal{X}_{x^2y^{l_2}}
\right]
y_1^{l_1}y_1^{l_2}+ \\
& \
\ \ \ \ \
+
\sum_{l_1,l_2,l_3}
\left[
\mathcal{Y}_{y^{l_1}y^{l_2}y^{l_3}}^j
-
\delta_{l_1}^j\,
3\,\mathcal{X}_{xy^{l_2}y^{l_3}}
\right]
y_1^{l_1}y_1^{l_2}y_1^{l_3}
+ \\
& \ 
\ \ \ \ \
+
\sum_{l_1,l_2,l_3,l_4}
\left[
-
\delta_{l_1}^j\,
\mathcal{X}_{y^{l_2}y^{l_3}y^{l_4}}
\right]
y_1^{l_1}y_1^{l_2}y_1^{l_3}y_1^{l_4}
+
\sum_{l_1=1}^m
\left[
3\,\mathcal{Y}_{xy^{l_1}}^j
-
\delta_{l_1}^j\,
3\,\mathcal{X}_{x^2}
\right]
y_2^{l_1} 
\endaligned
\end{equation}
$$
\small
\aligned
+ \\
& \
\ \ \ \ \ 
+
\sum_{l_1,l_2=1}^m
\left[
3\,\mathcal{Y}_{y^{l_1}y^{l_2}}^j
-
\delta_{l_1}^j\,
3\,\mathcal{X}_{xy^{l_2}}
-
\delta_{l_2}^j\,
6\,\mathcal{X}_{xy^{l_1}}
\right]
y_1^{l_1}y_2^{l_2}
+ \\
& \
\ \ \ \ \
+
\sum_{l_1,l_2,l_3=1}^m
\left[
-
\delta_{l_1}^j\,
3\,\mathcal{X}_{y^{l_2}y^{l_3}}
-
\delta_{l_3}^j\,
3\,\mathcal{X}_{y^{l_1}y^{l_2}}
\right]
y_1^{l_1}y_1^{l_2}y_2^{l_3}
+ 
\sum_{l_1,l_2=1}^m
\left[
-
\delta_{l_3}^j\,
3\,\mathcal{X}_{y^{l_2}}
\right]
y_2^{l_1}y_2^{l_2}
+ \\
& \
\ \ \ \ \
+
\sum_{l_1=1}^m
\left[
\mathcal{Y}_{y^{l_1}}^j
-
\delta_{l_1}^j\,
3\,\mathcal{X}_x
\right]
y_3^{l_1}
+
\sum_{l_1,l_2=1}^m
\left[
-
\delta_{l_1}^j\,
\mathcal{X}_{y^{l_2}}
-
\delta_{l_2}^j\,
3\,\mathcal{X}_{y^{l_1}}
\right]
y_1^{l_1}y_3^{l_2}.
\endaligned
$$
Fourthly:
$$
\small
\aligned
{\bf Y}_4^j
&
=
\mathcal{ Y}_{x^4}^j
+
\sum_{l_1=1}^m
\left[
4\,\mathcal{Y}_{x^3y^{l_1}}^j
-
\delta_{l_1}^j
\mathcal{X}_{x^4}\,
\right]
y_1^{l_1}
+
\sum_{l_1,l_2=1}^m
\left[
6\,\mathcal{Y}_{x^2y^{l_1}y^{l_2}}^j
-
\delta_{l_1}^j\,
4\,
\mathcal{X}_{x^3y^{l_2}}
\right]
y_1^{l_1}y_1^{l_2}
+ \\
& \
\ \ \ \ \
+
\sum_{l_1,l_2,l_3=1}^m
\left[
4\,\mathcal{Y}_{xy^{l_1}y^{l_2}y^{l_3}}^j
-
\delta_{l_1}^j\,
6\,\mathcal{X}_{x^2y^{l_2}y^{l_3}}
\right]
y_1^{l_1}y_1^{l_2}y_1^{l_3}
+ \\
& \
\ \ \ \ \
+
\sum_{l_1,l_2,l_3,l_4=1}^m
\left[
\mathcal{Y}_{xy^{l_1}y^{l_2}y^{l_3}y^{l_4}}^j
-
\delta_{l_1}^j\,
4\,\mathcal{X}_{xy^{l_2}y^{l_3}y^{l_4}}
\right]
y_1^{l_1}y_1^{l_2}y_1^{l_3}y_1^{l_4}
+ \\
& \ 
\ \ \ \ \
+
\sum_{l_1,l_2,l_3,l_4,l_5=1}^m
\left[
-
\delta_{l_1}^j\,
\mathcal{X}_{y^{l_2}y^{l_3}y^{l_4}y^{l_5}}
\right]
y_1^{l_1}y_1^{l_2}y_1^{l_3}y_1^{l_4}y_1^{l_5}
+
\sum_{l_1=1}^m
\left[
6\,\mathcal{Y}_{x^2y^{l_1}}^j
-
\delta_{l_1}^j\,
4\,\mathcal{X}_{x^3}
\right]
y_2^{l_1}
+
\endaligned
$$
\def\theequation{4.8}\begin{equation}
\small
\aligned
& \
\ \ \ \ \ 
+
\sum_{l_1,l_2=1}^m
\left[
12\,\mathcal{Y}_{xy^{l_1}y^{l_2}}^j
-
\delta_{l_1}^j\,
6\,\mathcal{X}_{x^2y^{l_2}}
-
\delta_{l_2}^j\,
12\,\mathcal{X}_{x^2y^{l_1}}
\right]
y_1^{l_1}y_2^{l_2}
+ \\
& \
\ \ \ \ \ 
+
\sum_{l_1,l_2,l_3=1}^m
\left[
6\,\mathcal{Y}_{y^{l_1}y^{l_2}y^{l_3}}^j
-
\delta_{l_1}^j\,
12\,\mathcal{X}_{xy^{l_2}y^{l_3}}
-
\delta_{l_3}^j\,
12\,\mathcal{X}_{xy^{l_1}y^{l_2}}
\right]
y_1^{l_1}y_1^{l_2}y_2^{l_3}
+ \\
& \
\ \ \ \ \ 
+
\sum_{l_1,l_2,l_3,l_4=1}^m
\left[
-
\delta_{l_1}^j\,
6\,\mathcal{X}_{y^{l_2}y^{l_3}y^{l_4}}
-
\delta_{l_4}^j\,
4\,\mathcal{X}_{y^{l_1}y^{l_2}y^{l_3}}
\right]
y_1^{l_1}y_1^{l_2}y_1^{l_3}y_2^{l_4}
+ \\
\endaligned
\end{equation}
$$
\small
\aligned
& \
\ \ \ \ \
+
\sum_{l_1,l_2=1}^m
\left[
3\,\mathcal{Y}_{y^{l_1}y^{l_2}}^j
-
\delta_{l_1}^j\,
12\,\mathcal{X}_{xy^{l_2}}
\right]
y_2^{l_1}y_2^{l_2}
+ \\
& \
\ \ \ \ \
+
\sum_{l_1,l_2,l_3=1}^m
\left[
-
\delta_{l_1}^j\,
3\,\mathcal{X}_{y^{l_2}y^{l_3}}
-
\delta_{l_2}^j\,
12\,\mathcal{X}_{y^{l_1}y^{l_3}}
\right]
y_1^{l_1}y_2^{l_2}y_2^{l_3}
+
\sum_{l_1=1}^m
\left[
4\,\mathcal{Y}_{xy^{l_1}}^j
-
\delta_{l_1}^j\,
6\,\mathcal{X}_{x^2}
\right]
y_3^{l_1}
+ \\
& \
\ \ \ \ \
+
\sum_{l_1,l_2=1}^m
\left[
4\,\mathcal{Y}_{y^{l_1}y^{l_2}}^j
-
\delta_{l_1}^j\,
4\,\mathcal{X}_{xy^{l_2}}
-
\delta_{l_2}^j\,
12\,\mathcal{X}_{xy^{l_1}}
\right]
y_1^{l_1}y_3^{l_2}
+ \\
& \
\ \ \ \ \
+
\sum_{l_1,l_2,l_3=1}^m
\left[
-
\delta_{l_1}^j\,
4\,\mathcal{X}_{y^{l_2}y^{l_3}}
-
\delta_{l_3}^j\,
6\,\mathcal{X}_{y^{l_1}y^{l_2}}
\right]
y_1^{l_1}y_1^{l_2}y_3^{l_3}
+ \\
& \
\ \ \ \ \
+
\sum_{l_1,l_2=1}^m
\left[
-
\delta_{l_1}^j\,
4\,\mathcal{X}_{y^{l_2}}
-
\delta_{l_2}^j\,
6\,\mathcal{X}_{y^{l_1}}
\right]
y_2^{l_1}y_3^{l_2}
+ \\
& \
\ \ \ \ \
+
\sum_{l_1=1}^m
\left[
\mathcal{Y}_{y^{l_1}}^j
-
\delta_{l_1}^j\,
4\,\mathcal{X}_x
\right]
y_4^{l_1}
+
\sum_{l_1,l_2=1}^m
\left[
-
\delta_{l_1}^j\,
\mathcal{X}_{y^{l_2}}
-
\delta_{l_2}^j\,
4\,\mathcal{X}_{y^{l_1}}
\right]
y_1^{l_1}y_4^{l_2}.
\endaligned
$$

\subsection*{4.9.~Inductive elaboration of the
general formula} Now we compare the formula~\thetag{ 2.9} for ${\bf
Y}_4$ with the above formula~\thetag{ 4.8} for ${\bf Y}_4^j$. The goal
is to find the rules of transformation and of development by
inspecting several instances, in order to devise how to transform and
to develope the formula~\thetag{ 2.25} to several dependent variables.

First of all, we have to develope the general monomial $(y_{ \lambda_1
})^{ \mu_1} \cdots (y_{ \lambda_d })^{ \mu_d}$. In every monomial
present in the expressions of ${\bf Y}_1^j$, of ${\bf Y}_2^j$, of
${\bf Y}_3^j$ and of ${\bf Y}_4^j$ above, we see that the number
$\alpha$ of indices $l_\beta$ appearing in all the sums $\sum_{l_1,
\dots, l_\alpha = 1}^m$ is exactly equal to $\mu_1 + \dots+ \mu_d$. To
denote these $\mu_1 + \cdots + \mu_d$ indices $l_\beta$, we shall use
the notation:
\def\theequation{4.10}\begin{equation}
\underbrace{
\underbrace{
l_{1:1},\dots,l_{1:\mu_1}}_{\mu_1},
\dots,
\underbrace{
l_{d:1},\dots,l_{d:\mu_d}}_{\mu_d}}_{
\mu_1+\cdots+\mu_d},
\end{equation} 
inspired by Convention~3.33. With such a choice of notation, we may
avoid long subscripts in the indices $l_\beta$, like $l_{\mu_1+\cdots+
\mu_d}$. It follows that the development of the general monomial $(y_{
\lambda_1 })^{ \mu_1} \cdots (y_{ \lambda_d })^{ \mu_d}$ to several
dependent variables yields $m^{\mu_1+ \cdots + \mu_d}$ possible
choices:
\def\theequation{4.11}\begin{equation}
\prod_{1\leq\nu_1\leq\mu_1}
y_{\lambda_1}^{l_{1:\nu_1}}
\ \cdots\cdots
\prod_{1\leq\nu_d\leq\mu_d}
y_{\lambda_d}^{l_{d:\nu_d}},
\end{equation}
where the indices $l_{1 :1}, \dots, l_{1: \mu_1 }, \dots, l_{ d:1},
\dots, l_{d: \mu_d }$ take their values in the set $\{ 1, 2, \dots, m
\}$. Consequently, the general expression of ${\bf Y }_\kappa^j$ must
be of the form:
\def\theequation{4.12}\begin{equation}
\small
\aligned
{\bf Y}_\kappa^j
&
=
\mathcal{Y}_{x^\kappa}^j
+
\sum_{d=1}^{\kappa+1}
\
\sum_{1\leq\lambda_1<\cdots<\lambda_d\leq\kappa}
\
\sum_{\mu_1\geq 1,\dots,\mu_d\geq 1}
\
\sum_{\mu_1\lambda_1+\cdots+\mu_d\lambda_d\leq\kappa+1} 
\
\\
& \
\ \ \ \ \ \ \ \ \ \ \ \ \ \ \
\sum_{l_{1:1}=1}^m
\cdots
\sum_{l_{1:\mu_1}=1}^m
\cdots\cdots
\sum_{l_{d:1}=1}^m
\cdots
\sum_{l_{d:\mu_d}=1}^m
\
\text{\bf [?]}
\\
& \
\ \ \ \ \ \ \ \ \ \ \ \ \ \ \ \ \ \ \
\prod_{1\leq\nu_1\leq\mu_1}
y_{\lambda_1}^{l_{1:\nu_1}}
\ \cdots\cdots
\prod_{1\leq\nu_d\leq\mu_d}
y_{\lambda_d}^{l_{d:\nu_d}},
\endaligned
\end{equation}
where the term in brackets {\bf [?]} is still unknown. To determine
it, let us examine a few instances.

According to~\thetag{ 4.8} (fourth line), the term $\left[ 6\,
\mathcal{ Y}_{ x^2y} - 4\, \mathcal{ X}_{ x^3} \right] y_2$ of ${\bf
Y}_4$ developes as $\sum_{ l_1 =1}^m \, \left[ 6\, \mathcal{ Y}_{ x^2
y^{l_1}}^j - \delta_{ l_1}^j \, 4 \, \mathcal{ X}_{ x^3} \right] y_2^{
l_1}$ in ${\bf Y}_4^j$. Here, 
$6\, \mathcal{ Y}_{ x^2y}$ just becomes
$6\, \mathcal{ Y}_{ x^2y^{l_1}}^j$. Thus, we suspect that the term
$\frac{\kappa\cdots(\kappa-\mu_1\lambda_1-\cdots-\mu_d\lambda_d+1)}
{(\lambda_1!)^{\mu_1}\,\mu_1! \cdots (\lambda_d! )^{\mu_d}\,\mu_d! }
\cdot \mathcal{ Y}_{ x^{ \kappa-\mu_1\lambda_1-\cdots-\mu_d\lambda_d}
\, y^{\mu_1+\cdots+\mu_d} } $ of the second line of~\thetag{ 2.25}
should simply be developed as
\def\theequation{4.13}\begin{equation}
\small
\aligned
& 
\frac{\kappa(\kappa-1)\cdots(\kappa-\mu_1\lambda_1-\cdots
-\mu_d\lambda_d+1)}{
(\lambda_1!)^{\mu_1}\ \mu_1!\cdots
(\lambda_d!)^{\mu_d}\ \mu_d!}
\cdot
\\
& \
\ \ \ \ \ \ \ \
\cdot
\frac{\partial^{\kappa-\mu_1\lambda_1-\cdots
-\mu_d\lambda_d+\mu_1+\cdots+\mu_d}\mathcal{Y}^j}{
(\partial x)^{\kappa-\mu_1\lambda_1-\cdots
-\mu_d\lambda_d}
\partial y^{l_{1:1}}
\cdots
\partial y^{l_{1:\mu_1}}
\cdots
\partial y^{l_{d:1}}
\cdots
\partial y^{l_{d:\mu_d}}
}.
\endaligned
\end{equation}
This rule is confirmed by inspecting all the other monomials of ${\bf
Y}_1^j$, of ${\bf Y}_2^j$, of ${\bf Y}_3^j$ and of ${\bf Y}_4^j$.

It remains to determine how we must develope the term in $\mathcal{
X}$ appearing in the last two lines of~\thetag{ 2.25}. To begin with,
let us rewrite in advance this term in the slightly different
shape, emphasizing a factorization:
\def\theequation{4.14}\begin{equation}
\small
\aligned
\frac{\kappa\cdots(
\kappa-\mu_1\lambda_1-\cdots-\mu_d\lambda_d+2)}
{(\lambda_1!)^{\mu_1}\,\mu_1!
\cdots
(\lambda_d!)^{\mu_d}\,\mu_d!
}
\left[
(\mu_1\lambda_1+\cdots+\mu_d\lambda_d)
\mathcal{X}_{
x^{\kappa-\mu_1\lambda_1-\cdots-\mu_d\lambda_d+1}
\,
y^{\mu_1+\cdots+\mu_d-1}
}
\right].
\endaligned
\end{equation}
Then we examine four instances extracted from the complete expression
of ${\bf Y}_4^j$:
\def\theequation{4.15}\begin{equation}
\small
\left\{
\aligned
&
\sum_{l_1,l_2,l_3=1}^m
\left[
4\,\mathcal{Y}_{xy^{l_1}y^{l_2}y^{l_3}}^j
-
\delta_{l_1}^j\,
6\,\mathcal{X}_{x^2y^{l_2}y^{l_3}}
\right]
y_1^{l_1}y_1^{l_2}y_1^{l_3}, 
\\
&
\sum_{l_1,l_2=1}^m
\left[
12\,\mathcal{Y}_{xy^{l_1}y^{l_2}}^j
-
\delta_{l_1}^j\,
6\,\mathcal{X}_{x^2y^{l_2}}
-
\delta_{l_2}^j\,
12\,\mathcal{X}_{x^2y^{l_1}}
\right]
y_1^{l_1}y_2^{l_2},
\\
&
\sum_{l_1,l_2,l_3,l_4=1}^m
\left[
-
\delta_{l_1}^j\,
6\,\mathcal{X}_{y^{l_2}y^{l_3}y^{l_4}}
-
\delta_{l_4}^j\,
4\,\mathcal{X}_{y^{l_1}y^{l_2}y^{l_3}}
\right]
y_1^{l_1}y_1^{l_2}y_1^{l_3}y_2^{l_4},
\\
&
\sum_{l_1,l_2,l_3=1}^m
\left[
-
\delta_{l_1}^j\,
4\,\mathcal{X}_{y^{l_2}y^{l_3}}
-
\delta_{l_3}^j\,
6\,\mathcal{X}_{y^{l_1}y^{l_2}}
\right]
y_1^{l_1}y_1^{l_2}y_3^{l_3},
\endaligned\right.
\end{equation}
and we compare them to the corresponding terms of ${\bf Y}_4$:
\def\theequation{4.16}\begin{equation}
\small
\left\{
\aligned
&
\left[
4\,\mathcal{Y}_{xy^3}
-
6\,\mathcal{X}_{x^2y^2}
\right]
(y_1)^3, 
\\
&
\left[
12\,\mathcal{Y}_{xy^2}
-
18\,\mathcal{X}_{x^2y}
\right]
y_1y_2,
\\
&
\left[
-
10\,\mathcal{X}_{y^3}
\right]
(y_1)^3y_2,
\\
&
\left[
-
10\,\mathcal{X}_{y^2}
\right]
(y_1)^2y_3.
\endaligned\right.
\end{equation}
In the development from~\thetag{ 4.16} to~\thetag{ 4.15}, we see that
the four integers just before $\mathcal{ X}$, namely $6 = 6$, $18 = 6
+ 12$, $10 = 6 + 4$ and $10 = 4 + 6$, are split in a certain
manner. Also, a single Kronecker symbol $\delta_{l_\alpha}^j$ is added
as a factor. {\it What are the rules}?

In the second splitting $18 = 6 + 12$, we see that the relative weight
of $6$ and of $12$ is the same as the relative weight of $1$ and $2$
in the splitting $3 = 1 + 2$ issued from the lower indices of the
corresponding monomial $y_1^{l_1} y_2^{l_2}$. Similarly, in the third
splitting $10 = 6 + 4$, the relative weight of $6$ and of $4$ is the
same as the relative weight of $1+1+1$ and of $2$ issued from the
lower indices of the corresponding monomial $y_1^{l_1} y_1^{l_2}
y_1^{l_3} y_2^{l_4}$. This rule may be confirmed by inspecting all
the other monomials of ${\bf Y}_2$, ${\bf Y}_2^j$, of ${\bf Y}_3$,
${\bf Y}_3^j$ and of ${\bf Y}_4$, ${\bf Y}_4^j$. For a general $\kappa
\geq 1$, the splitting of integers just amounts to decompose the sum
appearing inside the brackets of~\thetag{ 4.14} as $\mu_1\lambda_1,
\mu_2\lambda_2, \dots, \mu_d\lambda_d$. In fact, when we 
wrote~\thetag{ 4.14}, we emphasized in advance the decomposable
factor $(\mu_1 \lambda_1 + \cdots + \mu_d \lambda_d)$.

Next, we have to determine what is the subscript $\alpha$ in the
Kronecker symbol $\delta_{l_\alpha}^j$. We claim that in the four
instances~\thetag{ 4.15}, the subscript $\alpha$ is intrinsically
related to weight splitting. Indeed, recall that in the second line
of~\thetag{ 4.15}, the number $6$ of the splitting $18 = 6 + 12$ is
related to the number $1$ in the splitting $3 = 1 + 2$ of the lower
indices of the monomial $y_1^{l_1} y_2^{l_2}$. It follows that the
index $l_\alpha$ {\it must be}\, the index $l_1$ of the monomial
$y_1^{l_1}$. Similarly, also in the second line of~\thetag{ 4.15}, the
number $12$ of the splitting $18 = 6 + 12$ being related to the number
$2$ in the splitting $3 = 1 + 2$ of the lower indices of the monomial
$y_1^{l_1} y_2^{l_2}$, it follows that the index $l_\alpha$ attached
to the second $\mathcal{ X}$ term must be the index $l_2$ of the
monomial $y_2^{l_2}$.

This rule is still ambiguous. Indeed, let us examine the third line
of~\thetag{ 4.15}. We have the splitting $10 = 6 + 4$, homologous to
the splitting of relative weights $5 = (1+1+1) + 2$ in the monomial
$y_1^{ l_1} y_1^{ l_2} y_1^{ l_3} y_2^{ l_4}$. Of course, it is clear
that we must choose the index $l_4$ for the Kronecker symbol
associated to the second $\mathcal{ X}$ term $-4\, \mathcal{ X}_{
y^3}$, thus obtaining $-\delta_{ l_4}^j \, 4 \, \mathcal{ X}_{ y^{
l_1} y^{ l_2} y^{ l_3}}$. However, since the monomial $y_1^{ l_1}
y_1^{ l_2} y_1^{ l_3}$ has three indices $l_1$, $l_2$ and $l_3$, there
arises a question: {\it what index $l_\alpha$ must we choose for the
Kronecker symbol $\delta_{ l_\alpha}^j$ attached to the first
$\mathcal{ X}$ term $6\,\mathcal{ X}_{y^3}${\rm :} the index $l_1$,
the index $l_2$ or the index $l_3$}?

The answer is simple: {\it any of the three indices $l_1$, $l_2$ or
$l_3$ works}. Indeed, since the monomial $y_1^{ l_1} y_1^{ l_2}
y_1^{ l_3}$ is symmetric with respect to all permutations of the set
of three indices $\{ l_1, l_2, l_3\}$, we have
\def\theequation{4.17}\begin{equation}
\small
\aligned
\sum_{ l_1, l_2, l_3, l_4 = 1}^m\,
\left[
-
\delta_{l_1}^j\,6\,\mathcal{X}_{y^{l_2}y^{l_3}y^{l_4}}
\right]
y_1^{l_1}y_1^{l_2}y_1^{l_3}y_2^{l_4}
&
=
\sum_{ l_1, l_2, l_3, l_4 = 1}^m\,
\left[
-
\delta_{l_2}^j\,6\,\mathcal{X}_{y^{l_1}y^{l_3}y^{l_4}}
\right]
y_1^{l_1}y_1^{l_2}y_1^{l_3}y_2^{l_4}
= 
\\
&
=
\sum_{ l_1, l_2, l_3, l_4 = 1}^m\,
\left[
-
\delta_{l_3}^j\,6\,\mathcal{X}_{y^{l_1}y^{l_2}y^{l_4}}
\right]
y_1^{l_1}y_1^{l_2}y_1^{l_3}y_2^{l_4}.
\endaligned
\end{equation}
In fact, we have systematically used such symmetries during the
intermediate computations (not exposed here) which we achieved
manually to obtain the final expressions of ${\bf Y}_1^j$, of ${\bf
Y}_2^j$, of ${\bf Y}_3^j$ and of ${\bf Y}_4^j$. To fix ideas, we have
always choosen the first index. Here, the first index is $l_1$; in the
first sum of line~9 of~\thetag{ 4.8}, the first index $l_\alpha$ for
the second weight $12$ is $l_2$.

This rule may be confirmed by inspecting all the monomials of ${\bf
Y}_2^j$, of ${\bf Y}_3^j$, of ${\bf Y}_4^j$ (and also of ${\bf
Y}_5^j$, which we have computed in a 
manuscript, but not copied in this Latex file).

From these considerations, we deduce that for the general formula, the
weight decomposition is simply $\mu_1\lambda_1, \dots, \mu_d\lambda_d$
and that the Kronecker symbol $\delta_\alpha^j$ is intrinsically
associated to the weights. In conclusion,
building on inductive reasonings, we have formulated the following
statement.

\def\thetheorem{4.18}\begin{theorem}
For one independent variable $x$, for several dependent variables
$(y^1, \dots, y^m)$ and for $\kappa \geq 1$, the general expression of
the coefficient ${\bf Y }_\kappa^j$ of the prolongation~\thetag{ 4.3} 
of a vector field is{\rm :}
\def\theequation{4.19}\begin{equation}
\boxed{
\small
\aligned
{\bf Y}_\kappa^j
&
=
\mathcal{Y}_{x^\kappa}^j
+
\sum_{d=1}^{\kappa+1}
\
\sum_{1\leq\lambda_1<\cdots<\lambda_d\leq\kappa}
\
\sum_{\mu_1\geq 1,\dots,\mu_d\geq 1}
\
\sum_{\mu_1\lambda_1+\cdots+\mu_d\lambda_d\leq\kappa+1} 
\\
& \
\ \ \ \ \
\sum_{l_{1:1}=1}^m
\cdots
\sum_{l_{1:\mu_1}=1}^m
\cdots\cdots
\sum_{l_{d:1}=1}^m
\cdots
\sum_{l_{d:\mu_d}=1}^m
\
\frac{\kappa(\kappa-1)\cdots
(\kappa-\mu_1\lambda_1+\cdots+\mu_d\lambda_d+2)}{
(\lambda_1!)^{\mu_1}\ \mu_1!\cdots(\lambda_d!)^{\mu_d}\ \mu_d!
}
\\
& \
\left[
\aligned
&
(\kappa-\mu_1\lambda_1-\cdots-\mu_d\lambda_d+1)
\frac{\partial^{\kappa-\mu_1\lambda_1-\cdots-\mu_d\lambda_d+
\mu_1+\cdots+\mu_d}\mathcal{Y}^j}{
(\partial x)^{\kappa-\mu_1\lambda_1-\cdots-\mu_d\lambda_d}
\partial y^{l_{1:1}}
\cdots
\partial y^{l_{1:\mu_1}}
\cdots
\partial y^{l_{d:1}}
\cdots
\partial y^{l_{d:\mu_d}}}
- \\
& \
-
\delta_{l_{1:1}}^j\,\mu_1\lambda_1\,
\frac{\partial^{\kappa-\mu_1\lambda_1-\cdots-\mu_d\lambda_d+
\mu_1+\cdots+\mu_d}\mathcal{X}}{
(\partial x)^{\kappa-\mu_1\lambda_1-\cdots-\mu_d\lambda_d+1}
\widehat{\partial y^{l_{1:1}}}
\cdots
\partial y^{l_{1:\mu_1}}
\cdots
\partial y^{l_{d:1}}
\cdots
\partial y^{l_{d:\mu_d}}}
-
\\
& 
\ \ \ \ \ \ \ \ \
-
\cdots
-
\\
& \
-
\delta_{l_{d:1}}^j\,\mu_d\lambda_d\,
\frac{\partial^{\kappa-\mu_1\lambda_1-\cdots-\mu_d\lambda_d+
\mu_1+\cdots+\mu_d}\mathcal{X}}{
(\partial x)^{\kappa-\mu_1\lambda_1-\cdots-\mu_d\lambda_d+1}
\partial y^{l_{1:1}}
\cdots
\partial y^{l_{1:\mu_1}}
\cdots
\widehat{\partial y^{l_{d:1}}}
\cdots
\partial y^{l_{d:\mu_d}}}
\endaligned
\right]
\cdot
\\
& \
\ \ \ \ \ \ \ \ \ \ \ \ \ \ \ \ \ \ \  
\ \ \ \ \ \ \ \ \ \ \ \ \ \ \ \ \ \ \ 
\cdot
\prod_{1\leq\nu_1\leq\mu_1}
y_{\lambda_1}^{l_{1:\nu_1}}
\ \cdots\cdots
\prod_{1\leq\nu_d\leq\mu_d}
y_{\lambda_d}^{l_{d:\nu_d}}.
\endaligned
}
\end{equation}
Here, the notation $\widehat{ \partial y^l}$ means 
that the partial derivative is dropped.
\end{theorem}

\subsection*{ 4.20.~Deduction of a multivariate Fa\`a
di Bruno formula} Let $m \in \N$ with $m\geq 1$, let $y = (y^1,\dots,
y^m) \in \K^m$, let $f = f( y^1, \dots, y^m)$ be a $\mathcal{
C}^\infty$-smooth function from $\K^m$ to $\K$, let $x \in \K$ and let
$g^1 = g^1(x), \dots, g^m = g^m( x)$ be $\mathcal{ C}^\infty$
functions from $\K$ to $\K$. The goal is to obtain an explicit formula
for the derivatives, with respect to $x$, of the composition $h :=
f\circ g$, namely $h(x) := f \left( g^1(x), \dots, g^m(x)
\right)$. For $\lambda \in \N$ with $\lambda \geq 1$, and for $j= 1,
\dots, m$, we shall abbreviate the derivative $\frac{ d^\lambda g^j}{
dx^\lambda}$ by $g_\lambda^j$ and similarly for $h_\lambda$. The
partial derivatives $\frac{ \partial^\lambda f}{ \partial
y^{l_1}\cdots \partial y^{l_\lambda}}$ will be abbreviated by $f_{l_1,
\dots, l_\lambda }$.

Appying the chain rule, we may compute:
\def\theequation{4.21}\begin{equation}
\small
\aligned
h_1
& 
=
\sum_{l_1=1}^m\,f_{l_1}\,g_1^{l_1}, 
\\
h_2
& 
=
\sum_{l_1,l_2=1}^m\,f_{l_1,l_2}\,g_1^{l_1}\,g_1^{l_2}
+
\sum_{l_1=1}^m\,f_{l_1}\,g_2^{l_1}, 
\\
h_3
& 
=
\sum_{l_1,l_2,l_3=1}^m\,f_{l_1,l_2,l_3}\,
g_1^{l_1}\,g_1^{l_2}\,g_1^{l_3}
+
3\,\sum_{l_1,l_2=1}^m\,f_{l_1,l_2}\,
g_1^{l_1}\,g_2^{l_2}
+
\sum_{l_1=1}^m\,f_{l_1}\,g_3^{l_1},
\\
h_4
& 
=
\sum_{l_1,l_2,l_3,l_4=1}^m\,f_{l_1,l_2,l_3,l_4}\,
g_1^{l_1}\,g_1^{l_2}\,g_1^{l_3}\,g_1^{l_4}
+
6\,\sum_{l_1,l_2,l_3=1}^m\,f_{l_1,l_2,l_3}\,
g_1^{l_1}\,g_1^{l_2}\,g_2^{l_3}
+ \\
& \
\ \ \ \ \
+
3\,\sum_{l_1,l_2=1}^m\,f_{l_1,l_2}\,
g_2^{l_1}\,g_2^{l_2}
+
4\,\sum_{l_1,l_2=1}^m\,f_{l_1,l_2}\,
g_1^{l_1}\,g_3^{l_2}
+
\sum_{l_1=1}^m\,
f_{l_1}\,g_4^{l_1}, 
\\
h_5
& 
=
\sum_{l_1,l_2,l_3,l_4,l_5=1}^m\,f_{l_1,l_2,l_3,l_4,l_5}\,
g_1^{l_1}\,g_1^{l_2}\,g_1^{l_3}\,g_1^{l_4}\,g_1^{l_5}
+
10\,\sum_{l_1,l_2,l_3,l_4=1}^m\,f_{l_1,l_2,l_3,l_4}\,
g_1^{l_1}\,g_1^{l_2}\,g_1^{l_3}\,g_2^{l_4}
+ \\
& \
\ \ \ \ \
+
15\,\sum_{l_1,l_2,l_3=1}^m\,f_{l_1,l_2,l_3}\,
g_1^{l_1}\,g_2^{l_2}\,g_2^{l_3}
+
10\,\sum_{l_1,l_2,l_3=1}^m\,f_{l_1,l_2,l_3}\,
g_1^{l_1}\,g_1^{l_2}\,g_3^{l_3}
+ \\
& \
\ \ \ \ \
+
10\,\sum_{l_1,l_2=1}^m\,f_{l_1,l_2}\,
g_2^{l_1}\,g_3^{l_2}
+
5\,\sum_{l_1,l_2=1}^m\,f_{l_1,l_2}\,
g_1^{l_1}\,g_4^{l_2}
+
\sum_{l_1=1}^m\,
f_{l_1}\,g_5^{l_1}.
\endaligned
\end{equation}
Introducing the derivations
\def\theequation{4.22}\begin{equation}
\small
\aligned
F^2
&
:= 
\sum_{l_1=1}^m\,g_2^{l_1}\,
\frac{\partial}{\partial g_1^{l_1}}
+
\sum_{l_1=1}^m\,g_1^{l_1}
\left(
\sum_{l_2=1}^m\,
f_{l_1,l_2}\,\frac{\partial}{\partial f_{l_2}}
\right), \\
F^3
&
:=
\sum_{l_1=1}^m\,g_2^{l_1}\,
\frac{\partial}{\partial g_1^{l_1}}
+
\sum_{l_1=1}^m\,g_3^{l_1}\,
\frac{\partial}{\partial g_2^{l_1}}
+
\sum_{l_1=1}^m\,g_1^{l_1}
\left(
\sum_{l_2=1}^m\,
f_{l_1,l_2}\,\frac{\partial}{\partial f_{l_2}}
+
\sum_{l_2,l_3=1}^m\,
f_{l_1,l_2,l_3}\,\frac{\partial}{\partial f_{l_2,l_3}}
\right), 
\\
& \
\ \ \ \ \
\text{\bf 
\dots\dots\dots\dots\dots\dots\dots\dots\dots
\dots\dots\dots\dots\dots\dots\dots\dots\dots
\dots\dots\dots\dots\dots\dots\dots\dots\dots
} \\
F^\lambda
&
:=
\sum_{l_1=1}^m\,g_2^{l_1}\,
\frac{\partial}{\partial g_1^{l_1}}
+
\sum_{l_1=1}^m\,g_3^{l_1}\,
\frac{\partial}{\partial g_2^{l_1}}
+
\cdots
+
\sum_{l_1=1}^m\,g_\lambda^{l_1}\,
\frac{\partial}{\partial g_{\lambda-1}^{l_1}}
+ \\
& \
\ \ \ \ \
+
\sum_{l_1=1}^m\,g_1^{l_1}
\left(
\sum_{l_2=1}^m\,
f_{l_1,l_2}\,\frac{\partial}{\partial f_{l_2}}
+
\sum_{l_2,l_3=1}^m\,
f_{l_1,l_2,l_3}\,\frac{\partial}{\partial f_{l_2,l_3}}
+
\cdots
+
\sum_{l_2,\dots,l_\lambda=1}^m\,
f_{l_1,l_2,\dots,l_\lambda}\,
\frac{\partial}{\partial f_{l_2,\dots,l_\lambda}}
\right),
\endaligned
\end{equation}
we observe that the following
induction relations hold:
\def\theequation{4.23}\begin{equation}
\aligned
h_2
&
=
F^2
\left(
h_1
\right), 
\\
h_3
&
=
F^3
\left(
h_2
\right), 
\\
\text{\bf 
\dots\dots
}
& \
\ \ \ \ \
\text{\bf 
\dots\dots\dots
}
\\
h_\lambda
&
=
F^\lambda
\left(
h_{\lambda-1}
\right).
\endaligned
\end{equation}
To obtain the explicit version of the Fa\`a di Bruno in the case of
one variable $x$ and several variables $(y^1, \dots, y^m)$, it
suffices to extract from the expression of ${\bf Y}_\kappa^j$ provided
by Theorem~4.18 only the terms corresponding to $\mu_1 \lambda_1 +
\cdots + \mu_d\lambda_d = \kappa$, dropping all the $\mathcal{ X}$
terms. After some simplifications and after a translation by means of
an elementary dictionary, we may formulate a statement.

\def\thetheorem{4.24}\begin{theorem} 
For every integer $\kappa \geq 1$, the $\kappa$-th partial derivative
of the composite function $h = h( x) = f \left( 
g^1(x), \dots, g^m(x) \right)$ with
respect to $x$ may be expressed as an explicit polynomial depending on
the partial derivatives of $f$, on the derivatives of $g$ and having
integer coefficients{\rm:}
\def\theequation{4.25}\begin{equation}
\boxed{
\aligned
\frac{d^\kappa h}{dx^\kappa}
&
=
\sum_{d=1}^\kappa
\
\sum_{1\leq \lambda_1 < \cdots < \lambda_d \leq \kappa}
\
\sum_{\mu_1\geq 1,\dots,\mu_d\geq 1}
\
\sum_{\mu_1\lambda_1+\cdots+\mu_d\lambda_d=\kappa}
\
\frac{\kappa!}{
(\lambda_1!)^{\mu_1}\ \mu_1! 
\cdots
(\lambda_d!)^{\mu_d}\ \mu_d!
}
\\
& \
\ \ \ \ \ \ \ \ \ \ \ \ \ \ \
\sum_{l_{1:1},\dots,l_{1:\mu_1}=1}^m
\ \cdots \
\sum_{l_{d:1},\dots,l_{d:\mu_d}=1}^m
\\
\\
& \
\ \ \ \ \
\frac{\partial^{\mu_1+\cdots+\mu_d}f}{
\partial y^{l_{1:1}} 
\cdots
\partial y^{l_{1:\mu_1}}
\cdots
\partial y^{l_{d:1}} 
\cdots
\partial y^{l_{d:\mu_d}}
}
\
\prod_{1\leq\nu_1\leq\mu_1}
\frac{d^{\lambda_1} g^{l_{1:\nu_1}}}{d x^{\lambda_1}}
\ \cdots
\prod_{1\leq\nu_d\leq\mu_d}
\frac{d^{\lambda_d} g^{l_{d:\nu_d}}}{d x^{\lambda_d}}.
\endaligned
}
\end{equation}
\end{theorem} 

\section*{\S5.~Several independent variables and
several dependent variables}

\subsection*{5.1.~Expression of ${\bf Y}_{i_1}^j$, 
of ${\bf Y}_{i_1,i_2}^j$ and of ${\bf Y}_{i_1,i_2,i_3}^j$} Applying
the induction~\thetag{1.31} and working out the obtained formulas
until they take a perfect shape, we obtain firstly:
\def\theequation{5.2}\begin{equation}
\small
{\bf Y}_{i_1}^j
=
\mathcal{Y}_{x^{i_1}}^j
+
\sum_{l_1=1}^m\ \sum_{k_1=1}^n
\left[
\delta_{i_1}^{k_1}\,\mathcal{Y}_{y^{l_1}}^j
-
\delta_{l_1}^j\,
\mathcal{X}_{x^{i_1}}^{k_1}
\right]
y_{k_1}^{l_1}
+ 
\sum_{l_1,l_2=1}^m\ \sum_{k_1,k_2=1}^n
\left[
-
\delta_{l_2}^j\,
\delta_{i_1}^{k_1}\,\mathcal{X}_{y^{l_1}}^{k_2}
\right]
y_{k_1}^{l_1}y_{k_2}^{l_2}.
\end{equation}
Secondly:
\def\theequation{5.3}\begin{equation}
\small
\aligned{\bf Y}_{i_1,i_2}^j
&
=
\mathcal{Y}_{x^{i_1}x^{i_2}}^j
+
\sum_{l_1=1}^m\ \sum_{k_1=1}^n
\left[
\delta_{i_1}^{k_1}\,\mathcal{Y}_{x^{i_2}y^{l_1}}^j
+
\delta_{i_2}^{k_1}\,\mathcal{Y}_{x^{i_1}y^{l_1}}^j
-
\delta_{l_1}^j\,
\mathcal{X}_{x^{i_1}x^{i_2}}^{k_1}
\right]
y_{k_1}^{l_1}
+ \\
& \
\ \ \ \ \
+
\sum_{l_1,l_2=1}^m\ \sum_{k_1,k_2=1}^n
\left[
\delta_{i_1, \ i_2}^{k_1,k_2}\,
\mathcal{Y}_{y^{l_1}y^{l_2}}^j
-
\delta_{l_2}^j\,\delta_{i_1}^{k_1}\,
\mathcal{X}_{x^{i_2}y^{l_1}}^{k_2}
-
\delta_{l_2}^j\,\delta_{i_2}^{k_1}\,
\mathcal{X}_{x^{i_1}y^{l_1}}^{k_2}
\right]
y_{k_1}^{l_1}y_{k_2}^{l_2}
+ \\
& \
\ \ \ \ \
+
\sum_{l_1,l_2,l_3=1}^m\ \sum_{k_1,k_2,k_3=1}^n
\left[
-
\delta_{l_3}^j\,\delta_{i_1,\ i_2}^{k_1,k_2}\,
\mathcal{X}_{y^{l_1}y^{l_2}}^{k_3}
\right]
y_{k_1}^{l_1}y_{k_2}^{l_2}y_{k_3}^{l_3}
+ \\
& \
\ \ \ \ \
+
\sum_{l_1=1}^m\ \sum_{k_1,k_2=1}^n
\left[
\delta_{i_1,\ i_2}^{k_1,k_2}\,
\mathcal{Y}_{y^{l_1}}^j
-
\delta_{l_1}^j\,\delta_{i_1}^{k_1}\,
\mathcal{X}_{x^{i_2}}^{k_2}
-
\delta_{l_1}^j\,\delta_{i_2}^{k_1}\,
\mathcal{X}_{x^{i_1}}^{k_2}
\right]
y_{k_1,k_2}^{l_1}
+ \\
& \
\ \ \ \ \
+
\sum_{l_1,l_2=1}^m\ \sum_{k_1,k_2,k_3=1}^n
\left[
-
\delta_{l_1}^j\,\delta_{i_1,\ i_2}^{k_2,k_3}\,
\mathcal{X}_{y^{l_2}}^{k_1}
-
\delta_{l_2}^j\,\delta_{i_1,\ i_2}^{k_3,k_1}\,
\mathcal{X}_{y^{l_1}}^{k_2}
-
\delta_{l_2}^j\,\delta_{i_1,\ i_2}^{k_1,k_2}\,
\mathcal{X}_{y^{l_1}}^{k_3}
\right]
y_{k_1}^{l_1}y_{k_2}^{l_2}y_{k_3}^{l_3}.
\endaligned
\end{equation}
Thirdly:
$$
\small
\aligned
{\bf Y}_{i_1,i_2,i_3}^j
&
=
\mathcal{Y}_{x^{i_1}x^{i_2}x^{i_3}}^j
+
\sum_{l_1=1}^m\ \sum_{k_1=1}^n
\left[
\delta_{i_1}^{k_1}\,\mathcal{Y}_{x^{i_2}x^{i_3}y^{l_1}}^j
+
\delta_{i_2}^{k_1}\,\mathcal{Y}_{x^{i_1}x^{i_3}y^{l_1}}^j
+
\delta_{i_3}^{k_1}\,\mathcal{Y}_{x^{i_1}x^{i_2}y^{l_1}}^j
-
\delta_{l_1}^j\,
\mathcal{X}_{x^{i_1}x^{i_2}x^{i_3}}^{k_1}
\right]
y_{k_1}^{l_1}
+ \\
& \
\ \ \ \ \
+
\sum_{l_1,l_2=1}^m\ \sum_{k_1,k_2=1}^n
\left[
\delta_{i_1, \ i_2}^{k_1,k_2}\,
\mathcal{Y}_{x^{i_3}y^{l_1}y^{l_2}}^j
+
\delta_{i_3, \ i_1}^{k_1,k_2}\,
\mathcal{Y}_{x^{i_2}y^{l_1}y^{l_2}}^j
+
\delta_{i_2, \ i_3}^{k_1,k_2}\,
\mathcal{Y}_{x^{i_1}y^{l_1}y^{l_2}}^j
- 
\right. 
\\
& \
\ \ \ \ \ \ \ \ \ \ \ \ \ \ \ \ \ \ \ \ 
\ \ \ \ \ \ \ \ \ \ \ \ \ \ \ \ \
\left.
-
\delta_{l_2}^j\,\delta_{i_1}^{k_1}\,
\mathcal{X}_{x^{i_2}x^{i_3}y^{l_1}}^{k_2}
-
\delta_{l_2}^j\,\delta_{i_2}^{k_1}\,
\mathcal{X}_{x^{i_1}x^{i_3}y^{l_1}}^{k_2}
-
\delta_{l_2}^j\,\delta_{i_3}^{k_1}\,
\mathcal{X}_{x^{i_1}x^{i_2}y^{l_1}}^{k_2}
\right]
y_{k_1}^{l_1}y_{k_2}^{l_2}
+ \\
& \
\ \ \ \ \
+
\sum_{l_1,l_2,l_3=1}^m\ \sum_{k_1,k_2,k_3=1}^n
\left[
\delta_{i_1, \ i_2, \ i_3}^{k_1,k_2,k_3}\,
\mathcal{Y}_{y^{l_1}y^{l_2}y^{l_3}}^j
-
\delta_{l_3}^j\,\delta_{i_1,\ i_2}^{k_1,k_2}\,
\mathcal{X}_{x^{i_3}y^{l_1}y^{l_2}}^{k_3}
-
\right.
\\
& \
\ \ \ \ \ \ \ \ \ \ \ \ \ \ \ \ \ \ \ \ \
\ \ \ \ \ \ \ \ \ \ \ \ \ \ \ \ \ \ \ \ \
\ \ \ \ 
\left.
-
\delta_{l_3}^j\,\delta_{i_1,\ i_3}^{k_1,k_2}\,
\mathcal{X}_{x^{i_2}y^{l_1}y^{l_2}}^{k_3}
-
\delta_{l_3}^j\,\delta_{i_2,\ i_3}^{k_1,k_2}\,
\mathcal{X}_{x^{i_1}y^{l_1}y^{l_2}}^{k_3}
\right]
y_{k_1}^{l_1}y_{k_2}^{l_2}y_{k_3}^{l_3}
+
\endaligned
$$
$$
\small
\aligned
& \
\ \ \ \ \ 
+
\sum_{l_1,l_2,l_3,l_4=1}^m\ \sum_{k_1,k_2,k_3,k_4=1}^n
\left[
-
\delta_{l_4}^j\,\delta_{i_1,\ i_2,\ i_3}^{k_1,k_2,k_3}\,
\mathcal{X}_{y^{l_1}y^{l_2}y^{l_3}}^{k_4}
\right]
y_{k_1}^{l_1}y_{k_2}^{l_2}y_{k_3}^{l_3}y_{k_4}^{l_4}
+ \\
& \
\ \ \ \ \
+
\sum_{l_1=1}^m\ \sum_{k_1,k_2=1}^n
\left[
\delta_{i_1,\ i_2}^{k_1,k_2}\,
\mathcal{Y}_{x^{i_3}y^{l_1}}^j
+
\delta_{i_3,\ i_1}^{k_1,k_2}\,
\mathcal{Y}_{x^{i_2}y^{l_1}}^j
+
\delta_{i_2,\ i_3}^{k_1,k_2}\,
\mathcal{Y}_{x^{i_1}y^{l_1}}^j
- 
\right. \\
& \
\ \ \ \ \ \ \ \ \ \ \ \ \ \ \ \ \ \ \ \
\ \ \ \ \ \ \ \ \ \ \ \ \ \
\left.
-
\delta_{l_1}^j\,\delta_{i_1}^{k_1}\,
\mathcal{X}_{x^{i_2}x^{i_3}}^{k_2}
-
\delta_{l_1}^j\,\delta_{i_2}^{k_1}\,
\mathcal{X}_{x^{i_1}x^{i_3}}^{k_2}
-
\delta_{l_1}^j\,\delta_{i_3}^{k_1}\,
\mathcal{X}_{x^{i_1}x^{i_2}}^{k_2}
\right]
y_{k_1,k_2}^{l_1}
+ 
\endaligned
$$
\def\theequation{5.4}\begin{equation}
\small
\aligned
& \
\ \ \ \ \
+
\sum_{l_1,l_2=1}^m\ \sum_{k_1,k_2,k_3=1}^n
\left[
\delta_{i_1,\ i_2,\ i_3}^{k_1,k_2,k_3}\,
\mathcal{Y}_{y^{l_1}y^{l_2}}^j
+
\delta_{i_1,\ i_2,\ i_3}^{k_3,k_1,k_2}\,
\mathcal{Y}_{y^{l_1}y^{l_2}}^j
+
\delta_{i_1,\ i_2,\ i_3}^{k_2,k_3,k_1}\,
\mathcal{Y}_{y^{l_1}y^{l_2}}^j
-
\right. 
\\
& \
\ \ \ \ \ \ \ \ \ \ \ \ \ \ \ \ \ \ \ \
\ \ \ \ \ \ \ \ \ \ \ \ \ \ \ \ \ \ \ \
\ \ 
\left.
-
\delta_{l_1}^j\,\delta_{i_1,\ i_2}^{k_2,k_3}\,
\mathcal{X}_{x^{i_3}y^{l_2}}^{k_1}
-
\delta_{l_1}^j\,\delta_{i_1,\ i_3}^{k_2,k_3}\,
\mathcal{X}_{x^{i_2}y^{l_2}}^{k_1}
-
\delta_{l_1}^j\,\delta_{i_2,\ i_3}^{k_2,k_3}\,
\mathcal{X}_{x^{i_1}y^{l_2}}^{k_1}
-
\right. 
\\
& \
\ \ \ \ \ \ \ \ \ \ \ \ \ \ \ \ \ \ \ \
\ \ \ \ \ \ \ \ \ \ \ \ \ \ \ \ \ \ \ \
\ \ 
\left.
-
\delta_{l_2}^j\,\delta_{i_1,\ i_2}^{k_3,k_1}\,
\mathcal{X}_{x^{i_3}y^{l_1}}^{k_2}
-
\delta_{l_2}^j\,\delta_{i_1,\ i_3}^{k_3,k_1}\,
\mathcal{X}_{x^{i_2}y^{l_1}}^{k_2}
-
\delta_{l_2}^j\,\delta_{i_2,\ i_3}^{k_3,k_1}\,
\mathcal{X}_{x^{i_1}y^{l_1}}^{k_2}
-
\right. 
\\
& \
\ \ \ \ \ \ \ \ \ \ \ \ \ \ \ \ \ \ \ \
\ \ \ \ \ \ \ \ \ \ \ \ \ \ \ \ \ \ \ \
\ \ 
\left.
-
\delta_{l_2}^j\,\delta_{i_1,\ i_2}^{k_1,k_2}\,
\mathcal{X}_{x^{i_3}y^{l_1}}^{k_3}
-
\delta_{l_2}^j\,\delta_{i_1,\ i_3}^{k_1,k_2}\,
\mathcal{X}_{x^{i_2}y^{l_1}}^{k_3}
-
\delta_{l_2}^j\,\delta_{i_2,\ i_3}^{k_1,k_2}\,
\mathcal{X}_{x^{i_1}y^{l_1}}^{k_3}
\right]
y_{k_1}^{l_1}y_{k_2,k_3}^{l_2}
+ \\
\endaligned
\end{equation}
$$
\small
\aligned
& \
\ \ \ \ \
+
\sum_{l_1,l_2,l_3=1}^m\ \sum_{k_1,k_2,k_3,k_4=1}^n
\left[
-
\delta_{l_3}^j\,\delta_{i_1,\ i_2,\ i_3}^{k_1,k_2,k_3}\,
\mathcal{X}_{y^{l_1}y^{l_2}}^{k_4}
-
\delta_{l_3}^j\,\delta_{i_1,\ i_2,\ i_3}^{k_2,k_3,k_1}\,
\mathcal{X}_{y^{l_1}y^{l_2}}^{k_4}
-
\delta_{l_3}^j\,\delta_{i_1,\ i_2,\ i_3}^{k_3,k_2,k_1}\,
\mathcal{X}_{y^{l_1}y^{l_2}}^{k_4}
- 
\right. 
\\
& \
\ \ \ \ \ \ \ \ \ \ \ \ \ \ \ \ \ \ \ \
\ \ \ \ \ \ \ \ \ \ \ \ \ \ \ \ \ \ \ \ 
\ \ \ \ \ \ \ \ \ \
\left.
-
\delta_{l_2}^j\,\delta_{i_1,\ i_2,\ i_3}^{k_3,k_4,k_1}\,
\mathcal{X}_{y^{l_1}y^{l_3}}^{k_2}
-
\delta_{l_2}^j\,\delta_{i_1,\ i_2,\ i_3}^{k_3,k_1,k_4}\,
\mathcal{X}_{y^{l_1}y^{l_3}}^{k_2}
-
\right. 
\\
& \
\ \ \ \ \ \ \ \ \ \ \ \ \ \ \ \ \ \ \ \
\ \ \ \ \ \ \ \ \ \ \ \ \ \ \ \ \ \ \ \ 
\ \ \ \ \ \ \ \ \ \ \ \ \ \ \ \ \ \ \ \
\ \ \ \ \ \ \ \ \ \ \ \ \ \ \ \ \ \ \ \
\ \ \ \ \ \ \ \ \ \ \ \ \ \ \ \ \ \ \ \ 
\left.
-
\delta_{l_2}^j\,\delta_{i_1,\ i_2,\ i_3}^{k_1,k_3,k_4}\,
\mathcal{X}_{y^{l_1}y^{l_3}}^{k_2}
\right]
y_{k_1}^{l_1}y_{k_2}^{l_2}y_{k_3,k_4}^{l_3}
+ \\
& \
\ \ \ \ \
+
\sum_{l_1,l_2=1}^m\ \sum_{k_1,k_2,k_3,k_4=1}^n
\left[
-
\delta_{l_2}^j\,\delta_{i_1,\ i_2,\ i_3}^{k_1,k_2,k_3}\,
\mathcal{X}_{y^{l_1}}^{k_3}
-
\delta_{l_2}^j\,\delta_{i_1,\ i_2,\ i_3}^{k_2,k_4,k_1}\,
\mathcal{X}_{y^{l_1}}^{k_3}
-
\delta_{l_2}^j\,\delta_{i_1,\ i_2,\ i_3}^{k_4,k_1,k_2}\,
\mathcal{X}_{y^{l_1}}^{k_3}
\right]
y_{k_1,k_2}^{l_1}y_{k_3,k_4}^{l_2}
+ \\
\endaligned
$$
$$
\small
\aligned
& \
\ \ \ \ \
+
\sum_{l_1=1}^m\ \sum_{k_1,k_2,k_3=1}^n
\left[
\delta_{i_1,\ i_2, \i_3}^{k_1,k_2,k_3}\,
\mathcal{Y}_{y^{l_1}}^j
-
\delta_{l_1}^j\,\delta_{i_1,\ i_2}^{k_1,k_2}\,
\mathcal{X}_{x^{i_3}}^{k_3}
-
\delta_{l_1}^j\,\delta_{i_1,\ i_3}^{k_1,k_2}\,
\mathcal{X}_{x^{i_2}}^{k_3}
-
\delta_{l_1}^j\,\delta_{i_2,\ i_3}^{k_1,k_2}\,
\mathcal{X}_{x^{i_1}}^{k_3}
\right]
y_{k_1,k_2,k_3}^{l_1}
+ \\
& \
\ \ \ \ \
+
\sum_{l_1,l_2=1}^m\ \sum_{k_1,k_2,k_3,k_4=1}^n
\left[
-
\delta_{l_2}^j\,\delta_{i_1,\ i_2,\ i_3}^{k_1,k_2,k_3}\,
\mathcal{X}_{y^{l_1}}^{k_4}
-
\delta_{l_2}^j\,\delta_{i_1,\ i_2,\ i_3}^{k_4,k_1,k_2}\,
\mathcal{X}_{y^{l_1}}^{k_3}
-
\delta_{l_2}^j\,\delta_{i_1,\ i_2,\ i_3}^{k_3,k_4,k_1}\,
\mathcal{X}_{y^{l_1}}^{k_2}
-
\right.
\\
& \
\ \ \ \ \ \ \ \ \ \ \ \ \ \ \ \ \ \ \ \
\ \ \ \ \ \ \ \ \ \ \ \ \ \ \ \ \ \ \ \
\ \ \ \ \ \ \ \
\left.
-
\delta_{l_1}^j\,\delta_{i_1,\ i_2,\ i_3}^{k_2,k_3,k_4}\,
\mathcal{X}_{y^{l_2}}^{k_1}
\right]
y_{k_1}^{l_1}y_{k_2,k_3,k_4}^{l_2}.
\endaligned
$$

\subsection*{5.5.~Final synthesis}
To obtain the general formula for ${\bf Y}_{ i_1, \dots, i_\kappa}^j$,
we have to achieve the synthesis between the two formulas~\thetag{
3.74} and~\thetag{ 4.19}. We start with~\thetag{ 3.74} and we modify
it until we reach the final formula for ${\bf Y}_{ i_1, \dots,
i_\kappa }^j$.

We have to add the $\mu_1+\cdots+\mu_d$ 
sums $\sum_{ l_{ 1:1} =1 }^m \cdots \sum_{
l_{ 1: \mu_1 } =1 }^m \cdots \cdots \sum_{ l_{ d:1} =1}^m \cdots
\sum_{l_{ d: \mu_d }= 1}^m$, together with various indices
$l_\alpha$. About these indices, the only point which is not obvious
may be analyzed as follows.

A permutation $\sigma \in \mathfrak{ F}_{
\mu_1\lambda_1 + \cdots + \mu_d \lambda_d }^{ (\mu_1, \lambda_1),
\dots, (\mu_d, \lambda_d)}$ yields the list:
\def\theequation{5.6}\begin{equation}
\aligned
& \
\sigma(1\!:\!1\!:\!1),\dots,\sigma(1\!:\!1\!:\!\lambda_1),
\dots
\sigma(1\!:\!\mu_1\!:\!1),\dots,\sigma(1\!:\!\mu_1\!:\!\lambda_1),
\dots
\\
& \
\ \ \ \ \ \ \ \ \ 
\dots,
\sigma(d\!:\!1\!:\!1),\dots,\sigma(1\!:\!1\!:\!\lambda_d),
\dots
\sigma(d\!:\!\mu_d\!:\!1),\dots,\sigma(d\!:\!\mu_d\!:\!\lambda_d),
\endaligned
\end{equation} 
In the end of the 
sixth line of~\thetag{ 3.74}, the last term $\sigma ( d\! : \!
\mu_d \! : \! \lambda_d )$ of the above list appears as the subscript
of the upper index $k_{ \sigma (d: \mu_d: \lambda_d )}$ of the term
$\mathcal{ X }^{ k_{ \sigma( d:\mu_d: \lambda_d) }}$. According to the
formal rules of Theorem~4.19, we have to multiply the partial
derivative of $\mathcal{ X }^{ k_{ \sigma( d:\mu_d: \lambda_d ) }}$ by
a certain Kronecker symbol $\delta_{ l_\alpha}^j$. The question is:
{\it what is the subscript $\alpha$ and how
to denote it}?

As explained before the statement of Theorem~4.19, the subscript
$\alpha$ is obtained as follows. The term $\sigma ( d\! : \! \mu_d \!
: \! \lambda_d )$ is of the form $(e\! : \! \nu_e \! : \! \gamma_e
)$, for some $e$ with $1\leq e \leq d$, for some $\nu_e$ with $1\leq
\nu_e \leq \mu_e$ and for some $\gamma_e$ with $1\leq \gamma_e \leq
\lambda_e$. The single pure jet variable
\def\theequation{5.7}\begin{equation}
\aligned
y_{k_{e:\nu_e:1},\dots,
k_{e:\nu_e:\gamma_e},\dots,
k_{e:\nu_e:\lambda_e}}^{l_{e:\nu_e}}
\endaligned
\end{equation}
appears inside the total monomial
\def\theequation{5.8}\begin{equation}
\aligned
\prod_{1\leq\nu_1\leq\mu_1}\,
y_{k_{1:\nu_1:1},\dots,k_{1:\nu_1:\lambda_1}}^{l_{1:\nu_1}}
\ \cdots \
\prod_{1\leq\nu_d\leq\mu_d}\,
y_{k_{d:\nu_d:1},\dots,k_{d:\nu_d:\lambda_d}}^{l_{d:\nu_d}},
\endaligned
\end{equation}
placed at the end of the formula for ${\bf Y }_{ i_1, \dots,
i_\kappa}^j$ ({\it see} in advance formula~\thetag{ 5.13} below; this
total monomial generalizes to several dependent variables the total
monomial appearing in the last line of~\thetag{ 3.74}). According to
the rule explained before the statement of Theorem~4.18, the index
$l_\alpha$ must be equal to $l_{e : \nu_e }$, since
$l_{e : \nu_e }$ is attached to the monomial~\thetag{ 5.7}.
Coming back to the term
$\sigma ( d\! : \! \mu_d \! : \! \lambda_d )$, we shall denote this
index by
\def\theequation{5.9}\begin{equation}
\aligned
l_{e:\nu_e}
=:
l_{\pi(e:\nu_e:\gamma_e)}
=:
l_{\pi\sigma(d:\mu_d:\lambda_d)},
\endaligned
\end{equation}
where the symbol $\pi$ denotes the projection
from the set 
\def\theequation{5.10}\begin{equation}
\aligned
\{
1\!:\!1\!:\!1,\dots,1\!:\!\mu_1\!:\!\lambda_1,
\dots\dots,
d\!:\!1\!:\!1,\dots,d\!:\!\mu_d\!:\!\lambda_d
\}
\endaligned
\end{equation}
to the set
\def\theequation{5.11}\begin{equation}
\aligned
\{
1\!:\!1,\dots,1\!:\!\mu_1,
\dots,
d\!:\!1,\dots,d\!:\!\mu_d
\}
\endaligned
\end{equation}
simply defined by $\pi(e \! : \! \nu_e \! : \! \gamma_e) := (e \! : \!
\nu_e)$.

In conclusion, by means of this formalism, we may write down the
complete generalization of Theorems~2.24, 3.73 and~4.18 to several
independent variables and several dependent variables

\def\thetheorem{5.12}\begin{theorem}
For $j = 1, \dots, m$, for every $\kappa \geq 1$ and for every choice
of $\kappa$ indices $i_1,\dots, i_\kappa$ in the set $\{ 1, 2, \dots,
n\}$, the general expression of ${\bf Y}_{i_1, \dots, i_\kappa }^j$ is
as follows{\rm :}
\def\theequation{5.13}\begin{equation}
\small
\boxed{
\aligned
{\bf Y}_{i_1, \dots, i_\kappa}^j
&
=
\mathcal{Y}_{x^{i_1}\cdots x^{i_\kappa}}^j
+
\sum_{d=1}^{\kappa+1}
\ \
\sum_{1\leq\lambda_1<\cdots<\lambda_d\leq\kappa}
\ \
\sum_{\mu_1\geq 1,\dots,\mu_d\geq 1} 
\
\sum_{
\mu_1\lambda_1
+
\cdots
+
\mu_d\lambda_d\leq \kappa+1} 
\\
& \
\ \ \ \ \ \ \ \ \ \ \ \ \ \ \ \ \ \ \ \ 
\ \ \ \ \ \ \ \ \ \ \ \ \ 
\sum_{l_{1:1}=1}^m
\cdots
\sum_{l_{1:\mu_1}=1}^m
\cdots\cdots
\sum_{l_{d:1}=1}^m
\cdots
\sum_{l_{d:\mu_d}=1}^m
\\
&
\sum_{k_{1:1:1},\dots,k_{1:1:\lambda_1}=1}^n
\cdots \
\sum_{k_{1:\mu_1:1},\dots,k_{1:\mu_1:\lambda_1}=1}^n
\cdots\cdots \
\sum_{k_{d:1:1},\dots,k_{d:1:\lambda_d}=1}^n
\cdots \
\sum_{k_{d:\mu_d:1},\dots,k_{d:\mu_d:\lambda_d}=1}^n
\\
&
\left[
\aligned
& 
\sum_{\sigma\in\mathfrak{F}_{\mu_1\lambda_1+\cdots+\mu_d\lambda_d}^{
(\mu_1,\lambda_1),\dots,(\mu_d,\lambda_d)}}
\
\sum_{\tau\in\mathfrak{S}_\kappa^{
\mu_1\lambda_1+\cdots+\mu_d\lambda_d}}\,
\delta_{i_{\tau(1)},\dots,i_{\tau(\mu_1\lambda_1)},\dots,
i_{\tau(\mu_1\lambda_1+\cdots+\mu_d\lambda_d)}}^{
k_{\sigma(1:1:1)},\dots,k_{\sigma(1:\mu_1:\lambda_1)},
\dots,k_{\sigma(d:\mu_d:\lambda_d)}}
\cdot
\\
& \
\ \ \ \ \
\cdot
\frac{\partial^{\kappa-\mu_1\lambda_1-\cdots-\mu_d\lambda_d+
\mu_1+\cdots+\mu_d}
\mathcal{Y}^j}{
\partial x^{i_{\tau(\mu_1\lambda_1+\cdots+\mu_d\lambda_d+1)}}\cdots
\partial x^{i_{\tau(\kappa)}}
\partial y^{l_{1:1}}\cdots\partial y^{l_{d:\mu_d}}}\
- \\
& 
-
\sum_{\sigma\in\mathfrak{F}_{\mu_1\lambda_1+\cdots+\mu_d\lambda_d}^{
(\mu_1,\lambda_1),\dots,(\mu_d,\lambda_d)}}
\
\sum_{\tau\in\mathfrak{S}_\kappa^{
\mu_1\lambda_1+\cdots+\mu_d\lambda_d-1}}\,
\delta_{i_{\tau(1)},\dots,i_{\tau(\mu_1\lambda_1)},\dots,
i_{\tau(\mu_1\lambda_1+\cdots+\mu_d\lambda_d-1)}}^{
k_{\sigma(1:1:1)},\dots,k_{\sigma(1:\mu_1:\lambda_1)},
\dots,k_{\sigma(d:\mu_d:\lambda_d-1)}}
\cdot
\\
& \
\ \ \ \ \
\cdot
\delta_{l_{\pi\sigma(d:\mu_d:\lambda_d)}}^j
\cdot
\frac{\partial^{\kappa-\mu_1\lambda_1-\cdots-\mu_d\lambda_d
+\mu_1+\cdots+\mu_d}\mathcal{X}^{k_{\sigma(d:\mu_d:\lambda_d)}}}{
\partial x^{i_{\tau(\mu_1\lambda_1+\cdots+\mu_d\lambda_d)}}\cdots
\partial x^{i_{\tau(\kappa)}}
\partial y^{l_{1:1}}\cdots
\widehat{\partial y^{l_{\pi\sigma(d:\mu_d:\lambda_d)}}}
\cdots\partial y^{l_{d:\mu_d}}}
\endaligned
\right]
\cdot
\\
& \
\ \ \ \ \ \ \ \ \ \ \ \ \ \ \ \ \ \ \
\cdot
\prod_{1\leq\nu_1\leq\mu_1}\,
y_{k_{1:\nu_1:1},\dots,k_{1:\nu_1:\lambda_1}}^{l_{1:\nu_1}}
\ \cdots \
\prod_{1\leq\nu_d\leq\mu_d}\,
y_{k_{d:\nu_d:1},\dots,k_{d:\nu_d:\lambda_d}}^{l_{d:\nu_d}}.
\endaligned
}
\end{equation}
\end{theorem}

In this formula, the coset $\mathfrak{ F }_{ \mu_1 \lambda_1 + \cdots
+ \mu_d \lambda_d }^{ ( \mu_1, \lambda_1 ),\dots, ( \mu_d,
\lambda_d)}$ was defined in equation~\thetag{ 3.71};
as in Theorem~3.73, we have made the
identification{\rm :}
\def\theequation{5.14}\begin{equation}
\{1,\dots,\kappa\}
\equiv
\{
1\!:\!1\!:\!1,
\dots,
1\!:\!\mu_1\!:\!\lambda_1,
\text{\rm \dots\dots}, 
d\!:\!1\!:\!1,
\dots,
d\!:\!\mu_d\!:\!\lambda_d
\}.
\end{equation}

\subsection*{ 5.15.~Deduction of the 
most general multivariate Fa\`a di Bruno formula} Let $n \in \N$ with
$n\geq 1$, let $x = (x^1,\dots, x^n) \in \K^n$, let $m\in \N$ with
$m\geq 1$, let $g^j = g^j ( x^1, \dots, x^n)$, $j=1, \dots, m$, be
$\mathcal{ C}^\infty$-smooth functions from $\K^n$ to $\K^m$, let $y =
(y^1, \dots, y^m) \in \K^m$ and let $f = f(y^1, \dots, y^m)$ be a
$\mathcal{ C}^\infty$ function from $\K^m$ to $\K$. The goal is to
obtain an explicit formula for the partial derivatives of the
composition $h := f\circ g$, namely
\def\theequation{5.16}\begin{equation}
h(x^1, \dots, x^n) := f
\left(
g^1(x^1,\dots, x^n), \dots, g^m(x^1,\dots, x^n) 
\right).
\end{equation}
For $j= 1,\dots, m$, for $\lambda \in \N$ with $\lambda \geq 1$ and
for arbitrary indices $i_1, \dots, i_\lambda = 1, \dots, n$, we shall
abbreviate the partial derivative $\frac{ \partial^\lambda g^j}{
\partial x^{i_1} \cdots \partial x^{i_\lambda}}$ by $g_{i_1,\dots,
i_\lambda }^j$ and similarly for $h_{i_1, \dots, i_\lambda}$. For
arbitrary indices $l_1, \dots, l_\lambda = 1, \dots, m$, the partial
derivative $\frac{ \partial^\lambda f}{ \partial y^{l_1} \cdots
\partial y^{l_\lambda }}$ will be abbreviated by $f_{ l_1, \dots,
l_\lambda }$.

Appying the chain rule, we may compute:
\def\theequation{5.17}\begin{equation}
\small
\aligned
h_{i_1}
&
=
\sum_{l_1=1}^m\,f_{l_1}
\left[
g_{i_1}^{l_1}
\right],
\\
h_{i_1,i_2}
&
=
\sum_{l_1,l_2=1}^m\,f_{l_1,l_2}
\left[
g_{i_1}^{l_1}\,g_{i_2}^{l_2}
\right]
+
\sum_{l_1=1}^m\,f_{l_1}
\left[
g_{i_1,i_2}^{l_1}
\right],
\\
h_{i_1,i_2,i_3}
&
=
\sum_{l_1,l_2,l_3=1}^m\,f_{l_1,l_2,l_3}
\left[
g_{i_1}^{l_1}\,g_{i_2}^{l_2}\,g_{i_3}^{l_3}
\right]
+
\sum_{l_1,l_2=1}^m\,f_{l_1,l_2}
\left[
g_{i_1}^{l_1}\,g_{i_2,i_3}^{l_2}
+
g_{i_2}^{l_1}\,g_{i_1,i_3}^{l_2}
+
g_{i_3}^{l_1}\,g_{i_1,i_2}^{l_2}
\right]
+ \\
& \
\ \ \ \ \
+
\sum_{l_1=1}^m\,f_{l_1}
\left[
g_{i_1,i_2,i_3}^{l_1}
\right],
\endaligned
\end{equation}
$$
\small
\aligned
h_{i_1,i_2,i_3,i_4}
&
=
\sum_{l_1,l_2,l_3,l_4=1}^m\,f_{l_1,l_2,l_3,l_4}
\left[
g_{i_1}^{l_1}\,g_{i_2}^{l_2}\,g_{i_3}^{l_3}\,g_{i_4}^{l_4}
\right]
+ \\
& \
\ \ \ \ \ 
\sum_{l_1,l_2,l_3=1}^m\,f_{l_1,l_2,l_3}
\left[
g_{i_2}^{l_1}\,g_{i_3}^{l_2}\,g_{i_1,i_4}^{l_3}
+
g_{i_3}^{l_1}\,g_{i_1}^{l_2}\,g_{i_2,i_4}^{l_3}
+
g_{i_1}^{l_1}\,g_{i_2}^{l_2}\,g_{i_3,i_4}^{l_3}
+
\right.
\\
& \
\ \ \ \ \ \ \ \ \ \ \ \ \ \ \ \ \ \ \ 
\ \ \ \ \ \ \ \ \ \ \ \ \ \ \ \ \ \
\left.
+
g_{i_1}^{l_1}\,g_{i_4}^{l_2}\,g_{i_2,i_3}^{l_3}
+
g_{i_2}^{l_1}\,g_{i_4}^{l_2}\,g_{i_3,i_1}^{l_3}
+
g_{i_3}^{l_1}\,g_{i_4}^{l_2}\,g_{i_1,i_2}^{l_3}
\right]
+ \\
& \
\ \ \ \ \
+
\sum_{l_1,l_2=1}^m\,f_{l_1,l_2}
\left[
g_{i_1,i_2}^{l_1}\,g_{i_3,i_4}^{l_2}
+
g_{i_1,i_3}^{l_1}\,g_{i_2,i_4}^{l_2}
+
g_{i_1,i_4}^{l_1}\,g_{i_2,i_3}^{l_2}
\right]
+ \\
& \
\ \ \ \ \
+
\sum_{l_1,l_2=1}^m\,f_{l_1,l_2}
\left[
g_{i_1}^{l_1}\,g_{i_2,i_3,i_4}^{l_2}
+
g_{i_2}^{l_1}\,g_{i_1,i_3,i_4}^{l_2}
+
g_{i_3}^{l_1}\,g_{i_1,i_2,i_4}^{l_2}
+
g_{i_4}^{l_1}\,g_{i_1,i_2,i_3}^{l_2}
\right]
+ \\
& \
\ \ \ \ \
+
\sum_{l_1=1}^m\,f_{l_1}
\left[
g_{i_1,i_2,i_3,i_4}^{l_1}
\right].
\endaligned
$$
Introducing the derivations
\def\theequation{5.18}\begin{equation}
\small
\aligned
F_i^2
&
:= 
\sum_{k_1=1}^n\,\sum_{l_1=1}^m\,g_{k_1,i}^{l_1}\,
\frac{\partial}{\partial g_{k_1}^{l_1}}
+
\sum_{l_1=1}^m\,g_i^{l_1}\left(
\sum_{l_2=1}^m\,f_{l_1,l_2}\,\frac{\partial}{\partial f_{l_2}}
\right), \\
F_i^3
&
:=
\sum_{k_1=1}^n\,\sum_{l_1=1}^m\,g_{k_1,i}^{l_1}\,
\frac{\partial}{\partial g_{k_1}^{l_1}}
+
\sum_{k_1,k_2=1}^n\,\sum_{l_1=1}^m\,g_{k_1,k_2,i}^{l_1}\,
\frac{\partial}{\partial g_{k_1,k_2}^{l_1}}
+ \\
& \
\ \ \ \ \ 
+
\sum_{l_1=1}^m\,g_i^{l_1}
\left(
\sum_{l_2=1}^m\,f_{l_1,l_2}\,\frac{\partial}{\partial f_{l_2}}
+
\sum_{l_2,l_3=1}^m\,f_{l_1,l_2,l_3}\,
\frac{\partial}{\partial f_{l_2,l_3}}
\right),
\\
\text{\bf 
\dots\dots}
& \
\ \ \ \ \
\text{\bf 
\dots\dots\dots\dots\dots\dots\dots\dots\dots
\dots\dots\dots\dots\dots\dots\dots\dots\dots
\dots
} \\
F_i^\lambda
&
:=
\sum_{k_1=1}^n\,\sum_{l_1=1}^m\,g_{k_1,i}^{l_1}\,
\frac{\partial}{\partial g_{k_1}^{l_1}}
+
\sum_{k_1,k_2=1}^n\,\sum_{l_1=1}^m\,g_{k_1,k_2,i}^{l_1}\,
\frac{\partial}{\partial g_{k_1,k_2}^{l_1}}
+ 
\cdots
+ \\
& \
\ \ \ \ \ 
+
\sum_{k_1,k_2,\dots,k_{\lambda-1}=1}^n\,\sum_{l_1=1}^m\,
g_{k_1,k_2,\dots,k_{\lambda-1},i}
\
\frac{\partial}{\partial g_{k_1,\dots,k_{\lambda-1}}^{l_1}}
+ \\
& \
\ \ \ \ \ 
+
\sum_{l_1=1}^m\,g_i^{l_1}
\left(
\sum_{l_2=1}^m\,f_{l_1,l_2}\,\frac{\partial}{\partial f_{l_2}}
+
\sum_{l_2,l_3=1}^m\,f_{l_1,l_2,l_3}\,
\frac{\partial}{\partial f_{l_2,l_3}}
+ 
\right. 
\\
& \
\ \ \ \ \ \ \ \ \ \ \ \ \ \ \ \ \ \ \ \ 
\ \ \ \ \ \ \ \ \ \
\left.
+
\cdots
+
\sum_{l_2,l_3,\dots,l_\lambda}\,
f_{l_1,l_2,l_3,\dots,l_\lambda}\,
\frac{\partial}{\partial f_{l_2,l_3,\dots,l_\lambda}}
\right),
\endaligned
\end{equation}
we observe that the following
induction relations hold:
\def\theequation{5.19}\begin{equation}
\aligned
h_{i_1,i_2}
&
=
F_{i_2}^2
\left(
h_{i_1}
\right), 
\\
h_{i_1,i_2,i_3}
&
=
F_{i_3}^3
\left(
h_{i_1,i_2}
\right), 
\\
\text{\bf 
\dots\dots
}
& \
\ \ \ \ \
\text{\bf 
\dots\dots\dots\dots\dots
}
\\
h_{i_1,i_2,\dots,i_\lambda}
&
=
F_{i_\lambda}^\lambda
\left(
h_{i_1,i_2,\dots,i_{\lambda-1}}
\right).
\endaligned
\end{equation}

To obtain the explicit version of the Fa\`a di Bruno in the case of
several variables $(x^1, \dots, x^n)$ and several variables $(y^1,
\dots, y^m)$, it suffices to extract from the expression of ${\bf Y
}_{ i_1,\dots, i_\kappa }^j$ provided by Theorem~5.12 only the terms
corresponding to $\mu_1 \lambda_1 + \cdots + \mu_d\lambda_d = \kappa$,
dropping all the $\mathcal{ X}$ terms. After some simplifications and
after a translation by means of an elementary dictionary, we obtain
the fourth and the most general multivariate Fa\`a di Bruno formula.

\def\thetheorem{5.20}\begin{theorem} 
For every integer $\kappa \geq 1$ and for every choice of indices
$i_1, \dots, i_\kappa$ in the set $\{ 1, 2, \dots, n\}$, the
$\kappa$-th partial derivative of the composite function 
\def\theequation{5.21}\begin{equation}
h = h( x^1,
\dots, x^n) = 
f \left( g^1 (x^1, \dots, x^n),\dots, g^m (x^1, \dots,
x^n) \right) 
\end{equation}
with respect to the variables $x^{i_1}, \dots,
x^{i_\kappa}$ may be expressed as an explicit polynomial depending on
the partial derivatives of $f$, on the partial derivatives of the
$g^j$ and having integer coefficients{\rm :}
\def\theequation{5.22}\begin{equation}
\boxed{
\aligned
\frac{\partial^\kappa h}{\partial x^{i_1}\cdots
\partial x^{i_\kappa}}
&
=
\sum_{d=1}^\kappa
\
\sum_{1\leq \lambda_1 < \cdots < \lambda_d \leq \kappa}
\
\sum_{\mu_1\geq 1,\dots,\mu_d\geq 1}
\
\sum_{\mu_1\lambda_1+\cdots+\mu_d\lambda_d=\kappa}
\\
& \
\ \ \ \ \
\sum_{l_{1:1},\dots,l_{1:\mu_1}=1}^m
\cdots
\sum_{l_{d:1},\dots,l_{d:\mu_d}=1}^m
\
\frac{\partial^{\mu_1+\cdots+\mu_d}f}{
\partial y^{l_{1:1}}\cdots
\partial y^{l_{1:\mu_1}}\cdots 
\partial y^{l_{d:1}}\cdots
\partial y^{l_{d:\mu_d}}
}
\\
& 
\ \ \ \ \ \ \ \ \ \ 
\left[
\aligned
&
\sum_{\sigma\in\mathfrak{F}_\kappa^{
(\mu_1,\lambda_1),\dots,(\mu_d,\lambda_d)}}
\
\prod_{1\leq\nu_1\leq\mu_1}
\
\frac{\partial^{\lambda_1} g^{l_{1:\nu_1}}}{\partial
x^{i_{\sigma(1:\nu_1:1)}}\cdots
\partial x^{i_{\sigma(1:\nu_1:\lambda_1)}}}
\
\text{\bf \dots} 
\\
& \
\ \ \ \ \ \ \ \ \ \ \ \ \ \ 
\ \ \ \ \ \ \ \ \ \ \ \ \ \
\text{\bf \dots} 
\prod_{1\leq\nu_d\leq\mu_d}
\
\frac{\partial^{\lambda_d} g^{l_{d:\nu_d}}}{\partial
x^{i_{\sigma(d:\nu_d:1)}}\cdots
\partial x^{i_{\sigma(d:\nu_d:\lambda_d)}}}
\endaligned
\right].
\endaligned
}
\end{equation}
\end{theorem}

\vfill\end{document}